\documentclass[11pt,reqno]{amsart}
\usepackage{amsmath,amsfonts,amsthm,color,amssymb, mathrsfs,natbib, todonotes}
\usepackage{mathtools}
\usepackage{xcolor}
\usepackage{dsfont}
\usepackage{stmaryrd} 
\SetSymbolFont{stmry}{bold}{U}{stmry}{m}{n}
\usepackage[all]{xy}
\usepackage[left=1in, right=1in, top=1.1in,bottom=1.1in]{geometry}
\usepackage{enumitem}

\usepackage{hyperref}
\hypersetup{
colorlinks   = true,
citecolor    = blue,
linkcolor=blue
}
\allowdisplaybreaks
\numberwithin{equation}{section}

\newtheorem{theorem}{Theorem}[section]
\newtheorem{lemma}[theorem]{Lemma}
\newtheorem{corollary}[theorem]{Corollary}
\newtheorem{definition}[theorem]{Definition}
\newtheorem{example}{Example}
\newtheorem{proposition}[theorem]{Proposition}

\newtheorem{remark}[theorem]{Remark}

\newtheorem{assumption}{Assumption}

\def\osigma{\overline{\sigma}}
\def\usigma{\underline{\sigma}}

\def\R{\mathbb R}

\def\rme{{\rme}}

\def\tr{\text{tr}}

\def\la{\left\langle}
\def\ra{\right\rangle}

\def\E{\mathbb E}

\def\dd{{\mathrm d}}

\def\Lip{\mathrm{Lip}}

\def\mathbi#1{\textbf{\em #1}} 

\begin{document}
\title{Generalized Wishart processes  under nonlinear expectation}

\author[X. Li]{Xinpeng Li}\address{Research Center for Mathematics and Interdisciplinary Sciences; Frontiers Science Center for Nonlinear Expectations, Ministry of Education, Shandong University, Qingdao 266237, China}
\email{lixinpeng@sdu.edu.cn}

\author[Y. Mi]{Yiran Mi}\address{Research Center for Mathematics and Interdisciplinary Sciences; Frontiers Science Center for Nonlinear Expectations, Ministry of Education, Shandong University, Qingdao 266237, China}
\email{202317071@mail.sdu.edu.cn}

\author[J. Song]{Jian Song}\address{Research Center for Mathematics and Interdisciplinary Sciences; Frontiers Science Center for Nonlinear Expectations, Ministry of Education, Shandong University, Qingdao 266237, China}
\email{txjsong@sdu.edu.cn}

\author[S. Yang]{Shuzhen Yang}\address{Zhongtai Securities Institute for Financial Studies, Shandong University, Jinan, Shandong, 250100, China.}
\email{yangsz@sdu.edu.cn}

\begin{abstract} 
We develop a theory of generalized Wishart processes within the framework of $G$-expectation, with $G$-Dyson Brownian motion and the $G$-Wishart process as two typical examples. In contrast to the classical setting, we show that $G$-Dyson Brownian motion may exhibit collisions and provide a sufficient condition ensuring the non-collision property of the eigenvalues of generalized Wishart processes. We then derive the equations governing the high-dimensional limits of the empirical spectral measures under $G$-expectation, leading to the corresponding semicircle law for $G$-Dyson Brownian motion and Marchenko--Pastur law for the $G$-Wishart process.


\end{abstract}

\maketitle
\tableofcontents

\section{Introduction}

Random matrix theory studies the spectral properties of matrices with random entries and was initiated by the pioneering works of Wishart \cite{Wishart1928} and Wigner \cite{Wigner1955,Wigner1958}. 
One of its main objectives is to understand the behavior of empirical spectral measures. Two foundational results in this direction are Wigner's semicircle law for symmetric or Hermitian matrices and the Marchenko--Pastur law for sample covariance matrices, established by Wigner \cite{Wigner1958} and Marchenko--Pastur \cite{MarchenkoPastur1967}, respectively.  Random matrix methods are now widely used in statistics, finance, communications, and machine learning, see, e.g., \cite{BaiSilverstein2010,Johnstone2006,LalouxEtAl1999,PlerouEtAl1999,TulinoVerdu2004,CouilletDebbah2011}.


A fundamental dynamical model in random matrix theory is Dyson Brownian motion~\cite{Dyson1962}, which describes the eigenvalue processes of symmetric or Hermitian matrix-valued Brownian motions. The eigenvalues evolve as interacting particles with pairwise repulsion generated by a logarithmic potential, thereby connecting matrix diffusions, Coulomb gases, and universality theory~\cite{ErdosYau2017}. Another important class is given by Wishart processes \cite{Bru1989,Bru1991}, which model evolving covariance matrices and give rise to nonnegative eigenvalue dynamics.

Graczyk and Małecki \cite{GraczykMalecki2013,GraczykMalecki2014} introduced the so-called generalized Wishart processes as solutions to a class of matrix-valued SDEs that encompass, in particular, Dyson Brownian motion and Wishart process. Subsequently, the high-dimensional spectral limits of these processes were investigated in \cite{syy20}, and the fluctuations around the limiting spectral distributions were studied in~\cite{syy21}. We also refer to the survey \cite{syy22} for an overview of matrix-valued stochastic processes.

The above-mentioned works
are formulated on a probability space equipped with a single probability measure. In many applications, however, the probabilistic model itself is uncertain. For example, the volatility or covariance structure of a financial system may not be known precisely and may vary over a prescribed family of plausible models. This motivates the use of sublinear expectations in order to capture model ambiguity \cite{ArtznerEtAl1999,ChenEpstein2002}.
Denis and Martini~\cite{Denis2006} established one of the earliest general duality results for model uncertainty in continuous-time finance. Their work laid important foundations for the subsequent development of robust finance and is closely related in spirit to the theory of $G$-expectation. Peng's $G$-expectation framework \cite{Peng2007,Peng2008} provides a canonical setting for modeling volatility uncertainty and gives rise to $G$-Brownian motion and a nonlinear stochastic calculus \cite{DenisHuPeng2011,SonerTouziZhang2011}. We refer to the monograph \cite{peng2019nonlinear} for a systematic treatment of the theory of $G$-expectation and its applications.

Compared with the extensive literature on  (backward) $G$-stochastic differential equations 
(see, e.g., \cite{Gao2009, gqy14, hjps14, li2016, gpp18, peng2019nonlinear}), random matrix theory under nonlinear expectations remains at an early stage. Initial work by Stihi, Boutabia, and Meradji \cite{StihiBoutabiaMeradji2018} and Meradji, Boutabia, and Stihi~\cite{MeradjiBoutabiaStihi2019} introduced matrix-valued $G$-Brownian motion and $G$-Wishart processes, deriving eigenvalue dynamics and noncollision conditions.  However, large-dimensional spectral limits and nonlinear analogues of classical laws remain largely open.

In this paper, we investigate the spectral dynamics of generalized Wishart processes within the framework of nonlinear expectation and derive the McKean--Vlasov equations governing the high-dimensional limits of their empirical spectral measures. We focus in particular on two important special cases, the symmetric $G$-Brownian matrix and the $G$-Wishart process, and establish the corresponding semicircle law and Marchenko--Pastur law under $G$-expectation.

\subsection{Main results on dynamics of spectral processes} We first present our main results on the equation satisfied by the eigenvalue processes of generalized Wishart processes in the framework of $G$-expectation.

Let $\mathbb S_d$ denote the space of $d\times d$ symmetric matrices. A $d$-dimensional $G$-Brownian motion $B=(B^1, \dots, B^d)$ is characterized by its generator $G: \mathbb{S}_d \to \mathbb{R}$,  which is a  monotone  sublinear function (see Section~\ref{sec:pre} for details).  In this work, we impose the following condition on $G$:
\begin{assumption}
\label{assum:G}
Assume 
$$G(A)=\mathfrak g(\operatorname{tr}(A))=\mathfrak g\Big(\sum_{i=1}^{d} A_{i i}\Big),$$
for
$A=(A_{ij})_{1\leq i,j\leq d}\in\mathbb{S}_d$,
where
\begin{equation}\label{e:G-tilde}
 \mathfrak g(x)
=
\frac{1}{2}
\left(
\overline{\sigma}^{\,2}x^+
-
\underline{\sigma}^{\,2}x^-
\right), \quad x\in\mathbb{R},
\end{equation}
with
$0<\underline{\sigma}\leq\overline{\sigma}<\infty$.
\end{assumption}

When $\usigma=\osigma=1$ in Assumption~\ref{assum:G}, the $G$-Brownian motion $B$ reduces to a classical $d$-dimensional Brownian motion with independent components. Moreover, by~\cite[Lemma~3.2]{hu2026g}, under Assumption~\ref{assum:G}, the quadratic variations and quadratic covariations of the components of $B$  satisfy, for $1\le i\neq j \le d$,
\begin{equation}\label{e:diagonal-prop}
\left\langle B^{1}\right\rangle_{t}=\left\langle B^{i}\right\rangle_{t} \text{ and } \left\langle B^{i}, B^{j}\right\rangle_{t}=0, \quad \text{quasi-surely for all } t>0.
\end{equation}
Conversely, if a $d$-dimensional $G$-Brownian motion $B$ satisfies~\eqref{e:diagonal-prop}, then $G$ satisfies Assumption~\ref{assum:G}. 
In Section~\ref{sec:pre-G-Dyson}, we provide a preliminary study of $G$-Dyson Brownian motion, which suggests that Assumption~\ref{assum:G} is a natural condition for studying the generalized $G$-Wishart process~\eqref{e:g-Wishart} below (see also Section~\ref{rem:choice1}).

To describe the underlying matrix dynamics, we consider the following \emph{generalized $G$-Wishart process}.   Let \(X=(X_t)_{t\geq0}\) be an \(\mathbb{S}_N\)-valued stochastic process
satisfying 
\begin{equation}\label{e:g-Wishart}   \mathrm{d}X^{N}_{t}=g_{N}\left(X^{N}_{t}\right) \mathrm{d} B_{t} h_{N}\left(X^{N}_{t}\right)+h_{N}\left(X^{N}_{t}\right) \mathrm{d} B_{t}^{\top} g_{N}\left(X^{N}_{t}\right)+b_{N}\left(X^{N}_{t}\right) \dd\langle B^{11} \rangle_t,
    \quad t\geq0,
\end{equation}
where \(g_N,h_N,b_N: \R\to\mathbb{R}\) are continuous functions acting on the spectrum of $X_t^N$  and
$ B_t=\bigl(B_t^{kh}\bigr)_{1\leq k,h\leq N}$
is a \(d\)-dimensional \(G\)-Brownian motion satisfying Assumption~\ref{assum:G} with $d=N^2$.

The generalized $G$-Wishart process \eqref{e:g-Wishart} includes two important special cases: the symmetric $G$-Brownian matrix, whose eigenvalue process is  $G$-Dyson Brownian motion, and the $G$-Wishart process, which has potential applications in statistics and finance. Indeed, if we take $g_N(x) = (2N)^{-1/2}$, $h_N(x) = 1$ and $b_N(x) =
  0$ in \eqref{e:g-Wishart}, then we get
\begin{equation}\label{e:X-Dyson}
\mathrm{d}X^{N}_{t}=\frac{1}{\sqrt{2N}} \left(\mathrm{d} B_{t} +\dd B_{t}^{\top}\right).
\end{equation}
Denote \begin{equation*}\mathrm{d}\tilde{B}_{kh}=
\begin{cases}
\dfrac{1}{\sqrt 2} \left(\dd B_{kh}+\dd B_{hk}\right), & \text{if } k < h; \\[5pt]
 \mathrm{d}B_{hh}, & \text{if } k = h.
\end{cases}
\end{equation*}
Under Assumption \ref{assum:G} with $d=N\times N$, we have that $\langle\tilde{B}_{kh} \rangle _{t}=\langle\frac{1}{\sqrt 2} \left(B_{kh}+B_{hk}\right)  \rangle _{t}=\left \langle B_{11}  \right \rangle _{t}$ for $k<h$. Then $\langle \tilde{B}_{kh},\tilde{B}_{k^{'}h^{'}}  \rangle _{t}=\mathrm{I}_{\{(k,h)=(k',h')\}}  \langle\tilde{B}_{11}  \rangle _{t}$, and by Corollary \ref{corollary:G-brownian motion} it is easy to verify that $\tilde{B}$ is   a $\frac{N\left ( N+ 1\right ) }{2}$-dimensional  $\bar{G}$-Brownian motion with $\bar{G}$ satisfying Assumption~\ref{assum:G} with $d=\frac{N\left ( N+ 1\right ) }{2}$.
Hence, $X^N$ given by~\eqref{e:X-Dyson} can be represented as the following $N\times N$ symmetric $G$-Brownian matrix whose upper-triangular entries satisfy
\begin{equation}\label{eq:BN}
    \dd X^{N} =
\begin{cases}
\dfrac{1}{\sqrt{N}} \mathrm{d}\tilde B^{kh}, & \text{if } k < h, \\[5pt]
\dfrac{\sqrt{2}}{\sqrt{N}} \mathrm{d}\tilde B^{hh}, & \text{if } k = h.
\end{cases}
\end{equation}
If we instead choose $g_N(x) = \sqrt{x}$, $h_N(x) = 1$, and $b_N(x) =p$, then $X_t^N$ is the $G$-Wishart process $B_t^\intercal B_t$ (see Lemma~\ref{lem:G-Wishart}), where $B_t$ is a $p\times N$ $G$-Brownian matrix.

The well-posedness of the generalized Wishart process \eqref{e:g-Wishart} is assumed throughout this paper and is not investigated in full generality.  Existence and uniqueness are immediate in the symmetric $G$-Brownian matrix case and are established for the $G$-Wishart process in Theorem~\ref{thm:wishart original existence}.

The following result characterizes the dynamics of the eigenvalue vector $\lambda=(\lambda_1, \dots, \lambda_N)$ of \(X^{N}\) under Assumption~\ref{assum:G} and establishes the quasi-sure non-collision property of  the eigenvalue processes. \nameref{proof:thm:eigen-dyn} is provided in Section~\ref{subsection:Generalized Wishart}.

\begin{theorem}\label{thm:eigenvalue-dynamics}
Let $X^{N}$ be given by~\eqref{e:g-Wishart},  assuming that the initial eigenvalues are distinct and that $X^N_t$ has a quasi-continuous version for each $t\in[0,T]$. Assume that $b_{N},g_{N}^{2},h_{N}^{2}$ are Lipschitz continuous functions, $g_{N}^{2}h_{N}^{2}$ is either convex or belongs to the class $C^{1,1}(\mathbb R)$, where $C^{1,1}(\mathbb{R})$  denotes the space of differentiable functions whose derivatives are Lipschitz continuous. We further assume that the underlying $G$-Brownian motion satisfies Assumption~\ref{assum:G} with $d=N\times N$.

 Let
$\lambda^N_1(t),\ldots,\lambda^N_N(t)$
denote the eigenvalues of \(X^N_t\), and define the first collision time by
\begin{equation}\label{equation:collision time}
    \tau:=\inf \left\{t\ge 0: \lambda^{N}_{i}(t)=\lambda^{N}_{j}(t) \text { for some } i \neq j\right\}.\end{equation}
Then, $\tau=\infty$ quasi-surely. Moreover, the eigenvalue processes, arranged in increasing order as $\lambda^N_1(t)< \lambda^N_2(t)<\cdots<\lambda^N_N(t)$, satisfy the following $G$-SDE: for $i=1,2,\dots,N$,
\begin{align}\label{eq:eigenvalue-dynamics}
    \mathrm{d} \lambda^{N}_{i}=2 g_{N}(\lambda^{N}_{i}) h_{N}(\lambda^{N}_{i}) \mathrm{d} \bar B^i_t+b_{N}(\lambda^{N}_{i})  \mathrm{d} \langle \bar B^i \rangle _{t}+\sum_{j:j \neq i} \frac{H_{N}(\lambda^{N}_{i}, \lambda^{N}_{j})}{\lambda^{N}_{i}-\lambda^{N}_{j}} \mathrm{d} \langle \bar B^i  \rangle _{t},  
\end{align}
where $\bar B=(\bar B^{1}, \dots, \bar B^{N})$ is an $N$-dimensional $\bar G$-Brownian motion with $\bar G$ satisfying Assumption~\ref{assum:G} with $d=N$, and \begin{equation}\label{e:HN}
    H_N(x,y)
    :=
    g_N^2(x)h_N^2(y)+g_N^2(y)h_N^2(x).
\end{equation}
\end{theorem}

The well-posedness of the $G$-SDE~\eqref{eq:eigenvalue-dynamics} is established below, of which \nameref{proof:thm:strong-sol} is given in Section~\ref{subsection:Generalized Wishart}. 

\begin{theorem}\label{thm:strong-sol}
Assume the same conditions as in Theorem~\ref{thm:eigenvalue-dynamics}. Then, there exists a unique strong solution to the $G$-SDE~\eqref{eq:eigenvalue-dynamics}.
\end{theorem}


As a consequence of  Theorem \ref{thm:eigenvalue-dynamics}, we have the following result for $G$-Dyson Brownian motion:
 
 \begin{theorem}\label{thm:dyson non-collision}
Let $X^{N}$ be given by~\eqref{eq:BN} and $\tau$ be given by~\eqref{equation:collision time}. Suppose that the
initial eigenvalues are distinct.    
We further assume that the
underlying G-Brownian motion satisfies Assumption~\ref{assum:G} with $d=\frac{N\left ( N+ 1\right ) }{2}$.
Then, $\tau=\infty$  quasi-surely. Moreover, the eigenvalue processes, arranged in increasing order as $\lambda _{1}^{N} \left ( t \right )<\lambda _{2}^{N} \left ( t \right )<\dots<\lambda _{N}^{N} \left ( t \right )$, satisfy the following G-SDE: for $i=1,\dots,N$,
\begin{equation}\label{e:G-Dyson}
\begin{aligned}
\mathrm{d}\lambda _{i}^{N}=\frac{\sqrt{2} }{\sqrt{N} } \mathrm{d}\bar{B}^{i}   +\frac{1}{N}  \sum_{j: j \neq i} \frac{1}{\lambda_{i}^{N}-\lambda_{j}^{N}} \mathrm{d} \langle \bar{B}^{i}  \rangle _{t},
\end{aligned}
\end{equation}
where $\bar{B}=(\bar{B}^1, \dots, \bar{B}^N)$ is an $N$-dimensional  $\bar{G}$-Brownian motion with $\bar{G}$ satisfying Assumption~\ref{assum:G} with $d=N$.
\end{theorem}


Theorem~\ref{thm:strong-sol} establishes the existence and uniqueness of a strong solution to the $G$-SDE~\eqref{e:G-Dyson} when the initial eigenvalues are distinct. In fact, this result remains valid even when the initial eigenvalues are allowed to coincide, as shown in the following proposition, of which \nameref{proof:prop:closed-chamber} is given in Section~\ref{subsection:Discussion of \texorpdfstring{$G$}{}-Dyson Brownian motion}. Set $\Delta_{N}:=\big\{x \in \mathbb{R}^{N}: x_{1}<\cdots<x_{N}\big\}$ and denote by $\overline\Delta_N$ its closure. 
 
\begin{proposition}
\label{prop:closed-chamber} Assume Assumption~\ref{assum:G} for the underlying $G$-Brownian motion. Then, there exists a unique solution
$\lambda^N\in C\bigl(\mathbb{R}_+;\overline{\Delta}_N\bigr)$ to the G-SDE~\eqref{e:G-Dyson} with initial
condition $\lambda^N(0)\in\overline{\Delta}_N$. Moreover, $\lambda^N(t)\in\Delta_N$ for all 
 $t>0$ quasi-surely.  
\end{proposition}


An $N\times N$ \emph{Wishart process} with $p$ degrees of freedom is  defined by
\begin{equation}\label{wishart definition}
    X^{N}_{t}:=Y_{t}^{\top} Y_{t} \in \mathbb{S}_{N}^+,
    \end{equation} 
where $Y_{t}$  is a  $p \times N$  $\tilde G$-Brownian matrix and $\mathbb{S}_{N}^+$ denotes the space of $N \times N$ nonnegative definite matrices. Assuming that $\tilde G$ satisfies Assumption~\ref{assum:G} with $d=p\times N$  and $p\ge N$, by Lemma~\ref{lem:G-Wishart}, $X^{N}_{t}$ satisfies the following  $G$-SDE:
 \begin{equation}\label{G-N-wishart SDE}
    \mathrm{d}X^{N}_t = \sqrt{X^{N}_t} \mathrm{d}
   B_t + \mathrm{d} B_t^{\top} \sqrt{X^{N}_t} + p \mathbf{I}_N\mathrm{d}\langle  B^{11} \rangle_t,
\end{equation}
where $B$ is an $N\times N$  $G$-Brownian matrix with $G$ satisfying Assumption~\ref{assum:G} with $d=N\times N$. We consider the following \emph{general 
$G$-Wishart process} with  parameter  $\alpha $:
\begin{equation}\label{e:G-Wishart}
\left\{
\begin{aligned}
\mathrm{d}X^{N}_{t} &=\sqrt{X^{N}_{t}} \,\mathrm{d} B_{t}+\mathrm{d} B_{t}^{\top}\sqrt{X^{N}_{t}}+\alpha \mathbf{I}_N\,\mathrm{d}\langle B^{11} \rangle_t \\
X^{N}_{0} &=x_{0}.
\end{aligned}
\right.
\end{equation}

The following result establishes the existence and uniqueness of a positive-definite solution to the $G$-SDE~\eqref{e:G-Wishart}, \nameref{proof:thm:wishart original existence} of which is given in Section~\ref{subsection:Discussion of \texorpdfstring{$G$}{}-Wishart process}.
\begin{theorem}\label{thm:wishart original existence}
Assuming $\alpha \ge N + 1$, the $G$-SDE~\eqref{e:G-Wishart} with  $x_0 \in \mathbb{S}_N^{++}$,   where $\mathbb{S}_{N}^{++}$ denotes the space of $N \times N$ positive definite matrices, admits a unique solution $X^{N}_t \in \mathbb{S}_N^{++}$ q.s.\  for $t\in[0,T]$. Moreover, we have $X_t^N\in (L_G^2(\Omega_t))^d$  with $d=N\times N$ for $t\in[0,T]$.
\end{theorem}

Theorem~\ref{thm:eigenvalue-dynamics} implies directly the following result for $G$-Wishart process.

\begin{theorem}\label{thm:wishart non-collision}
Let $X^{N}$ be given by~\eqref{e:G-Wishart}, where $x_{0} \in \mathbb{S}_N^{++}$ and $\tau$ be given by~\eqref{equation:collision time}. Suppose that the initial eigenvalues are distinct. We further assume that the
underlying G-Brownian motion satisfies Assumption \ref{assum:G} with $d=N\times N$.
Then, $\tau=\infty$  quasi-surely. Moreover, the eigenvalue processes, arranged in increasing order as $\lambda _{1}^{N} \left ( t \right )<\lambda _{2}^{N} \left ( t \right )<\dots<\lambda _{N}^{N} \left ( t \right )$, satisfy the following G-SDE: for $i=1,2,\dots,N$,
\begin{align}\label{eq:wishart eigenvalue}
    \mathrm{d} \lambda^{N}_{i}=2 \sqrt{\lambda^N_{i}}\mathrm{d} \bar B^i+\alpha \mathrm{d}\langle \bar B^i \rangle _{t}+\sum_{j:j \neq i} \frac{\lambda^{N}_{i}+\lambda^{N}_{j}}{\lambda^{N}_{i}-\lambda^{N}_{j}}\mathrm{d}\langle \bar B^i \rangle _{t},
\end{align}
where $\bar B=(\bar B^1, \dots, \bar B^N)$ is an  $N$-dimensional  $\bar{G}$-Brownian motion with $\bar{G}$ satisfying Assumption~\ref{assum:G} with $d=N$.
\end{theorem}

In particular, for an $N\times N$ Wishart process $X^N$ given by ~\eqref{wishart definition} with $p\ge N+1$, we have the following result, of which  \nameref{proof:thm: wishart initial zero} is given in Section \ref{subsection:Discussion of \texorpdfstring{$G$}{}-Wishart process}.
\begin{proposition}
    \label{thm: wishart initial zero}
Let $X^{N}$ be given by~\eqref{G-N-wishart SDE}, where $p \ge N+1$.  
Let
$\lambda_1^N(0)\le \cdots \le \lambda_N^N(0)$
be the eigenvalues of \(x_0\).
We further assume that the
underlying G-Brownian motion satisfies Assumption~\ref{assum:G} with $d=N\times N$. Then, quasi-surely, $0<\lambda_1^N(t)<\cdots<\lambda_N^N(t)$ for every $t>0$.
For $t\ge0$ and $i=1,2,\dots,N$,
\begin{align}\label{eq:wishart eigenvalue zero}
    \mathrm{d} \lambda^{N}_{i}=2 \sqrt{\lambda^N_{i}}\mathrm{d} \bar B^i+p \mathrm{d}\langle \bar B^i \rangle _{t}+\sum_{j:j \neq i} \frac{\lambda^{N}_{i}+\lambda^{N}_{j}}{\lambda^{N}_{i}-\lambda^{N}_{j}}\mathrm{d}\langle \bar B^i \rangle _{t},
\end{align}
where $\bar B=(\bar B^1, \dots, \bar B^N)$ is an  $N$-dimensional  $\bar{G}$-Brownian motion with $\bar{G}$ satisfying Assumption~\ref{assum:G} with $d=N$.
\end{proposition}

\subsection{Main results on high-dimensional spectral limits}
\label{Main results on high-dimensional spectral limits}
Now we present our results on high-dimensional limits of the empirical eigenvalue
distributions associated with generalized $G$-Wishart processes. Let \(C([0,T],\mathbf{P}(\mathbb{R}))\) denote the space of
continuous functions mapping from \([0,T]\) to \(\mathbf{P}(\mathbb{R})\), where
\(\mathbf{P}(\mathbb{R})\) is the space of probability measures on
\(\mathbb{R}\), equipped with the topology of weak convergence.
We study the high-dimensional limit of the empirical spectral measure processes: 
 \begin{equation}\label{e:LN}
    \Big\{L_N(t)
    :=
    \frac{1}{N}
    \sum_{i=1}^N
    \delta_{\lambda_i^N(t)},~ t\in[0,T]\Big\}_{N\in\mathbb N},
\end{equation}
which are viewed as  
$C\bigl([0,T],\mathbf P(\mathbb{R})\bigr)$-valued random elements. We first establish quasi-sure relative compactness of this sequence under the following assumption.

\begin{assumption}\label{assumption:generalized-wishart}
We assume that the following conditions hold.
\begin{enumerate}
    \item There exists a nonnegative function $\varphi(x)\in C^2(\mathbb{R})$ such that 
    $\lim\limits_{|x|\to\infty}\varphi(x)=\infty$, \(\varphi'(x)b_N(x)\) is bounded uniformly with
    respect to $x\in\mathbb{R}$ and $N\geq 1$, and
     $$\sum_{N=1}^{\infty}
    \Big(\frac{\|\varphi' g_N h_N\|_{\infty}^2}{N}\Big)^{l_1}
    <\infty$$
    for some positive integer ${l_1}$. One typical choice is $\varphi(x)=\ln (1+x^2)$.
    \item The function
     $N H_N(x,y)\frac{\varphi'(x)-\varphi'(y)}{x-y}$
    is bounded uniformly with respect to $x,y\in\mathbb{R} $ and
    $N\geq 1$. 
    \item The initial empirical measure satisfies
     $$C_0:=
    \sup_{N\geq 1}\langle\varphi,L_N(0)\rangle
    =
    \sup_{N\geq 1}\frac{1}{N}
    \sum_{i=1}^N\varphi\bigl(\lambda_i^N(0)\bigr)
    <\infty.$$
    \item There exists a sequence
    $\{\widetilde f_k\}_{k\in\mathbb{N}}$ of $ C^2(\mathbb{R})$ functions which is 
    dense in  the space $C_0(\mathbb{R})$ of continuous
    functions vanishing at infinity such that
    $$C_k:=\sup_{N}\left\{
\left\| \widetilde f_k' b_N \right\|_\infty
+ \left\| \widetilde f_k'' g_N^2 h_N^2 \right\|_\infty
+ \frac12 \left\|
N H_N(x,y) \frac{\widetilde f_k'(x)-\widetilde f_k'(y)}{x-y}
\right\|_\infty
\right\} < \infty,$$ and
    $$\psi(k):=
    \sum_{N=1}^{\infty}
    \Big(
        \frac{
            \|\widetilde f_k' g_N h_N\|_{\infty}^2
        }{N}
    \Big)^{l_2}
    <\infty,
    \quad k\in\mathbb{N}$$
    for some  integer ${l_2}\geq2$.
\end{enumerate}
\end{assumption}

To investigate the high-dimensional limits, we need to work with an infinite-dimensional $\mathbi G$-Brownian motion $B=(B^1, B^2, \dots)$ and its associated capacity $c_\infty$, which are constructed in Section~\ref{sec:Kolmog}.  For each $d\in \mathbb N$, the first $d$ coordinates $(B^1,\dots,B^d)$  of $B$ form a $d$-dimensional $G_d$-Brownian motion driving our system,  with $G_d$ satisfying Assumption \ref{assum:G}.   


The following result characterizes the relative compactness of the sequence $\left\{L_{N}(t), t \in[0, T]\right\}_{N \in \mathbb{N}}$, and \nameref{proof:thm: pre-compact} is given in Section ~\ref{subsection:Generalized Wishart Processes}.
\begin{theorem}\label{thm: pre-compact}
Suppose that \eqref{eq:eigenvalue-dynamics} has a  solution that is nonexploding and noncolliding for  $t \in[0, T]$.  Then under Assumption \ref{assumption:generalized-wishart}, the sequence  $\left\{L_{N}(t), t \in[0, T]\right\}_{N \in \mathbb{N}}$  given by~\eqref{e:LN} is relatively compact in  $C\left([0, T], \mathbf{P}(\mathbb{R})\right)$ quasi-surely.
\end{theorem}

The following result characterizes the spectral limit, of which \nameref{proof:thm:limit-equation} is given in Section~\ref{subsection:Generalized Wishart Processes}.

\begin{theorem}\label{thm:limit-equation}
Suppose that Assumption~\ref{assum:G} holds for each $d\
\in\mathbb N$ and assume the same conditions as in Theorem~\ref{thm: pre-compact}. Furthermore, assume that there exist continuous functions  $b(x)$  and  $H(x,y)$, such that  $b_{N}(x)$  converges to  $b(x)$  and  $N H_{N}(x, y)$  converges to  $H(x,y)$  uniformly as $N\to\infty$, and that
\begin{equation*}
    \left\|\frac{b(x)}{1+x^{2}}\right\|_{\infty}<\infty, \quad\left\|\frac{H(x,y)}{(1+|x|)\left(1+y^{2}\right)}\right\|_{\infty}<\infty. \end{equation*}

Suppose that \(L_N(0)\) converges weakly to \(\mu_0\), quasi-surely, as
\(N\to\infty\), and 
that
$\{L_{N}\}$ has a limit measure $\mu \in C([0,T],\mathbf{P}(\mathbb{R}))$.
Then, for every \(z\in\mathbb{C}\setminus\mathbb{R}\),
\begin{equation}\label{eq:limit-resolvent-equation}
\begin{aligned}
\int \frac{\mu_{t}(\mathrm{d}x)}{z-x}&=  \int \frac{\mu_{0}(\mathrm{d}x)}{z-x}+\int_{0}^{t}\left[\int \frac{b(x)}{(z-x)^{2}} \mu_{s}(\mathrm{d}x)\right]\mathrm{d} \langle \bar{B}^1  \rangle _{s}  \\
&\quad +\int_{0}^{t}\left[\iint \frac{H(x,y)}{(z-x)(z-y)^{2}} \mu_{s}(\mathrm{d}x) \mu_{s}(\mathrm{d}y)\right] \mathrm{d}\langle \bar{B}^1 \rangle _{s}.
\end{aligned}
\end{equation}
\end{theorem}

\begin{remark}
Due to the uncertainty, the high-dimensional limit of the empirical distribution of the eigenvalues is a random measure, noting that the quadratic variation process $\langle  \bar B^1 \rangle _t$ appearing in \eqref{eq:limit-resolvent-equation} is no longer deterministic in the framework of G-expectation. Indeed, $\langle  \bar B^1 \rangle _t$ is maximally distributed (see Definition~\ref{def:maximal}). More precisely, for each $\varphi\in C_{l,Lip}(\R)$ we have $\hat{\mathbb E}[\varphi(\langle  \bar B^1 \rangle _t)]=\sup_{ \underline{\sigma} \le b \le \overline{\sigma}} \varphi(b^2 t)$ and the function $v(t,x) = \hat{\mathbb E}[\varphi(x+\langle  \bar B^1 \rangle _t)]$   satisfies the nonlinear equation  $\partial_t v= 2\mathfrak g(\partial_x v)$ with the initial condition $v(0,x) =\varphi(x)$. 
\end{remark}

\begin{remark}
Let $\langle \mu,f\rangle:=\int_{\mathbb R}f(x)\,\mu(\mathrm dx)$. Equivalently, the McKean--Vlasov equation corresponding to~\eqref{eq:limit-resolvent-equation} can be written as
\begin{equation}\label{eq:general-MV}
\begin{aligned}
\langle\mu_t,f\rangle
={}&\langle\mu_0,f\rangle
+\int_0^t\int_{\mathbb R}
b(x)f'(x)\,\mu_s(\mathrm dx)
\,\mathrm d\langle\bar B^1\rangle_s\\
&+\frac12\int_0^t
\iint_{\mathbb R^2}
H(x,y)\frac{f'(x)-f'(y)}{x-y}
\,\mu_s(\mathrm dx)\mu_s(\mathrm dy)
\,\mathrm d\langle\bar B^1\rangle_s,
\end{aligned}
\end{equation}
for sufficiently smooth test functions $f$.  
Here, we use the convention
\[
\left.
\frac{f'(x)-f'(y)}{x-y}
\right|_{x=y}
:=f''(x).
\]

If we assume that $\mu_t$ admits a density
$\rho_t$, namely, $\mu_t(\mathrm{d}x)=\rho_t(x)\,\mathrm{d}x$, then formally we have the following density
formulation:
\begin{equation}\label{eq:limit-density-equation}
\mathrm{d}\rho_t(x)
=
-\partial_x\left[
\rho_t(x)\left(
b(x)
+
\operatorname{p.v.}\int_{\mathbb{R}}
\frac{H(x,y)}{x-y}\rho_t(y)\,\mathrm{d}y
\right)
\right]
\mathrm{d}\langle\bar{B}^1\rangle_t,
\end{equation}
where
\[
\operatorname{p.v.}\int_{\mathbb{R}}
\frac{H(x,y)}{x-y}\rho_t(y)\,\mathrm{d}y
:=
\lim_{\varepsilon\downarrow0}
\int_{\{y\in\mathbb{R}:|x-y|>\varepsilon\}}
\frac{H(x,y)}{x-y}\rho_t(y)\,\mathrm{d}y,
\]
 denotes the principal value. 


In general, uniqueness of solutions to~\eqref{eq:limit-resolvent-equation}, \eqref{eq:general-MV} and \eqref{eq:limit-density-equation} has not been established. Nevertheless, for the special cases of $G$-Dyson Brownian motion and the $G$-Wishart process, we establish uniqueness of solutions to the corresponding limiting equations of~\eqref{eq:limit-resolvent-equation} in Proposition~\ref{dyson converge} and Proposition~\ref{wishart converge}, respectively.
\end{remark}

As a corollary of Theorems~\ref{thm: pre-compact} and~\ref{thm:limit-equation}, we have the following result for $G$-Dyson Brownian motion. Moreover, the uniqueness of the solution to the equation of the limiting measure can be obtained, \nameref{proof:Proposition  dyson converge} of which is given in Section~\ref{subsection:Dyson Brownian motion}.

\begin{proposition}\label{dyson converge}
Suppose that Assumption~\ref{assum:G} holds for each $d\
\in\mathbb N$ and assume that \eqref{e:G-Dyson} has a  solution $ \lambda^{N}(t)$ which is nonexploding and noncolliding for $t\in[0,T]$ with initial condition $ \lambda^{N}(0)$ satisfying
the corresponding condition (3) in Assumption~\ref{assumption:generalized-wishart}.
Assume that   $L_{N}(0) $ converges weakly as  $N$ goes to infinity to $\nu \in \mathbf{P}(\mathbb{R})$.  Then, $\{L_{N}(t),t \in[0, T]\}$ converges quasi-surely in $ C\left([0, T], \mathbf{P}(\mathbb{R})\right)$ as $N\to\infty$. Its limit is the unique measure-valued process $\{\mu_{t},t \in[0, T]\}$  so that $\mu_{0}=\nu$  and the function
\begin{equation}\label{eq:G}V_{t}(z)=\int(z-x)^{-1}  \mu_{t}( \mathrm{d}x)\end{equation}
satisfies the equation
\begin{equation}\label{eq:Gt}
V_{t}(z)=V_{0}(z)-\int_{0}^{t} V_{s}(z) \partial_{z} V_{s}(z) \mathrm{d} \langle \bar B^1 \rangle _{s},\end{equation}
for  $z \in \mathbb{C} \backslash \mathbb{R}$ .
\end{proposition}

\begin{remark} For $G$-Dyson Brownian motion, the McKean--Vlasov equation~\eqref{eq:general-MV} becomes
\begin{equation*}\label{eq:Dyson-MV}
\begin{aligned}
\langle\mu_t,f\rangle
=\langle\mu_0,f\rangle+\frac12\int_0^t
\iint_{\mathbb R^2}
\frac{f'(x)-f'(y)}{x-y}
\,\mu_s(\mathrm dx)\mu_s(\mathrm dy)
\,\mathrm d\langle\bar B^1\rangle_s,
\end{aligned}
\end{equation*}
and the density formulation~\eqref{eq:limit-density-equation} is
\begin{equation*}\label{eq:Dyson-density-equation}
\mathrm{d}\rho_t(x)
=
-\partial_x\left[
\rho_t(x)\,
\operatorname{p.v.}\int_{\mathbb{R}}
\frac{\rho_t(y)}{x-y}\,\mathrm{d}y
\right]
\mathrm{d}\langle\bar B^1\rangle_t.
\end{equation*}
\end{remark}

\begin{remark} \label{remark:5.3}  
For $\mu_t$ satisfying~\eqref{eq:Gt}, we  have $\mu_t = \mu_0 \boxplus S_{\langle \bar B^1\rangle _t}$, where $\boxplus$ denotes free additive convolution and 
\begin{equation}\label{eq:Sa}
    S_a(\mathrm dx)=\frac{1}{2\pi a} \sqrt{4a-x^2}\mathbf 1_{\{|x|\le 2\sqrt a\}}\mathrm d x,\,\text{for } a>0\, ;S_0=\delta_0\end{equation}
is the semicircle law with variance $a$. If $\mu_0$ is the point mass at 0, then $\mu_t$ coincides with the semicircle law with variance $\langle\bar B^1 \rangle _t $, i.e.,  $\mu_t=S_{\langle \bar B^1 \rangle _t}.$
\end{remark}

As a corollary of Theorems~\ref{thm: pre-compact} and~\ref{thm:limit-equation}, we have the following result for $G$-Wishart process, where the uniqueness of the solution to the equation of the limiting measure is also obtained.  \nameref{proof:Proposition  wishart converge} is provided in Section~\ref{subsection:Wishart process}.

Let $X_{t}^{N}:=\frac{1}{N} B_{t}^{\top} B_{t}$, where  $B_{t}$  is a $ p \times N$  $G$-Brownian matrix satisfying Assumption~\ref{assum:G} with $p\ge N+1$.
Then by Proposition \ref{thm: wishart initial zero}, it is easy to know that the eigenvalue processes $\lambda^{N}$ of $X^N$ satisfy the following  system corresponding to ~\eqref{eq:eigenvalue-dynamics} for $g_{N}(x)=\sqrt{x}$, $h_{N}(x)=1 / \sqrt{N} $, and  $b_{N}(x)=\frac{p}{N}$:
\begin{equation}\label{equation of wishart process}
    \mathrm{d} \lambda^{N}_{i}=\frac{2}{\sqrt N}\sqrt{\lambda^{N}_{i}}\mathrm{d} \bar B_i+\frac{p}{N}\mathrm{d}\langle \bar B_i \rangle _{t}+\frac{1}{N}\sum_{j:j \neq i} \frac{\lambda^{N}_{i}+\lambda^{N}_{j}}{\lambda^{N}_{i}-\lambda^{N}_{j}}\mathrm{d}\langle \bar B_i \rangle _{t}.
\end{equation}
We further assume that
\begin{equation}\label{e:c}
\lim_{N\to\infty}\frac{p}{N}=c \geq 1.
\end{equation}
Then $b_{N}(x)$ converges to  $b(x)=c$ and $N H_{N}(x, y)$  converges to  $H(x,y)=x+y$.


\begin{proposition}\label{wishart converge}
Suppose that Assumption~\ref{assum:G} holds for each $d\
\in\mathbb N$ and assume that~\eqref{equation of wishart process} has a  solution $ \lambda^{N}(t)$ which is nonexploding and noncolliding for  $t \in[0, T]$ with initial condition $ \lambda^{N}(0)$ satisfying
the corresponding condition (3) in Assumption~\ref{assumption:generalized-wishart}.
Assume that the empirical measure  $L_{N}(0) $ converges weakly as  $N\to\infty$ to $\nu \in \mathbf{P}(\mathbb{R}) $.
Then,  $\{L_{N}(t),t \in[0, T]\}$ converges quasi-surely in $C\left([0, T], \mathbf{P}(\mathbb{R})\right)$. Its limit is the unique measure-valued process $\{\mu_{t},t \in[0, T]\}$  so that $ \mu_{0}=\nu$  and the function 
\begin{equation}V_{t}(z)=\int(z-x)^{-1}  \mu_{t}( \mathrm{d}x)\end{equation}
satisfies the equation
\begin{equation}\label{wishart:z-x,G}V_{t}(z)=V_{0}(z)-\int_{0}^{t}\Big[V_{s}(z)^{2}+(c-1)\partial_{z} V_{s}(z)+2 z V_{s}(z) \partial_{z} V_{s}(z)\Big] \mathrm{d}\langle \bar B^1 \rangle _{s},\end{equation}
for  $z \in \mathbb{C} \backslash \mathbb{R}$, where $c\ge 1$ is given in~\eqref{e:c}.
\end{proposition}

\begin{remark}
For the $G$-Wishart process, the corresponding McKean--Vlasov equation of~\eqref{wishart:z-x,G} is
\begin{equation}\label{eq:Wishart-MV}
\begin{aligned}
\langle\mu_t,f\rangle
={}&\langle\mu_0,f\rangle
+c\int_0^t\int_{\mathbb R}
f'(x)\,\mu_s(\mathrm dx)
\,\mathrm d\langle\bar B^1\rangle_s\\
&+\frac12\int_0^t
\iint_{\mathbb R^2}
(x+y)\frac{f'(x)-f'(y)}{x-y}
\,\mu_s(\mathrm dx)\mu_s(\mathrm dy)
\,\mathrm d\langle\bar B^1\rangle_s.
\end{aligned}
\end{equation}
Since the measures are supported
on $[0,\infty)$,  the spatial integrals may be
restricted accordingly. Assuming that $\mu_t$ has a density $\rho_t(x)$, then the density formulation is
\begin{equation*}
\mathrm{d}\rho_t(x)
=
-\partial_x\left[
\rho_t(x)\left(
c-1
+
2x\,\operatorname{p.v.}\int_{\mathbb{R}}
\frac{\rho_t(y)}{x-y}\,\mathrm{d}y
\right)
\right]
\mathrm{d}\langle\bar B^1\rangle_t.
\end{equation*}
\end{remark}

\begin{remark}\label{remark: wishart semicircle law}
If $\mu_0$ is the point mass at 0, then we have
\begin{align*}V_{t}(z)
&=\frac{z + (1-c)\langle \bar B^1 \rangle_t - \sqrt{(z - (c+1)\langle \bar B^1 \rangle_t)^2 - 4c\langle  \bar B^1 \rangle_t^2}}{2 z \langle \bar B^1 \rangle_t},
\end{align*}
 and consequently, $\mu_t$ has the density
\begin{align*}
\rho_t(x)
 = \frac{\sqrt{\left( \langle \bar B^1 \rangle_t(1+\sqrt{ c })^2 - x \right) \left( x - \langle \bar B^1 \rangle_t(\sqrt{ c}-1)^2 \right)}}{2\pi   \langle \bar B^1 \rangle_t x} \mathbf{1}_{\left[ \langle \bar B^1 \rangle_t(\sqrt{ c}-1)
 ^2, \, \langle \bar B^1\rangle_t(1+\sqrt{ c})^2 \right]}(x),\end{align*}
 which is the Marchenko--Pastur law with mean $c \langle \bar B^1\rangle_t$. 
\end{remark}


\subsection{\texorpdfstring{On discrete $G$-Gaussian Wigner and Wishart matrices}{}}

Let $\xi^{(N)}=\bigl(\xi_{ij}^{(N)}\bigr)_{1\leq i\leq j\leq N}$
be a $G$-normal vector of dimension
$d=\frac{N(N+1)}{2}$.
An $N\times N$ $G$-Gaussian Wigner matrix $W_N$ is defined by
$$W_N(i,j)
=
\begin{cases}
\frac{1}{\sqrt N}\xi_{ij}^{(N)}, & i<j,\\[4pt]
\frac{\sqrt2}{\sqrt N}\,\xi_{ii}^{(N)}, & i=j,\\[4pt]
W_N(j,i), & i>j.
\end{cases}$$
Let $L_N^W = \frac{1}{N}\sum_{i=1}^N \delta_{\lambda_i(W_N)}$
be the empirical spectral measure of $W_N$.
For $a>0$, denote  $S_a$ the semicircle law with variance $a$ given by ~\eqref{eq:Sa}.



As a corollary of $G$-Dyson Brownian motion  for $t=1$, we have the following result for the discrete $G$-Gaussian Wigner matrix.

\begin{proposition}
Suppose that Assumption~\ref{assum:G} holds for each $d\
\in\mathbb N$ and assume that $G$ satisfies Assumption \ref{assum:G} with
$d = \frac{N(N+1)}{2}$. 
Then, $L_N^W$ converges weakly to $S_{\langle \bar B^1\rangle_1}$ quasi-surely. 
\end{proposition}



A discrete $N\times N$ $G$-Gaussian Wishart matrix $Q_N$ is defined by
$$Q_N=\frac{1}{N}Y_N^{\mathsf T}Y_N,$$ where
$Y_N=\left(\eta_{ij}^{(N)}\right)_{\substack{1\leq i\leq p,\\1\leq j\leq N}}$
is a $p\times N$ random matrix
whose vector of entries $\eta^{(N)}=\bigl(\eta_{ij}^{(N)}\bigr)_{\substack{1\leq i\leq p,\\1\leq j\leq N}}$ is a  $G$-normal vector of dimension $d=p\times N$.
Let $L_N^Q
=
\frac{1}{N}\sum_{i=1}^{N}
\delta_{\lambda_i(Q_N)}$
denote the empirical spectral distribution of $Q_N$.
For $a>0$, define
$\ell_{c,a}=a(1-\sqrt c)^2,
\,r_{c,a}=a(1+\sqrt c)^2$.
Denote $\Pi_{c,a}$ the Marchenko--Pastur law with parameters $c$ and $a$  given by
$$\Pi_{c,a}(\mathrm{d}x)=(1-c)^+\delta_0(\mathrm{d}x)+\frac{\sqrt{(r_{c,a}-x)(x-\ell_{c,a})}}{2\pi a  x}
\mathbf{1}_{(\ell_{c,a},r_{c,a})}(x)\,\mathrm{d}x.$$

As a corollary of $G$-Wishart process for $t=1$, we have the following result for the discrete $G$-Gaussian Wishart matrix.
\begin{proposition}
Suppose that Assumption~\ref{assum:G} holds for each $d\
\in\mathbb N$  and assume that $G$ satisfies Assumption \ref{assum:G} with
$d = p\times N$, $p\ge N+1$ and $\frac{p}{N}\to c$. Then, $L_N^Q$ converges weakly  to $\Pi_{c,\langle \bar B^1\rangle_1}$ quasi-surely.
\end{proposition}




\subsection{Summary of contributions}
In this work, we investigate generalized Wishart processes  defined by~\eqref{e:g-Wishart} within the framework of $G$-expectation, focusing on eigenvalue non-collision property and high-dimensional spectral limits. Volatility uncertainty affects both the diffusion  and the drift of the eigenvalue dynamics. Consequently, classical non-collision results need not carry over to this setting, and spectral collisions may occur under certain covariance structures. Under these covariance structures, the GOE property may also fail; see Section~\ref{sec:pre-G-Dyson} for details.

We identify a sufficient condition for the spectral non-collision property, stated in Assumption~\ref{assum:G}, and study the empirical spectral measures of generalized $G$-Wishart processes. We characterize their high-dimensional limits by an equation for the Stieltjes transform or, equivalently, a McKean--Vlasov equation for the limiting measure. In particular, we establish the semicircle law for $G$-Dyson Brownian motion and the Marchenko--Pastur law for the $G$-Wishart process. 


In general, the limiting spectral distribution \(\mu_t\) of a generalized \(G\)-Wishart process is random rather than deterministic. This randomness reflects the volatility uncertainty of the underlying \(G\)-expectation, encoded by the maximally distributed random variable \(\langle \bar B^1\rangle_t\). These results show how volatility uncertainty affects spectral dynamics and extend classical spectral limit laws to the framework of nonlinear expectation.

\

The rest of the paper is organized as follows. In Section~\ref{sec:pre}, we review the basic theory of nonlinear expectation and derive several corollaries related to the $G$-Lévy characterization theorem. Section~\ref{sec:pre-G-Dyson} examines alternative assumptions under which the GOE property and the non-collision property of $G$-Dyson Brownian motion may fail. In Section~\ref{sec:non-collision}, we derive the spectral dynamics of generalized Wishart processes and establish the non-collision result under Assumption~\ref{assum:G}, as stated in Theorem~\ref{thm:eigenvalue-dynamics}. We also investigate $G$-Dyson Brownian motion and the $G$-Wishart process as two important special cases. In Section~\ref{sec:hign limits}, we establish the relative compactness of the empirical spectral measure processes and characterize their high-dimensional limits, as stated in Theorem~\ref{thm: pre-compact} and Theorem~\ref{thm:limit-equation}, respectively. Furthermore, for $G$-Dyson Brownian motion and the $G$-Wishart process, we prove uniqueness of solutions to the corresponding limiting equations, thereby obtaining the semicircle and Marchenko--Pastur laws under $G$-expectation, respectively.

\section{Preliminaries on nonlinear expectations}\label{sec:pre}

In this section, we provide some preliminaries on nonlinear expectations, with particular emphasis on the theory of $G$-expectation and $G$-Brownian motion. We refer to \cite{DenisHuPeng2011,hu2019levy,Peng2007,Peng2008,peng2019nonlinear,gpp18,liu2020exit} for further details on Sections~\ref{sec:2.1}--\ref{sec:2.5}. In Section~\ref{sec:Kolmog}, we construct an infinite-dimensional $\mathbi G$-Brownian motion and establish a Kolmogorov extension theorem for capacities adapted to the present setting. We refer to \cite{Ibr13,dkn18} for related results.


\subsection{\texorpdfstring{Sublinear expectation and $G$-normal distribution}{}}\label{sec:2.1}

Let $\Omega$ be a given set and $\mathcal{H}$ be a linear space of real-valued functions on $\Omega$ such that $c \in \mathcal{H}$ for all constants $c$ and $|X| \in \mathcal{H}$ if $X \in \mathcal{H}$. We further suppose that if $X_{1}, \cdots, X_{n} \in \mathcal{H}$, then $\varphi\left(X_{1}, \cdots, X_{n}\right) \in \mathcal{H}$ for each $\varphi \in C_{\text {b,Lip}}\left(\mathbb{R}^{n}\right)$, where $C_{\text{b,Lip}}\left(\mathbb{R}^{n}\right)$ denotes the space of all bounded  Lipschitz-continuous functions on $\mathbb{R}^{n}$.

\begin{definition}\label{def:Gexp}
A functional $\E: \mathcal{H} \rightarrow \mathbb{R}$ is   a sublinear expectation if the following properties are satisfied: for each $X, Y \in \mathcal{H}$, 
\begin{itemize}
\item[$(1)$] $X \geq Y$ implies $\E[X] \geq \E[Y]$ (monotonicity);

\item[$(2)$]  $\E[c]=c$ for $c \in \mathbb{R}$ (constant preserving);

\item[$(3)$]  $\E[X+Y] \leq \E[X]+\E[Y]$ (sub-additivity);

\item[$(4)$]  $\E[\lambda X]=\lambda \E[X]$ for $\lambda>0$ (positive homogeneity).
\end{itemize}
\noindent The triple $(\Omega, \mathcal{H}, \E)$ is called a sublinear expectation space.
\end{definition}

 When only (1) and (2) are satisfied, $\E$ is called a nonlinear expectation and  $(\Omega, \mathcal{H}, \E)$ is therefore called a nonlinear expectation space. If equality holds in (3), then $\mathbb E$ reduces to a linear expectation which is denoted by $E$.

For  $p \geq 1$, suppose $|X|^{p} \in \mathcal{H}$ for each $X \in \mathcal{H}$, and  let $\mathbb{L}^{p}(\Omega)$ be the completion of $\mathcal{H}$ under the norm $\|X\|_{\mathbb{L}^{p}}:=\left(\E\left[|X|^{p}\right]\right)^{1 / p}.$ Then, we have $\mathbb{L}^{p}(\Omega)\subset \mathbb{L}^{p'}(\Omega)$ for  $p\ge p'$. Clearly $\big(\Omega, \mathbb{L}^{1}\left(\Omega\right),  \E\big)$ forms a sublinear expectation space, and we call it
a complete sublinear expectation space.

\begin{definition}
Let $X_{1}$ and $X_{2}$ be two $n$-dimensional random vectors on complete sublinear expectation spaces $\left(\Omega_{1}, \mathbb{L}^{1}\left(\Omega_{1}\right), \E_{1}\right)$ and $\left(\Omega_{2}, \mathbb{L}^{1}\left(\Omega_{2}\right), \E_{2}\right)$, respectively. They are called identically distributed, denoted by $X_{1} \stackrel{d}{=} X_{2}$, if $\E_{1}\left[\varphi\left(X_{1}\right)\right]=\E_{2}\left[\varphi\left(X_{2}\right)\right]$ for each $\varphi \in C_{\mathrm{b,Lip }}\left(\mathbb{R}^{n}\right)$.
\end{definition}

\begin{definition}
Let $\left(\Omega, \mathbb{L}^{1}(\Omega), \E\right)$ be a complete sublinear expectation space. An $m$-dimensional random vector $Y$ is said to be independent from an $n$-dimensional random vector $X$ if   $\E[\varphi(X, Y)]=\E\big[\E[\varphi(x, Y)]\big|_{x=X}\big],
$ for each $\varphi \in C_\mathrm{b,Lip}\left(\mathbb{R}^{n+m}\right)$.
\end{definition}

Denote by $\mathbb{S}_d$  the space of all $d \times d$ symmetric matrices, and  let $G: \mathbb{S}_d \rightarrow \mathbb{R}$ be a  monotone sublinear function which satisfies, for $A, A^{\prime} \in \mathbb{S}_d$,
\begin{equation*}
\left\{
\begin{array}{l}
G(A) \geq G\left(A^{\prime}\right) \text { if } A \geq A^{\prime}, \\
G\left(A+A^{\prime}\right) \leq G(A)+G\left(A^{\prime}\right), \\
G(\lambda A)=\lambda G(A) \text { for } \lambda \geq 0.
\end{array}
\right.
\end{equation*}
By the Hahn-Banach theorem, there exists a bounded, convex, and closed subset $\Gamma \subset \mathbb{S}^+_d$ such that
\begin{equation*}
G(A)=\frac{1}{2} \sup _{\gamma \in \Gamma} \operatorname{tr}[\gamma A], \text { for } A \in \mathbb{S}_d,
\end{equation*}
where $\mathbb{S}^{+}_d$ denotes the space of all nonnegative definite matrices in $\mathbb{S}_d$. 


The following $G$-normal distribution  and maximal distribution arise naturally from the central limit theorem and the law of large numbers in the framework of sublinear expectations, respectively.

\begin{definition}[$G$-normal distribution]
Let $(\Omega, \mathbb{L}^1(\Omega), \E)$ be a complete sublinear expectation space and let $G:\mathbb S_d\to \R$ be a monotone sublinear function.
A $d$-dimensional random vector $X$ is called $G$-normally distributed, if for each $\varphi \in C_{\mathrm{b,Lip}}(\mathbb{R}^d)$, $\E[\varphi(X)] = u^{\varphi}(1, 0)$,
where $u^{\varphi}$ is the viscosity solution of the following $G$-heat equation:
\begin{equation}\label{e:G-heat}
\begin{cases}
\partial_t u(t,x) - G\left(D_x^2 u(t,x)\right) = 0, \quad t>0, \, x\in\R^d, \\[3pt]
u(0, x) = \varphi(x),
\end{cases}
\end{equation}
where $D_x^2 := (\partial^2_{x_ix_j} )_{d\times d}$ denotes the Hessian matrix. 
\end{definition}
\begin{definition}[maximal distribution]\label{def:maximal}
Let $(\Omega, \mathbb{L}^1(\Omega), \E)$ be a complete sublinear expectation space. A $d$-dimensional random vector $\eta$ is called maximally distributed, if there exists a bounded, closed and convex subset $\Gamma \subset \mathbb{R}^d$ such that 
\begin{equation*}
    \mathbb{E}[\varphi(\eta)] = \max_{y\in \Gamma} \varphi(y), \quad \varphi \in C_\mathrm{b,Lip}(\mathbb{R}^d).
\end{equation*}
\end{definition}

\subsection{\texorpdfstring{$G$-Brownian motion and L\'evy's characterization}{}} \label{sec:G-levy}

\begin{definition}[$G$-Brownian motion]
Let $(\Omega, \mathbb{L}^{1}(\Omega), \E)$ be a complete sublinear expectation space and let $G:\mathbb S_d\to \R$ be a monotone sublinear function. A $d$-dimensional process $\left(B_{t}\right)_{t \geq 0}$ with $B_{t} \in\left(\mathbb{L}^{1}(\Omega)\right)^{d}$ for each $t \geq 0$ is called a G-Brownian motion if the following properties are satisfied:
\begin{itemize}
    \item[$(1)$] $B_{0}=0$;

\item[$(2)$] for  $t,s\ge 0$, $B_{t+s}-B_t$  is independent from $\left(B_{t_{1}}, \cdots, B_{t_{n}}\right)$ for each $n \in \mathbb{N}$ and $0 \leq t_{1} \leq \cdots \leq t_{n} \leq t$;

\item[$(3)$] for $t, s \geq 0$, $B_{t+s}-B_{t} \stackrel{d}{=} \sqrt{s}\xi$, where $\xi$ is G-normally distributed. 
\end{itemize}

\noindent Moreover, if $\E[B_t] = \E[-B_t] = 0$, then we call $(B_t)_{t \geq 0}$  a symmetric $G$-Brownian motion.
\end{definition}

\begin{theorem}\label{thm:1.6}
Let $\left(B_{t}\right)_{t \geq 0}$ be a $d$-dimensional  symmetric $G$-Brownian motion.  Then, for each fixed $\varphi \in C_\mathrm{b, Lip}(\mathbb{R}^{d})$, the function
$$u(t, x):=\E\left[\varphi\left(x+B_{t}\right)\right], \quad(t, x) \in[0, \infty) \times \mathbb{R}^{d}$$
is the viscosity solution to the $G$-heat equation \eqref{e:G-heat}, 
with $G$  given by
$$G(A)=\frac{1}{2} \E\left[\left\langle A B_{1}, B_{1}\right\rangle\right], \quad A \in \mathbb{S}_d.$$
\end{theorem}

When $d>1$, one can verify that the first $(d-1)$ components of $B$, denoted by  $\tilde B=(B^1, \dots, B^{d-1})$, is a $(d-1)$-dimensional $\tilde G$-Brownian motion, where $\tilde G:\mathbb S_{d-1}\to \R$ is a monotone sublinear function given by 
\[
\tilde G(A):=G\left(
\begin{pmatrix}
A & 0\\
0 & 0
\end{pmatrix}
\right) \text{ for } A\in\mathbb{S}_{d-1}.
\]
A similar result holds for the projection to the first $k$ components of $B$ for any positive integer $k<d$.

\begin{proposition}\label{prop: transform-G}
Let $B$ be a $d$-dimensional $G$-Brownian motion. For any deterministic $n\times d$ matrix $Q$, define
$X_t = QB_t$. Then $X_t$ is a $G_Q$-Brownian motion with the generator given by $G_Q(A)=G(Q^\top A Q)$ for $A\in \mathbb S_n$.
\end{proposition}
\begin{proof}
Let $$G_Q(A):=\frac{1}{2}\E[\langle AQB_1, QB_1 \rangle], \quad A \in \mathbb{S}_n.$$
It is easily seen that $X_t$ is a $G_Q$-Brownian motion, and we can verify that
$$G_Q(A)=\frac{1}{2}\E[\langle Q^{\top} AQB_1, B_1 \rangle] = G(Q^{\top} AQ).$$
    \end{proof}

For each $t \geq 0$, set
\begin{equation*}
\Lip\left(\Omega_{t}\right):=\left\{\varphi\left(B_{t_{1}}, B_{t_{2}}, \cdots, B_{t_{n}}\right): n \in \mathbb{N}, 0 \leq t_{1}<\cdots<t_{n} \leq t, \varphi \in C_\text{b,Lip}(\mathbb{R}^{n \times d})\right\};
\end{equation*}
and 
\begin{equation*}
\Lip(\Omega):=\bigcup_{m=1}^{\infty} \Lip\left(\Omega_{m}\right). 
\end{equation*}
For each $p \geq 1$, we denote by $L_{G}^{p}\left(\Omega_{t}\right)$ (resp., $L_{G}^{p}(\Omega)$) the completion of $\Lip\left(\Omega_{t}\right)$ (resp., $\Lip(\Omega)$) under the norm $\|X\|_{L_{G}^{p}}:=\left(\E\left[|X|^{p}\right]\right)^{\frac{1}{p}}$ for each $t \geq 0$.  The sublinear expectation $\E: L_G^1(\Omega)\to \R$ is now called a $G$-expectation.

For $\varphi\left(B_{t_{1}}, B_{t_{2}}-B_{t_{1}}, \cdots, B_{t_{n}}-B_{t_{n-1}}\right) \in \Lip(\Omega)$,  its conditional expectation at time $t_{i}$ is defined  by
\begin{equation*}
\E_{t_{i}}\left[\varphi\left(B_{t_{1}}, B_{t_{2}}-B_{t_{1}}, \cdots, B_{t_{n}}-B_{t_{n-1}}\right)\right]:=\psi\left(B_{t_{1}}, \cdots, B_{t_{i}}-B_{t_{i-1}}\right),
\end{equation*}
for $i=1, \dots, n$, where
$$\psi\left(x_{1}, \cdots, x_{i}\right):=\E\left[\varphi\left(x_{1}, \cdots, x_{i}, B_{t_{i+1}}-B_{t_{i}}, \cdots, B_{t_{n}}-B_{t_{n-1}}\right)\right].$$ 
Then,  the mapping $\E_{t_{i}}:\Lip(\Omega) \to \Lip (\Omega_{t_i})$ satisfies  conditions (1)--(4) in  Definition~\ref{def:Gexp} and can be continuously extended to the mapping which is still denoted by  $\E_{t_i}: L_{G}^{1}(\Omega)\to  L_{G}^{1}\left(\Omega_{t_{i}}\right)$. We call $(\Omega, L_{G}^{1}(\Omega),(L_{G}^{1}(\Omega_{t}))_{t \geq 0},(\E_{t})_{t \geq 0})$ a  $G$-expectation space.

\begin{definition}
For $p\ge 1$, let $\mathbf M_G^{p,0}(0, T)$ denote the space of processes of the following form:
$$\eta_{t}(\omega)=\sum_{i=0}^{n-1} \xi_{i}(\omega) \mathbf I_{\left[t_{i}, t_{i+1}\right)}(t),$$
where  $\left\{t_{0}, t_1, \dots, t_{n}\right\}$ is a partition of $[0, T]$ and $\xi_{i} \in L^p_G\left(\Omega_{t_{i}}\right)$ for $i\in\{0,1, \dots, n-1\}$. We denote by $\mathbf M_G^{p}(0, T)$ the completion of $\mathbf M_G^{p, 0}(0, T)$ under the norm $$\|\eta\|_{\mathbf M_G^{p}}:=\left(\E\left[\int_{0}^{T}\left|\eta_{t}\right|^{p} \mathrm{d} t\right]\right)^{\frac{1}{p}}.$$
\end{definition}

\begin{definition}
A $d$-dimensional process $M_{t}=\left(M_{t}^{1}, \cdots, M_{t}^{d}\right)^{\top}$ is called a martingale if for any $t \in[0, \infty)$, $M^i_{t} \in L_{G}^{1}\left(\Omega_{t}\right)$ and $\E_s\left[M_{t}^{i}\right]=M_{s}^{i}$ for all $s \leq t$ and $i \in\{1,\cdots, d\}$. Furthermore, a $d$-dimensional martingale $M$ is called symmetric if $\E\left[M_{t}^{i}\right]=-\E\left[-M_{t}^{i}\right]$ for all $t$ and  $i$.
\end{definition}

Clearly, the $G$-Brownian motion $B$ is a symmetric martingale. For a symmetric martingale $M$ satisfying 
\begin{equation}\label{e:con-M}
\text{$M_0=0, M_t \in L_G^2(\Omega_t)$ and } \E_t[|M_{t+s}-M_t|^2]\le C_M s,
\end{equation}
for some constant $C_M$ for all $s,t \ge0$, we can define stochastic integrals with respect to $M$ as follows.

For $\eta_{t}(\omega)=\sum_{i=0}^{n-1} \xi_{i}(\omega) I_{\left[t_{i}, t_{i+1}\right)}(t)\in \mathbf M_G^{2,0}(0,T)$, define
\begin{equation}\label{e:def-ito-int}
I(\eta)= \int_0^T\eta_s \dd M_s := \sum_{i=0}^{n-1} \xi_i (M_{t_{i+1}}-M_{t_i}). 
\end{equation}
Then the mapping $I: \mathbf M_G^{2, 0}(0,T)\to L_G^2(\Omega_T)$ can be continuously extended to the mapping from $\mathbf M_G^{2}(0,T)$ to $L_G^2(\Omega_T)$. Moreover, for each $\eta\in \mathbf M_G^2(0,T)$, the stochastic integral $\int_0^t \eta_s \dd M_s$ is a symmetric martingale and 
\[\E \left[\left|\int_0^T \eta_s \dd M_s\right|^2\right]\le C_M \E \left[\int_0^T |\eta_s|^2 \dd s\right].\]


Similar to the classical case, we have the following $G$-It\^o isometry and $G$-Burkholder-Davis-Gundy inequality  for integrals against $G$-Brownian motion. 
Let $B=(B^1,\ldots,B^d)^\top$ be a $d$-dimensional
$G$-Brownian motion. For each fixed $\mathbf a\in\mathbb R^d$,
denote $B_t^{\mathbf{a}}:= \langle \mathbf{a}, B_t \rangle=
\sum_{i=1}^d a_i B_t^i$.
\begin{theorem}[$G$-It\^o isometry] 
Let $\eta\in \mathbf M_G^2(0,T)$. Then,  we have  $$\mathbb{E}\left[\left(\int_0^T \eta_t \mathrm{d}B^{\mathbf{a}}_t\right)^2\right]=\mathbb{E}\left[\int_0^T |\eta_t|^2 \mathrm{d}\langle B^{\mathbf{a}}\rangle_t\right].$$
\end{theorem}
\begin{theorem}[$G$-Burkholder-Davis-Gundy inequality]
For any $\eta\in \mathbf M_G^2(0,T)$ and $p>0$, there exist constants $c_p$ and $C_p$ with $0<c_p<C_p<\infty$, depending only on $p$, such that
\begin{equation*}
c_p \mathbb{E}\left[ \left( \int_{0}^{T} |\eta_s|^2 \mathrm{d}\langle B^{\mathbf{a}} \rangle_s \right)^{p/2} \right]\le \mathbb{E}\left[ \sup_{0\le t\le T} \left| \int_{0}^{t} \eta_s \mathrm{d}B^{\mathbf{a}}_s \right|^p \right]\le C_p \mathbb{E}\left[ \left( \int_{0}^{T} |\eta_s|^2 \mathrm{d}\langle B^{\mathbf{a}} \rangle_s \right)^{p/2} \right].
\end{equation*}
\end{theorem}

For a $d$-dimensional martingale $M_{t}=\left(M_{t}^{1}, \cdots, M_{t}^{d}\right)^{\top}$, we denote by $\la M \ra_t = (\la M^i , M^j\ra_t)_{d\times d}$ its quadratic covariation  matrix, where $\la M^i , M^j\ra_t:=\lim_{|\pi_n|\to 0}\sum_{k=0}^{n-1} (M^i_{t_{k+1}}-M^i_{t_{k}})(M^j_{t_{k+1}}-M^j_{t_{k}})$, with $|\pi_n|:=\max\{|t_{k+1}-t_k|: k=0,1,\dots,n-1\}$, is the quadratic covariation between $M^i$ and $M^j$. When $i=j$, we write $\la M^i\ra_t:=\la M^i, M^i\ra_t$ for the corresponding quadratic variation.   We remark that $\int_0^t \eta_s \dd \la M\ra_s$ is well-defined for $\eta\in \mathbf M_G^1(0,T)$ if the condition \eqref{e:con-M} is satisfied. 

The following is  L\'{e}vy's characterization of $G$-Brownian motion established in  \cite{hu2019levy}.

\begin{theorem}\label{thm:levy2}

On a $G$-expectation space $(\Omega, L_{G}^{1}(\Omega),(L_{G}^{1}(\Omega_{t}))_{t \geq 0},(\E_{t})_{t \geq 0})$, let $(M_{t})_{t \geq 0}$ be a $d$-dimensional symmetric martingale satisfying $M_{0}=0$ and $M_{t} \in\left(L_{G}^{2}\left(\Omega_{t}\right)\right)^{d}$ for each $t \geq 0$. Given a monotone sublinear function $\bar G: \mathbb{S}_d \to\mathbb{R}$, the following conditions are equivalent:
\begin{itemize}
    \item[$(1)$] $(M_{t})_{t\ge 0}$ is a $d$-dimensional $\bar G$-Brownian motion;


\item[$(2)$] $\E_{t}\left[\left|M_{t+s}-M_{t}\right|^{2}\right] \leq C s$ for all  $s,t \geq 0$ for some $C>0$,  and the process $\frac{1}{2} \operatorname{tr}\left[A\langle M\rangle_{t}\right]-\bar G(A) t$ is a martingale for each $A \in \mathbb{S}_d$. 
\end{itemize}
\end{theorem}

The following result provides a criterion for a symmetric martingale to be a  $G$-Brownian motion satisfying Assumption~\ref{assum:G}.

\begin{corollary}\label{corollary:G-brownian motion}
Let $(M_t)_{t\ge 0}$ be a $d$-dimensional symmetric martingale as in Theorem~\ref{thm:levy2}. Let $(B_t)_{t\ge 0}$ be a $1$-dimensional $\mathfrak g$-Brownian motion, where $\mathfrak g$ is given by~\eqref{e:G-tilde}. If $\left \langle M_{i},M_{j}  \right \rangle_{t}=\mathbf 1_{\{i=j\}}\left \langle B \right \rangle _{t}$ for all $i, j\in\{1, \dots, d\}$ and $t\geq 0$, then $(M_t)_{t\ge 0}$ is a $d$-dimensional $\bar G$-Brownian motion, with $\bar G$ satisfying Assumption~\ref{assum:G}.
\end{corollary}
\begin{proof}
Noting  that for $s,t>0$ we have
$$ \E_{t}\left[\left|M \left ( t+s\right )-M\left ( t \right )\right|^{2}\right]= \E_{t}\left[\sum_{i=1}^{d} \left|M_{i}\left ( t+s\right )-M_{i} \left ( t \right )\right|^{2}\right]\le\sum_{i=1}^{d}\E_{t}\left[ \langle B\rangle_{t+s}-\langle B\rangle_{t}\right]\leq d \,\overline{\sigma}^{2} s,$$ by Theorem~\ref{thm:levy2}, it suffices to show that  $\frac{1}{2} \operatorname{tr}\left[A\langle M\rangle_{t}\right]-\bar G(A) t$ is a martingale for each  $A \in \mathbb{S}_{d}$.  

 Indeed, for $s<t$, we have 
 \begin{align*}
  &\E_{s} \left[\frac{1}{2} \operatorname{tr}[A\langle M\rangle_t]-\bar G(A)t\right] =\E_{s}\bigg[ \frac12\sum_{i=1}^d a_{ii}\langle M_{i}\rangle_t- \bar G(A)t  \bigg]\\
&= \E_{s}\bigg[ \frac12 \sum_{i=1}^d a_{ii} \big(\langle B\rangle_t - \langle B\rangle_s\big)   \bigg] + \frac12\bigg( \sum_{i=1}^d a_{ii}\langle M_i\rangle_s \bigg) - \bar G(A)t \\[4pt]
&= \frac12\Bigg( \bar{\sigma}^2\bigg( \sum_{i=1}^d a_{ii} \bigg)^+ - \underline{\sigma}^2\bigg( \sum_{i=1}^d a_{ii} \bigg)^- \Bigg)(t-s)  + \frac12\bigg( \sum_{i=1}^d a_{ii}\langle M_i\rangle_s \bigg) - \bar G(A)t \\[4pt]
&= \frac12\bigg( \sum_{i=1}^d a_{ii}\langle M_i\rangle_s \bigg) -\bar G(A)s= \frac{1}{2} \operatorname{tr}[A\langle M\rangle_s]-\bar G(A)s.
\end{align*}
This completes the proof. 
\end{proof}

Similar to Assumption~\ref{assum:G}, it is also natural to impose the following diagonal structure assumption on the generator $G$. 

\begingroup
\addtocounter{assumption}{-2} 
  \renewcommand{\theassumption}{\arabic{assumption}'} 
  \renewcommand{\theHassumption}{\arabic{assumption}'} 
  
\begin{assumption}
\label{assum:G'}
Assume 
$G(A)
=
\sum_{i=1}^d \mathfrak g(A_{ii})$ for 
$A=(A_{ij})_{1\leq i,j\leq d}\in\mathbb{S}_d$,
where $\mathfrak g$ is given by~\eqref{e:G-tilde}. 
\end{assumption}
\endgroup
\addtocounter{assumption}{1}
\begin{corollary}
\label{G property}
Let $(B_{t})_{t\ge 0}$ be a $d$-dimensional $G$-Brownian motion  satisfying Assumption~\ref{assum:G'}, and  $\alpha=(\alpha^1,\alpha^2,\dots,\alpha^d)\in \mathbf M_G^{2}(0, T)^d$ be a $d$-dimensional process satisfying $\left\|\alpha_t\right\|^2=(\alpha^1_t)^2+(\alpha^2_t)^2+\dots+(\alpha^d_t)^2=1$ for all $t\in [0,T]$. Then,  the process $M_t:=\int_0^t\alpha_r \mathrm{d}B_{r}$ is a one-dimensional $\mathfrak g$-Brownian motion.
\end{corollary}
\begin{proof}
Obviously, $M$ is a symmetric martingale, noting the boundedness of  the integrand $\alpha$. 
To show  that $M$ is a $\mathfrak g$-Brownian motion, by Theorem~\ref{thm:levy2}, 
it suffices to  verify that both $\frac{1}{2}\langle M \rangle_t-\frac{1}{2}\overline{\sigma}^2t$ and $\frac{1}{2}\underline{\sigma}^2t-\frac{1}{2}\langle M \rangle_t$ are $\mathfrak g$-martingales.

Indeed, for $0\leq s<t$,
\begin{align*}
\E_s[\langle M\rangle_t-\overline{\sigma}^2t]=&\E_s[\langle M \rangle_t-\langle M \rangle_s]+\langle M\rangle_s-\overline{\sigma}^2t\\
=&\E_s\left[\sum_{i=1}^d\int_s^t(\alpha_r^i)^2\mathrm{d}\langle B^i\rangle_r\right]+\langle M\rangle_s-\overline{\sigma}^2t\\
=&\overline{\sigma}^2\E_s\left[\sum_{i=1}^d\int_s^t(\alpha_r^i)^2\mathrm{d}r\right]+\langle M\rangle_s-\overline{\sigma}^2t\\
=&\langle M\rangle_s-\overline{\sigma}^2s,
\end{align*}
where in the last step we use the assumption $\|\alpha_r\|^2=1$.  Similarly, we show that $\frac{1}{2}\underline{\sigma}^2t-\frac{1}{2}\langle M \rangle_t$ is also a $\mathfrak g$-martingale. The proof is completed.
\end{proof}

\begin{remark}\label{rem:Levy}    
Let $(B_{t})_{t\ge 0}$ be a $d$-dimensional $G$-Brownian motion and $\alpha^{(n)}=(\alpha_1^{(n)},\alpha_2^{(n)},\dots,\alpha_n^{(n)})^\top$ be a $(n\times d)$-matrix-valued  process, where $\alpha_i^{(n)}\in \mathbf M_G^2(0,T)^d, i=1, \dots, n$, are $d$-dimensional processes  satisfying $\|\alpha_i^{(n)}(t)\|^2=1$ and $\langle \alpha^{(n)}_i(t), \alpha^{(n)}_j(t)\rangle  =0$ for each $t\in [0,T]$ and $i\neq j$. 
   
   Under Assumption \ref{assum:G'}, Corollary \ref{G property} implies that each component of $\int_0^t \alpha^{(n)}(r) \,\mathrm{d}B_r$ is a $\mathfrak g$-Brownian motion. However, in general $W_t:=\int_0^t \alpha^{(n)}(r) \,\mathrm{d}B_r$
    is not a $G$-Brownian motion for any $G$,  as illustrated in the following example for $d=3,n=2$ when $\underline{\sigma}<\overline{\sigma}$: 
 \begin{example}\label{example: not same $G$}
Let $B$ be a $3$-dimensional $G$-Brownian motion and  $\alpha^{(2)}(t):=A$ be a constant process with 
\begin{equation}\label{e:A}
A:=\dfrac{1}{\sqrt{2}}
\begin{pmatrix}
1 & -1& 0\\[2pt]
1& 1 & 0
\end{pmatrix}.
\end{equation}
Then,   $W_t=\int_0^t \alpha^{(2)}(r) \,\mathrm{d}B_r=Q(B_t^1, B_t^2)^\top$, where 
$Q
=
\frac{1}{\sqrt{2}}
\begin{pmatrix}
1 & -1\\
1 & 1
\end{pmatrix}$ is an orthogonal matrix. 
Under Assumption \ref{assum:G'}, clearly $\left(B_{t}^{1},B_{t}^{2}\right)^\top$ is a $2$-dimensional $G$-Brownian motion with the same form of $G$. However,  Proposition~\ref{prop: transform-G} implies that $Q(B_t^1, B_t^2)^\top$ is a $2$-dimensional $G$-Brownian motion with generator $G_Q$, where $G_Q(M) = G(Q^{\top}MQ)$ for $M\in\mathbb S_2$. Note that $G_Q\neq G$ since  $G_Q(D) =0\neq \frac{1}{2}
\left(
\overline{\sigma}^{\,2}
-
\underline{\sigma}^{\,2}
\right)=G(D)$ for 
$D=\begin{pmatrix}
1 & 0\\
0 & -1
\end{pmatrix}$. Now, let $\alpha^{(2)}(t) = A\mathbf I_{[0,T/2)}(t) + C\mathbf I_{[T/2, T]}(t)$ with $A$ given 
by \eqref{e:A} and $C:= \begin{pmatrix}
1 & 0& 0\\
0& 1 & 0
\end{pmatrix}
$, then $W$ is a $G_Q$-Brownian motion on $[0,T/2)$ while $(W_{\frac T2+s}-W_s)_{0\le s\le \frac{T}{2}}$ is a $G$-Brownian motion.
\end{example}  

By contrast, if the generator $G$ satisfies Assumption~\ref{assum:G}, then Corollary~\ref{corollary:G-brownian motion}, together with the fact that the quadratic covariation is given by
$$
\left\langle
\int_0^{\cdot} \alpha^{(n)}_i(r) \,\mathrm{d}B_r,
\int_0^{\cdot} \alpha^{(n)}_j(r) \,\mathrm{d}B_r
\right\rangle_t
=
\int_0^t
\alpha^{(n)}_i(r)
\big(\alpha^{(n)}_j(r)\big)^\top
\,\mathrm{d}\left\langle B^{1} \right\rangle_r
=
\mathbf I_{\{i=j\}}
\left\langle B^{1} \right\rangle_t,
$$
where the second equality is due to~\eqref{e:diagonal-prop},  implies that $\int_0^t \alpha^{(n)}(r)\mathrm{d}B_r$ is an $n$-dimensional $\bar G$-Brownian motion with $\bar G$ still satisfying Assumption~\ref{assum:G}.
\end{remark}

\subsection{\texorpdfstring{$G$-expectation}{} on the path space}\label{sec:G-exp} Let  $\Omega:=C_{0}\left(\mathbb{R}^{+};\R^d\right)$ be the space of $\mathbb{R}^{d}$-valued continuous paths $\left(\omega_{t}\right)_{t\ge 0}$ with $\omega_{0}=0$, endowed with the distance
\begin{equation*}\label{e:rho}
\rho\left(\omega^{1}, \omega^{2}\right):=\sum_{i=1}^{\infty} 2^{-i}\Big[\big(\max _{t \in[0, i]}\left|\omega_{t}^{1}-\omega_{t}^{2}\right|\big) \wedge 1\Big].
\end{equation*}
Then $(\Omega, \rho)$ is a Polish space, and we use $\mathcal{B}(\Omega)$ to denote its Borel $\sigma$-algebra.  For each fixed $T \in[0, \infty)$,  denote $\Omega_{T}= \left\{\omega_{. \wedge T}: \omega \in \Omega\right\}$. The canonical process $(B_{t}(\omega)=\omega_{t})_{t \ge0}$ for $\omega \in \Omega$ is then a $G$-Brownian motion on the $G$-expectation space $(\Omega, L_{G}^{1}(\Omega),(L_{G}^{1}(\Omega_{t}))_{t \geq 0},(\mathbb{E}_{t})_{t \geq 0})$.

In this subsection, the generator $G$ is given by  
\begin{equation}\label{e:G-Theta}
G(A):=\frac12\sup_{\gamma\in\Theta}\tr[\gamma \gamma^\top A], \quad A\in \mathbb S_d.
\end{equation}
where $\Theta$ is a bounded and closed subset in $\mathbb{R}^{d \times d}$. It is easy to verify that one can choose 
\begin{equation}\label{e:Theta1}
\Theta:=\{\sigma I_d: \underline{\sigma}\le \sigma\le \overline{\sigma}\}
\end{equation}
for $G$ satisfying Assumption~\ref{assum:G} and 
\begin{equation}\label{e:Theta2}
\Theta:=\{\mathrm{diag}(\sigma_1, \dots, \sigma_d): \underline{\sigma}\le \sigma_i\le \overline{\sigma}, i=1,\dots, d\}
\end{equation}
for $G$ satisfying Assumption~\ref{assum:G'}.

Let $(\Omega, \mathcal{F}, \mathbb F:=(\mathcal F_t), P)$ be the classical Wiener space, on which  $\left(W_{t}\right)_{t \geq 0}=\left(W_{t}^{i}\right)_{i=1, t \geq 0}^{d}$ is a $d$-dimensional Wiener process.  We denote by $\mathcal{A}_{t, T}^{\Theta}$, the space of all $\Theta$-valued $\mathbb F$-adapted processes on an interval $[t, T]$. For each fixed $\theta \in \mathcal{A}_{0, \infty}^{\Theta}:=\cup_{n=1}^\infty\mathcal A^\Theta_{0,n}$,  we denote by $P_\theta$ the law of the martingale $\int^{\cdot}_{0} \theta_{s} \mathrm{d} W_{s}$. Denote $\mathcal P=\left\{P_{\theta}: \theta \in \mathcal{A}_{0, \infty}^{\Theta}\right\}$ and let $\mathscr{P}=\overline{\mathcal P}$ be the closure of $\mathcal P$ under the topology of weak convergence. Then,  $\mathcal P$ is tight and hence $\mathscr{P}$ is weakly compact.

Let  $L^0(\Omega)$ denote the space of all $\mathcal B(\Omega)$-measurable random variables.  For   $X \in  L^0(\Omega)$ satisfying $E_{P}[X]$ exists for each $P \in \mathscr{P}$,  define the \emph{upper expectations} as follows
\begin{equation}\label{e:upper-exp}
\bar{\mathbb{E}}[X]:=\sup _{P \in \mathcal P} E_{P}[X] \text{ and }
 \hat{\mathbb{E}}[X]:=\sup _{P \in \mathscr{P}} E_{P}[X].
\end{equation}
Then, we have 
\begin{equation}\label{e:equal E}
\mathbb{E}[X]  =\bar \E[X]= \hat{\mathbb{E}}[X], \quad \text{ for all } X \in L_{G}^{1}(\Omega),
\end{equation}
where we recall that the space $L_G^1(\Omega)$ is defined in Section~\ref{sec:G-levy}. Similarly, we define the \emph{upper probability} 
\begin{equation*}
c(A):=\sup _{P \in \mathscr{P}} P(A)=\hat \E(\mathbf I_A), \quad A \in \mathscr{B}(\Omega),
\end{equation*}
which is indeed a \emph{Choquet capacity} since $\mathscr P$ is weakly compact.

\begin{definition}
For a set $A\in\mathcal B(\Omega)$, it is called polar if $c(A)=0$. We say that a property holds \emph{quasi-surely (q.s.)} if it holds outside a polar set.
\end{definition}

The following Borel--Cantelli lemma is useful for identifying polar sets.
\begin{lemma}\label{lem:B-C}
Let $\left(A_{n}\right)_{n \in \mathbb{N}}$ be a sequence of Borel sets satisfying $\sum_{n=1}^{\infty} c\left(A_{n}\right)<\infty,$ then $\lim \sup _{n \rightarrow \infty} A_{n}$ is a polar set.
\end{lemma}

\begin{definition}
A random variable  $X\in L^0(\Omega)$ is said to be
\emph{quasi-continuous (q.c.)}, if
for any $\varepsilon>0$, there exists an open set $O$ with $c(O) < \varepsilon$ such that $X\mathbf \mid_{O^c}:O^c\to \R$ is a continuous mapping.

We say that a random variable $X$ has a q.c.\ version if there exists
a q.c.\ random variable $Y$ such that $X= Y$ q.s.
\end{definition}


We have the following characterizations for $L_G^1(\Omega_T)$ and  $L_G^1(\Omega)$, respectively:
\begin{equation}\label{e:LG1}
\begin{aligned}
L_{G}^{1}\left(\Omega_{T}\right) & =\left\{X \in L^0(\Omega_T): X \text { has a q.c.\ version, } \lim _{n \rightarrow \infty} \hat{\mathbb{E}}\left[|X| \mathbf{I}_{\{|X|>n\}}\right]=0\right\}; \\
L_{G}^{1}(\Omega) & =\left\{X \in L^0(\Omega): X \text { has a q.c.\ version, } \lim _{n \rightarrow \infty} \hat{\mathbb{E}}\left[|X| \mathbf{I}_{\{|X|>n\}}\right]=0\right\}.
\end{aligned}
\end{equation}
The above  equations remain valid with $\hat{\E}$  replaced by $\bar \E$.

\subsection{It\^o's formula}\label{sec:extended-G} It\^o's formula was first proved for $C^2$ functions whose second-order derivatives have polynomial growth (see, e.g., \cite{Peng2007}). The growth condition on the second-order derivatives was imposed due to the lack of a theory of stopping times, which made localization techniques unavailable. We remark that the stochastic analysis recalled in Sections~\ref{sec:G-levy} and \ref{sec:G-exp} mainly concerns the space $L_G^1(\Omega)$, whose random variables are q.c., whereas stopping times do not, in general, possess the q.c.\ property. 

In \cite{lp11},  It\^o's formula was extended to $C^2$ functions without any restriction on the growth of their second derivatives by analyzing the upper expectation $\hat\E$ defined in \eqref{e:upper-exp}. Note that $\hat \E$ coincides with the $G$-expectation $\E$ on the space $L_G^1(\Omega)$ (see eq.\ \eqref{e:equal E}).

For $p \geq 1$, let $\mathbf M^{b, 0}(0, T)$ denote  the space of processes of the form:
\begin{align*}
 \eta_{t}(\omega)=\sum_{j=0}^{n-1} \xi_{j}(\omega) \mathbf{I}_{\left[t_{j}, t_{j+1}\right)}(t).
\end{align*}
In the above, $\xi_{j} \in B_{b}(\Omega_{t_{j}})$ for $j\in\{0, \ldots, n-1\}$, where $B_b(\Omega_t)$ denotes the space of all bounded $\mathcal B(\Omega_t)$-measurable random variables.   We denote by $\mathbf M_*^{p}(0, T)$ the completion of $\mathbf M^{b, 0}(0, T)$ under the norm
$$\|\eta\|_{\mathbf M_*^{p}}=\left(\hat{\mathbb{E}}\left[\int_{0}^{T}\left|\eta_{t}\right|^{p} \mathrm{~d} t\right]\right)^{1 / p},$$
where $\hat \E$ is given in~\eqref{e:upper-exp}. The following $G$-Burkholder-Davis-Gundy inequality also holds:
\begin{theorem}[Generalized $G$-Burkholder-Davis-Gundy inequality]
For any $\eta \in M_*^2(0,T)$ and $p>0$, there exist constants
$c_p$ and $C_p$ with $0<c_p<C_p<\infty$, depending only on $p$,
such that
\[
\begin{aligned}
c_p\,
\hat{\mathbb{E}}\left[
    \left(\int_0^T |\eta_s|^2\,\mathrm d\langle B^a\rangle_s\right)^{p/2}
\right]\leq
\hat{\mathbb{E}}\left[
    \sup_{t\in[0,T]}
    \left|\int_0^t \eta_s\,\mathrm dB_s^{\mathbf{a}}\right|^p
\right]\leq
\,C_p\,
\hat{\mathbb{E}}\left[
    \left(\int_0^T |\eta_s|^2\,\mathrm d\langle B^a\rangle_s\right)^{p/2}
\right].
\end{aligned}
\]
\end{theorem}

We denote by $L_*^p(\Omega)$ (resp.\ $L_*^p(\Omega_t)$) the completion of $B_b(\Omega)$ (resp.\ $B_b(\Omega_t)$) under the norm $\|X\|_p:=(\hat\E[|X|^p])^{1/p}$. Then, $L_*^1(\Omega)$ has the following characterization:
\[L_{*}^{1}(\Omega)  =\left\{X \in L^0(\Omega): \lim _{n \rightarrow \infty} \hat{\mathbb{E}}\left[|X| \mathbf{I}_{\{|X|>n\}}\right]=0\right\}.
\]
In particular, random variables in $L_*^1(\Omega)$ are not necessarily q.c.\ (see \eqref{e:LG1} for a comparison with $L_G^1(\Omega)$).   It is also worth noting that, given a general weakly compact family $\mathscr P$ of probability measures, a random variable $X$ satisfying $\hat\E[|X|]=\sup_{P\in\mathscr P}E_P[|X|]<\infty$ need not satisfy the condition $ \lim _{n \rightarrow \infty} \hat{\mathbb{E}}\left[|X| \mathbf{I}_{\{|X|>n\}}\right]=0$ (see \cite[Examples~20 and 21]{DenisHuPeng2011}).

We refer to $(\Omega, L_{*}^{1}(\Omega),(L_{*}^{1}(\Omega_{t}))_{t \geq 0}
)$ as the \emph{extended $G$-expectation space}. Clearly, its restriction to $L_{G}^{1}(\Omega)$ coincides with the $G$-expectation space $(\Omega, L_{G}^{1}(\Omega),(L_{G}^{1}(\Omega_{t}))_{t \geq 0}
)$.

The mapping $I: \mathbf M^{b,0}(0,T)\to L_*^2(\Omega_T)$ defined as in~\eqref{e:def-ito-int} with $M$ replaced by $G$-Brownian motion $B$ can be continuously extended to $I: \mathbf M_*^{2}(0,T)\to L_*^2(\Omega_T)$, which also defines a stochastic integral of $\eta\in \mathbf M_*^{2}(0,T)$. We remark that, in this setting, conditional expectation has not yet been properly defined, and hence a martingale theory is not available. Noting that now $\eta_t$ is not necessarily a q.c.\ random variable, this extends the definition given in Section~\ref{sec:G-levy}. In particular, this extended definition of stochastic integrals allows us to define stochastic integrals stopped at a stopping time $\tau$:
$$\int_{0}^{t \wedge \tau} \eta_{s}\mathrm{~d} B_{s} = \int_{0}^{t} \mathbf{I}_{[0, \tau]}(s) \eta_{s}\mathrm{~d} B_{s} \quad \text{q.s.},$$
where the notion of a stopping time is defined below.

\begin{definition}\label{def:stopping time}
A stopping time $\tau$ with respect to the filtration $(\mathcal B(\Omega_t))_{t\ge 0}$ is a map on $\Omega$ with values in $[0, \infty]$ such that  $\{\tau \leq t\} \in \mathcal B(\Omega_t)$
for all $t\ge 0$.
\end{definition}

For $p\ge 1$, let $\mathbf M_{\mathrm{loc}}^p(0,T)$ denote the space of stochastic processes $\eta$ satisfying  that there exists a sequence $(\tau_m)_{m\in\mathbb N^+}$ of increasing stopping times  with $\tau_m\uparrow T$ q.s., such that $\eta\mathbf I_{[0,\tau_m]}\in\mathbf M_*^p(0,T)$ and $\int_0^T |\eta_s|^p \dd s<\infty$,  q.s.
Then, the definition of $\int_0^T \eta_s \dd B_s$ and $\int_0^T \xi_s \dd \la B\ra_s$ can be extended to $\eta\in\mathbf M_{\mathrm{loc}}^2(0,T)$ and $\xi\in \mathbf M_{\mathrm{loc}}^1(0,T)$, respectively.

We can now apply the localization technique to remove the growth condition on the second-order derivatives and obtain the following It\^o's formula (see \cite[Theorem 5.4]{lp11}).

\begin{proposition}[It\^o's formula]\label{Prop:18}
Let $\Phi \in C^{1,2}\left([0,T]\times\mathbb{R}^{n}\right)$ and $X=\left(X^{1}, \ldots, X^{n}\right)$ be an $n$-dimensional process on $[0, T]$ of the form: $$X_{t}^{m}=X_{0}^{m}+\int_{0}^{t} \alpha_{s}^{m} \mathrm{d} s+ \int_{0}^{t} \eta_{s}^{m i j} \mathrm{~d}\left\langle B^{i}, B^{j}\right\rangle_{s}+ \int_{0}^{t} \beta_{s}^{m i} \mathrm{~d} B_{s}^{i}, \quad m=1, \dots, n,$$ 
where $\alpha^{m}, \eta^{m i j} \in \mathbf M_{\mathrm{loc}}^{1}(0, T)$ and  $\beta^{m i} \in \mathbf M_{\mathrm{loc}}^{2}(0, T)$. Then,  for each $t \in[0, T]$, we have q.s.,
\begin{align*}
\Phi\left(t,X_{t}\right)-\Phi\left(0,X_{0}\right)= & \int_{0}^{t} \partial_{x^{m}} \Phi\left(u, X_{u}\right) \beta_{u}^{m i} \dd B_{u}^{i}+\int_{0}^{t} \left[\partial_t\Phi(u,X_u)+\partial_{x^{m}} \Phi\left(u, X_{u}\right) \alpha_{u}^{m} \right]\dd  u \\
& +\int_{0}^{t}\left[\partial_{x^{m}} \Phi\left(u, X_{u}\right) \eta_{u}^{m i j}+\frac{1}{2} \partial_{x^{m} x^{k}}^{2} \Phi\left(u, X_{u}\right) \beta_{u}^{m i} \beta_{u}^{k j}\right] \mathrm{d}\left\langle B^{i}, B^{j}\right\rangle_{u}.
\end{align*}
In the above equations we use the Einstein summation convention. 
\end{proposition}

\subsection{Weak and strong solutions}\label{sec:2.5}
Consider the following $G$-SDE:
\begin{equation}\label{eq:GSDE}
X_t
=
X_0
+\int_0^t b(s,X_s)\,ds
+\int_0^t h_{ij}(s,X_s)\,
d\langle B^i,B^j\rangle_s
+\int_0^t \sigma_j(s,X_s)\,dB_s^j,
\quad t\in[0,T],
\end{equation}
where the initial value $X_0$ is a  constant.

\begin{definition}[Weak solution]

A triple $(\mathfrak{G}, B, X)$ is called a weak solution to the G-SDE~\eqref{eq:GSDE} if $\mathfrak{G}:=(\Omega,L_G^1(\Omega), (L_G^1(\Omega_t))_{t\in[0,T]},(\E_t)_{t\in[0,T]})$ is a $G$-expectation space, $B$ is a $G$-Brownian motion on $\mathfrak G$, $X_t$ is  $\mathcal B(\Omega_t)$-measurable for each $t\in[0,T]$, and \eqref{eq:GSDE} holds q.s.
\end{definition}

Consider \eqref{eq:GSDE} on a prescribed $G$-expectation space $\mathfrak{G}:=(\Omega,L_G^1(\Omega), (L_G^1(\Omega_t))_{t\in[0,T]},(\E_t)_{t\in[0,T]})$, with a given $G$-Brownian motion $B$ on $\mathfrak G$. We define a strong solution as follows.

\begin{definition}[Strong solution]
A process $X$ is called a strong solution  to the G-SDE~\eqref{eq:GSDE} relative to $B$  if $(\mathfrak{G}, B, X)$ is a weak solution to~\eqref{eq:GSDE} and $X_t$ is $\sigma\{B_s:s\in[0,t]\}$-measurable for each $t\in[0,T]$.
\end{definition}

\begin{definition}[Pathwise uniqueness]\label{def:path-unique}
We say that pathwise uniqueness holds for the $G$-SDE~\eqref{eq:GSDE} if, whenever $X^1$ and $X^2$ are two (weak) solutions of \eqref{eq:GSDE} on the same $G$-expectation space, driven by the same $G$-Brownian motion $B$ and starting from the same initial condition, we have
$X_t^1=X_t^2$ for $t\in[0,T]$, q.s.
\end{definition}

\subsection{Kolmogorov extension theorem}\label{sec:Kolmog}

In order to investigate high-dimensional spectral limits considered in Section \ref{Main results on high-dimensional spectral limits}, we need to work on  the infinite-dimensional path space
\begin{equation}\label{eq:Omega-infinity}
\quad
\Omega^\infty:=C_0(\R^+; \R)^{\mathbb N}=\big\{\omega=(\omega^1,\omega^2,\ldots):\omega^i\in C_0(\R^+; \R) \text{ for }i\in\mathbb N\big\},
\end{equation}
 where an infinite-dimensional 
$\mathbi G$-Brownian motion $B=(B^1, B^2,\dots)$ can be constructed through the canonical process $B_t(\omega):=\omega(t),
t\in\R^+$. Here, $\mathbi G:=\{G_d\}_{d\in\mathbb N}$ is a family of monotone sublinear functions $G_d:\mathbb S_d\to \R$ satisfying Assumption~\ref{assum: consistent G} below. 
We denote by $\pi_d(\omega):=(\omega^1,\ldots,\omega^d)$ the projection from $\Omega^\infty$ to 
$\Omega^d:=C_0(\R^+;\mathbb R^d)$.


\begin{assumption}\label{assum: consistent G}
Let \textnormal{$\mathbi{G}:=\{G_d\}_{d\in\mathbb{N}}$} be a family of monotone sublinear functions where $G_d:\mathbb{S}_d\to\mathbb{R}$. We assume that the family \textup{$\mathbi G$} is consistent under
dimension extension: 
\[
G_{d+1}\!\left(
\begin{pmatrix}
A & 0\\
0 & 0
\end{pmatrix}
\right)
= G_d(A), \text{ for every $d\in\mathbb{N}$ and every $A\in\mathbb{S}_d$.}
\]
\end{assumption}
Throughout this subsection, we impose
Assumption~\ref{assum: consistent G}.
In particular, this assumption is satisfied if each \(G_d\)
satisfies Assumption~\ref{assum:G} with the same volatility
bounds in every dimension.

For each \(d\in\mathbb N\), let  \({\mathbb E}^d\) denote the corresponding \(G_d\)-expectation on $\Omega^d$.   Then, by the construction of the $G_d$-Brownian motion,  it is easy to verify that the following consistency condition is satisfied:
\begin{equation}\label{eq:cylinder-consistency}
{\mathbb E}^{d+1}[X]={\mathbb E}^d[X], \text{ for }X\in \Lip\big(\Omega^d\big).
\end{equation}
Here, $\Lip\big(\Omega^d\big)=\bigcup_{T\in\mathbb N} \Lip\big(\Omega^{d}_T\big)$, where
$$\Lip\big(\Omega^{d}_T\big):=\left\{\varphi\left(\pi_d(B_{t_{1}}, B_{t_{2}}, \cdots, B_{t_{n}})\right): n \in \mathbb{N}, 0 \leq t_{1}<\cdots<t_{n} \leq T, \varphi \in C_\text{b,Lip}(\mathbb{R}^{n \times d})\right\}.$$ 
Then, we can define  sublinear expectations $\mathbb E^\infty$ and $\mathbb E^\infty_t$ for $t\ge 0$ on the space $\mathcal H_\infty:=\bigcup_{d\in \mathbb N} \Lip(\Omega^d)$ of cylinder
functions over $\Omega^\infty$
by 
$$\mathbb E^\infty[X]:=\mathbb E^d[X] \text{ and } \mathbb E_t^\infty[X]:=\mathbb E_t^d[X], \text{ if } X\in\operatorname{Lip}(\Omega^d).$$  The consistency condition~\eqref{eq:cylinder-consistency} ensures that
this definition is independent of the finite-dimensional
representation of the cylinder function. Denote by $L_{\mathbi G}^1(\Omega^\infty)$ the completion of
$\mathcal H_\infty$ under the norm
$\mathbb E^\infty[|X|]$.
Similarly, for each $t\geq0$, we define
$L_{\mathbi G}^1(\Omega_t^\infty)$ as the completion of
$\bigcup_{d\in\mathbb N}\operatorname{Lip}(\Omega_t^d)$
under the same norm. We have thus constructed a \(\mathbi G\)-expectation space
$\left(
\Omega^\infty,\,
L_{\mathbi G}^1(\Omega^\infty),\,
(L_{\mathbi G}^1(\Omega_t^\infty))_{t\ge0},\,
(\mathbb E_t^\infty)_{t\ge0}
\right)$
on the infinite-dimensional canonical path space
\(\Omega^\infty\), carrying the infinite-dimensional
\(\mathbi G\)-Brownian motion \(B=(B^1,B^2,\ldots)\).
For every \(d\in\mathbb N\), the coordinate projection
\(\pi_d(B)=(B^1,\ldots,B^d)\) is a \(d\)-dimensional
\(G_d\)-Brownian motion on this space.

In order to characterize properties that hold quasi-surely as stated in Section~\ref{Main results on high-dimensional spectral limits},  we now investigate a Kolmogorov extension theorem for capacities. We aim to  define a capacity on $\Omega^\infty$ which is consistent with the infinite-dimensional $\mathbi G$-Brownian motion.


For each \(d\in\mathbb N\),  choose a nonempty weakly compact family \(\mathscr P_d\subset\mathbf P(\Omega^d)\) representing \({\mathbb E}^d\). 
Define, for bounded Borel function \(X\) on $\Omega^d$ and \(A\in\mathcal B(\Omega^d)\),
\begin{equation*}\label{eq:Ed-cd}
\mathcal E_d(X):=\sup_{P\in\mathscr P_d} E_P[X],
\quad
c_d(A):=\sup_{P\in\mathscr P_d}P(A).
\end{equation*}
By~\eqref{e:equal E} and~\eqref{eq:cylinder-consistency}, we have 
\begin{equation}\label{e:E-d+1}
\mathcal E_{d+1}(X)=\mathcal E_d(X) \text{ for every 
$X\in\operatorname{Lip}(\Omega^d_T)$}.
\end{equation}
Moreover, we have the following result.

\begin{lemma}\label{lem:Cb-consistency}
For every \(X\in C_b(\Omega^d)\),
$\mathcal E_{d+1}(X)=\mathcal E_d(X)$.
\end{lemma}
\begin{proof}
Let \(\mathcal{A}:=\mathscr P_d\cup \tilde{\mathscr P_d}\), where \(\tilde{\mathscr P_d}:=\mathscr P_{d+1}|_{\Omega^d}\) denotes the family of marginal laws on \(\Omega^d\) induced by the measures in \(\mathscr P_{d+1}\). Then \(\mathcal{A}\) is weakly compact and hence tight. Consequently, for any given $\delta\in(0,1)$, there exists a compact set \(K\subset\Omega^d\) such that $\sup_{Q\in\mathcal{A}}Q(K^c)<\delta.$

Consider 
$\Lip(K):=\{Y|_K:Y\in\Lip(\Omega^d)\}\subset C(K)$.
Since $\Lip(K)$ contains the constant functions and separates points of $K$, the Stone--Weierstrass theorem implies that it is dense in $C(K)$.
Thus, given $X\in C_b(\Omega^d)$, for any $\varepsilon>0$, there exists $Y_0\in\Lip\left(\Omega^d\right)$ such that
$\sup_{\omega\in K}|X(\omega)-Y_0(\omega)|<\varepsilon$. Set $Y:=(-M)\vee(Y_0\wedge M)$, where $M:=\|X\|_{\infty}$. Then, $Y\in\Lip\left(\Omega^d\right)$ with
$\|Y\|_{\infty}\leq M$, and 
$$\sup_{\omega\in K}|X(\omega)-Y(\omega)|
\leq \sup_{\omega\in K}|X(\omega)-Y_0(\omega)|<\varepsilon.$$
For any $Q\in\mathcal A$, we now have
\begin{align*}
\int_{\Omega^d}|X-Y|\,\mathrm dQ
=
\int_K|X-Y|\,\mathrm dQ+\int_{K^c}|X-Y|\,\mathrm dQ\leq \varepsilon+2M\delta.
\end{align*} 
Consequently, 
 \begin{align*}
&\Big|\sup_{Q\in\mathscr P_d}\int X\,\mathrm dQ-\sup_{Q\in\tilde{\mathscr P_d}}\int X\,\mathrm dQ\Big|\\
\leq &\Big|
\sup_{Q\in \mathscr P_d}\int X\,\mathrm dQ-
\sup_{Q\in \mathscr P_d}\int Y\,\mathrm dQ
\Big|+\Big|
\sup_{Q\in \mathscr P_d}\int Y\,\mathrm dQ-
\sup_{Q\in \tilde{\mathscr P_d}}\int Y\,\mathrm dQ
\Big|\\
&+\Big|
\sup_{Q\in \tilde{\mathscr P_d}}\int Y\,\mathrm dQ-
\sup_{Q\in \tilde{\mathscr P_d}}\int X\,\mathrm dQ
\Big|\\
\leq &2\varepsilon+4M\delta,
\end{align*}
where the term $\big|
\sup_{Q\in \mathscr P_d}\int Y\,\mathrm dQ-
\sup_{Q\in \tilde{\mathscr P_d}}\int Y\,\mathrm dQ
\big|$ vanishes by ~\eqref{e:E-d+1}. This implies the desired result, noting that  $\varepsilon>0$ and $\delta>0$ can be chosen arbitrarily small.
\end{proof}


\begin{proposition}\label{prop:projectivity}
For every \(d\in\mathbb N\), let $\overline{\operatorname{co}}^{\,w}(\mathscr P_{d+1})|_{\Omega^d}$ be the set of marginal laws on \(\Omega^d\) induced by the measures in \(\overline{\operatorname{co}}^{\,w}(\mathscr P_{d+1})\), where $\overline{\operatorname{co}}^{\,w}(\mathscr P_{d})$ denotes the weakly closed convex hull of $\mathscr P_{d}$. Then, we have
\begin{equation}\label{eq:family-consistency}
\overline{\operatorname{co}}^{\,w}(\mathscr P_{d+1})|_{\Omega^d}=\overline{\operatorname{co}}^{\,w}(\mathscr P_d).
\end{equation}
\end{proposition}
\begin{proof}
We first establish the following characterization:
\begin{equation}\label{eq:dominated-family}
\overline{\operatorname{co}}^{\,w}(\mathscr P_d)
=
\Big\{P\in\mathbf P(\Omega^d):
\int X\,\dd P\le\mathcal E_d(X)
\text{ for every }X\in C_b(\Omega^d)\Big\}.
\end{equation}
The inclusion from left to right is immediate. Conversely, if 
$P_0\notin \overline{\operatorname{co}}^{\,w}(\mathscr P_d)$, then by  the Hahn–Banach separation theorem  there exists $X_0\in C_b(\Omega^d)$ such that $\int X_0\,\dd P_0>\sup_{P\in\overline{\operatorname{co}}^{\,w}(\mathscr P_{d})}\int X_0\,\dd P=\mathcal E_d(X_0)$, and hence $P_0$ does not belong to the right-hand side. This proves~\eqref{eq:dominated-family}.

Denote \(\mathcal Q_d:=\overline{\operatorname{co}}^{\,w}(\mathscr P_{d+1})|_{\Omega^d}\) . 
By Lemma~\ref{lem:Cb-consistency}, 
$\sup_{Q\in \mathcal Q_d}\int X\,\dd Q
=\mathcal E_{d+1}(X)
=\mathcal E_d(X)$, for every $ X\in C_b(\Omega^d)$.
Hence, by ~\eqref{eq:dominated-family},  $\mathcal Q_d\subseteq\overline{\operatorname{co}}^{\,w}(\mathscr P_d)$.
 Since \(\mathcal Q_d\) is weakly closed  and convex, the same separation argument as above yields the reverse inclusion.
\end{proof}

\begin{theorem}\label{thm:extension}
Define
 $$\tilde{\mathscr P}_\infty
:=\left\{P\in\mathbf P(\Omega^{\infty}):
P|_{\Omega^d}\in\overline{\operatorname{co}}^{\,w}(\mathscr P_{d})\text{ for all }d\in\mathbb N\right\},$$ 
where $P|_{\Omega^d}$ denotes the
marginal law on \(\Omega^d\) induced by $P\in\mathbf P(\Omega^{\infty})$. 
Then \(\tilde{\mathscr P}_\infty\) is nonempty and weakly compact, and
\begin{equation}\label{eq:all-marginal-laws}
\{P|_{\Omega^d}:P\in\tilde{\mathscr P}_\infty\}=\overline{\operatorname{co}}^{\,w}(\mathscr P_{d})
\text{ for every } d\in\mathbb N.
\end{equation}

Define, for any bounded Borel function $X$ on $\Omega^\infty$ and $A\in\mathcal B(\Omega^\infty)$, 
\begin{equation*}\label{eq:E-c-infinity}
\mathcal E_\infty(X):=\sup_{P\in\tilde{\mathscr P}_\infty}\int X\,\dd P,
\quad 
c_\infty(A):=\mathcal E_\infty(\mathbf 1_A)=\sup_{P\in\tilde{\mathscr P}_\infty}P(A).
\end{equation*}
Then \(c_\infty\) is a Choquet capacity. Moreover, for every
\(d\in\mathbb N\),  we have
\begin{align}
\mathcal E_\infty(X)
&=\mathcal E_d(X),\,
c_\infty(A)=c_d(A)
\label{eq:all-Borel}
\end{align}
for any bounded measurable function $X$ on $\Omega^d$ and  \(A\in\mathcal B(\Omega^d)\).
\end{theorem}

\begin{proof}
Fix \(d\in\mathbb N\) and \(Q_d\in\overline{\operatorname{co}}^{\,w}(\mathscr P_{d})\). By ~\eqref{eq:family-consistency} we can extend $Q_d$ successively to higher dimensions and take its marginals in lower dimensions. This gives a consistent family $\{Q_d\}$ with $Q_d\in\overline{\operatorname{co}}^{\,w}(\mathscr P_{d})$ for all $d$. The Kolmogorov extension theorem yields  \(Q\in\tilde{\mathscr P}_\infty\) with $Q|_{\Omega^d}=Q_d$, which proves nonemptiness and ~\eqref{eq:all-marginal-laws}.

Let $\mathcal R_i$ be the set of the
$i$-th coordinate marginals of $\overline{\operatorname{co}}^{\,w}(\mathscr P_{i})$.
Then, each $\mathcal R_i$ is weakly compact. Given any \(\varepsilon>0\), we can choose compact sets \(K_i\subset C_0(\R^+; \R)\) such that
$\sup_{\nu\in\mathcal R_i}\nu(K_i^c)<\varepsilon 2^{-i}$.
Then \(K:=\prod_{i=1}^{\infty}K_i\) is compact and $\sup_{P\in\tilde{\mathscr P}_\infty}P(K^c)\le\sum_{i=1}^{\infty}\sup_{\nu\in\mathcal R_i}\nu(K_i^c)
<\varepsilon$.
Thus, \(\tilde{\mathscr P}_\infty\) is tight. 
 It is also weakly closed, noting that each $\overline{\operatorname{co}}^{\,w}(\mathscr P_{d})$ is weakly closed. Thus, \(\tilde{\mathscr P}_\infty\) is weakly compact by Prokhorov’s theorem. Consequently, the upper probability $c_\infty$ is a Choquet capacity. 


Finally, by ~\eqref{eq:all-marginal-laws} and Lemma \ref{lem:barycenters},  for every bounded Borel function $X$ on $\Omega^d$, 
$$\sup_{P\in\tilde{\mathscr P}_\infty}\int X\,\dd P=
\sup_{P\in\overline{\operatorname{co}}^{\,w}(\mathscr P_{d})}\int X\,\dd P=\sup_{P\in\mathscr P_{d}}\int X\,\dd P.$$
This proves ~\eqref{eq:all-Borel} and the proof is completed. 
\end{proof}
Let $\{A_d\}_{d\ge1}$ be a sequence of Borel sets with $A_d\in\mathcal B(\Omega^d)$  such that $c_d(A_d)=0$ for each $d\in\mathbb N$. Then, we have $c_\infty\left(\bigcup_{d=1}^{\infty} A_d\right)
\le\sum_{d=1}^{\infty}c_d(A_d)=0$.

\begin{lemma}\label{lem:barycenters}
Let \(S\) be a Polish space and let \(\mathscr P\subset\mathbf P (S)\) be a nonempty and weakly compact set. Let
 $\overline{\operatorname{co}}^{\,w}(\mathscr P)$ be the weakly closed convex hull of $\mathscr P$.
 Then \(\overline{\operatorname{co}}^{\,w}(\mathscr P)\) is weakly compact and for every bounded Borel function \(f\),
\begin{equation}\label{eq:convexification}
    \sup_{P\in\overline{\operatorname{co}}^{\,w}(\mathscr P)}\int_S f\,\dd P
=\sup_{P\in\mathscr P}\int_S f\,\dd P.\end{equation}
In particular, the two families induce the same upper probabilities on every Borel set.
\end{lemma}
\begin{proof} Clearly \(\overline{\operatorname{co}}^{\,w}(\mathscr P)\) is weakly compact, since $\mathscr P$ is weakly compact. Now we prove \eqref{eq:convexification}. By  \cite[Theorem 461I]{Fremlin2013}, we have that
$\overline{\operatorname{co}}^{\,w}(\mathscr{P})=\{\bar \nu:\nu\in\mathbf P(\mathscr{P})\}$,
where $\bar \nu\in \mathbf P(S)$ is the barycenter of $\nu$ (see \cite[Definitions 461A]{Fremlin2013}) which is uniquely characterized by 
\begin{equation}
\label{eq:barycenter}
\int_S g\,d\bar \nu
=
\int_{\mathscr{P}}
    \left(\int_S g\,dP\right)\nu(dP),
\text{ for all } g\in C_b(S).
\end{equation}
It is easy to verify that $\bar \nu(A)=\int_{\mathscr{P}}P(A)\,\nu(dP),
\,A\in\mathcal{B}(S)$.
Thus, eq.~\eqref{eq:barycenter} actually holds for all bounded measurable functions. Consequently, for every bounded Borel measurable
$f:S\to\mathbb{R}$,
\begin{equation*}
\int_S f d\bar \nu=
\int_{\mathscr{P}}
    \left(\int_S f\,dP\right)\nu(dP)\le\sup_{P\in\mathscr{P}}\int_S f\,\dd P.
\end{equation*}
The reverse inequality follows from
$\mathscr{P}\subseteq\overline{\operatorname{co}}^{\,w}(\mathscr{P})$.
This completes the proof.
\end{proof}

\section{A preliminary study on \texorpdfstring{$G$}{}-Dyson Brownian motion}\label{sec:pre-G-Dyson}

In this section, we provide a preliminary analysis on  $G$-Dyson Brownian motion, i.e., the eigenvalue processes $\lambda^N=(\lambda_1^N, \dots, \lambda_N^N)$ of the symmetric $G$-Brownian matrix $X^N$ given by~\eqref{eq:BN}.

For each fixed  $t$ and $i$,  $\lambda _{i}^{N} \left ( t \right )$ is a function of $d$ variables $\left (X_{kh}^{N}(t), 1\le k\le h\le N\right )$, which is denoted by $\Phi_i(X^N(t))$. Then, we have (see, e.g.,   \cite[Appendix]{song2023eigenvalues})
\begin{equation}\label{e:derivative-matrix}
\begin{cases} 
 & \frac{\partial \Phi_i(X^N)}{\partial X^{N}_{k h}}  =2 U_{k i} U_{h i} \mathbf{I}_{[k \neq h]}+U_{k i}^{2} \mathbf{I}_{[k=h]}; \\[10pt]
& \frac{\partial^{2} \Phi_i(X^N)}{\partial\left( X_{k h}^{N}\right)^{2}}  =2 \sum_{j: j \neq i} \frac{\left|U_{k i} U_{h j}+U_{h i} U_{k j}\right|^{2}}{\lambda_{i}^{N}-\lambda_{j}^{N}} \mathbf{I}_{[k \neq h]}+2 \sum_{j: j \neq i} \frac{\left|U_{k i} U_{k j}\right|^{2}}{\lambda_{i}^{N}-\lambda_{j}^{N}} \mathbf{I}_{[k=h]};\\[10pt]
 &\frac{\partial^{2} \Phi_i(X^{N})}{\partial X^{N}_{k h} \partial X^{N}_{k^{\prime} h^{\prime}}} 
 =2 \sum_{j: j \neq i} \frac{\left(U_{k i} U_{h j}+U_{h i} U_{k j}\right)\left(\mathbf{I}_{[k \neq h]}+\mathbf{I}_{[k=h]} / 2\right)\left(U_{k^{\prime} i} U_{h^{\prime} j}+U_{h^{\prime} i} U_{k^{\prime} j}\right)\left(\mathbf{I}_{\left[k^{\prime} \neq h^{\prime}\right]}+\mathbf{I}_{\left[k^{\prime}=h^{\prime}\right]} / 2\right)}{\lambda^{N}_{i}-\lambda^{N}_{j}}.
\end{cases}
\end{equation}

\begin{remark}\label{remark:smooth eigenvalue and eigenvector}
For  a stochastic process $X^{N}$ taking values in  $\mathbb{S}_{N}$, denote by $\Lambda^{N}_{t}=\operatorname{diag}\left[\lambda^{N}_{i}(t)\right]$ the diagonal matrix of eigenvalues of $X^{N}_{t}$, arranged in non-decreasing order. Let $U_t$ be the matrix whose columns are the corresponding orthonormal eigenvectors of $X^{N}_{t}$.  Then, $X^N$ has the orthogonal decomposition $X^N={U}\Lambda^{N} U^\top.$

Assume further that  $X^{N}_{0}$ has $N$ distinct  eigenvalues.  By ~\cite[Chapter~II, Theorem~1.8]{Kato1995}, in a neighbourhood of any real symmetric matrix with simple spectrum, the ordered eigenvalues and the corresponding spectral projections depend analytically, and hence smoothly, on the matrix entries. Moreover, the  orthogonal eigenvector matrix constructed by spectral projections can be chosen locally smoothly in the matrix argument. Therefore,  matrices $ \Lambda^{N}$  and $U$  can be chosen  as smooth functions of  $X^{N}$  until the first collision time
\begin{equation}\label{e:tau}
    \tau=\inf \left\{t: \lambda^{N}_{i}(t)=\lambda^{N}_{j}(t) \text { for some } i \neq j\right\}.\end{equation}

Now we show that the random time $\tau$ is a stopping time in the sense of Definition~\ref{def:stopping time}. Denote $F=\left \{ x\in  \mathbb {R}^{N}:\exists i \neq j,   x_{i}=x_{j} \right \}$ and  $d(x,F):=\inf_{y\in F}|x-y|$.
Then, $\tau=\inf \left\{t>0:   \lambda^{N}(t)\in F\right\}$. Since $X^N(t)$
is $\mathcal B(\Omega_t)$-measurable,  so are $\lambda^N(t)$
and $d(\lambda^N(t),F)$. By the  Hoffman-Wielandt inequality (see, e.g., \cite{agz10}),  
$\lambda^N$ is continuous with respect to $X^N$, and hence  $\lambda^N(t)$ is a  continuous function of $t$. Thus, for each $t>0$, $$ \{\tau\le t\}=
\bigcap\limits_{n=1}^{\infty}
\bigcup\limits_{\substack{s\in\mathbb{Q}\\0\leq s\leq t}}
\Big\{
d\bigl(\lambda^N(s),F\bigr)<\frac{1}{n}
\Big\}\in \mathcal{B}(\Omega_{t}),$$
and hence $\tau$ is a stopping time.
\end{remark}



The following result provides  the SDEs  satisfied by $G$-Dyson Brownian motion: 
\begin{theorem}\label{thm:general-DBm}
 Let $X^N$  be given by~(\ref{eq:BN}). Assume that the initial eigenvalues of $X^N$ are distinct and 
strictly ordered as  $\lambda_1^N(0)<\lambda_2^N(0)<\cdots<\lambda_N^N(0)$. We further assume that the
underlying G-Brownian motion satisfies Assumption~\ref{assum:G'} with $d=\frac{N\left ( N+ 1\right ) }{2}$. For $R>0$, define
\begin{equation}
\tau_R:=\inf\Big\{t>0:\min_{i\neq j}\left|\lambda_i^N(t)-\lambda_j^N(t)\right|<R^{-1}\Big\}.
\end{equation}  
Then,  for $t<\tau_R$, 
\begin{equation}\label{eq:dyson}
\begin{aligned}
\mathrm{d}\lambda _{i}^{N}( t ) &=\frac{1}{N}  \sum_{j: j \neq i} \frac{1}{\lambda_{i}^N-\lambda_{j}^N} \Big( \sum_{1\le k< h\le N} \left|U_{k i} U_{h j}+U_{h i} U_{k j}\right|^{2}\mathrm{d}\left \langle B_{kh} \right \rangle _{t} +\sum_{k=1}^{N} 2\left|U_{k i} U_{k j}\right|^{2} \mathrm{d}\left \langle B_{kk}  \right \rangle _{t}\Big)\\
&\qquad +\frac{\sqrt{2} }{\sqrt{N} }\Big(\sum_{1\le k<  h\le N}\sqrt{2}U_{ki}U_{hi}\mathrm{d}B_{kh}+\sum_{k=1}^{N}U_{ki}^{2}\mathrm{d}B_{kk}\Big), \quad i=1,\dots,N.
\end{aligned}
\end{equation}
\end{theorem}

\begin{proof}
For $t<\tau_R$,  $\Phi_{i} \in C^{2}\left(\mathbb{R}^{n}\right)$ for $i=1,\dots,N$, where $n=\frac{N(N+1)}{2}$.
Then, by Proposition~\ref{Prop:18} and \eqref{e:derivative-matrix},  we get
 \begin{align*}\mathrm{d}\lambda _{i}^{N} \left ( t \right ) &=\mathrm{d}\Phi _{i} \left ( X^{N} \left ( t \right )  \right )\\
&=\sum_{1\le k<  h\le N}\frac{\partial\Phi _{i} \left ( X^{N}   \right ) }{\partial x_{kh} }\frac{1}{\sqrt{N} }\mathrm{d}B_{kh}+\sum_{k=1}^{N}\frac{\partial\Phi _{i} \left ( X^{N}  \right ) }{\partial x_{kk} }\frac{\sqrt{2} }{\sqrt{N} }\mathrm{d}B_{kk} \\
&\qquad +\frac{1}{2} \Big(   \sum_{1\le k<  h\le N} \frac{\partial^{2} \Phi _{i} \left ( X^{N }  \right ) }{\partial x_{kh}^{2}} \frac{1}{N} \mathrm{d}\left \langle B_{kh}  \right \rangle _{t}+
\sum_{k=1}^{N}\frac{\partial^{2} \Phi _{i} \left ( X^{N }  \right ) }{\partial x_{kk} ^{2}} \frac{2}{N} \mathrm{d}\left \langle B_{kk}  \right \rangle _{t}\Big ) \\
&=\sum_{1\le k<  h\le N}2U_{ki}U_{hi}\frac{1}{\sqrt{N} }\mathrm{d}B_{kh}+\sum_{k=1}^{N}U_{ki}^{2}\frac{\sqrt{2} }{\sqrt{N} }\mathrm{d}B_{kk}\\
& \qquad + \frac{1}{2} \Big(\sum_{1\le k< h\le N}2 \sum_{j: j \neq i} \frac{\left|U_{k i} U_{h j}+U_{h i} U_{k j}\right|^{2}}{\lambda^N_{i}-\lambda^N_{j}} \frac{1}{N} \mathrm{d}\left \langle B_{kh}  \right \rangle _{t}+\sum_{k=1}^{N} 2\sum_{j: j \neq i} \frac{\left|U_{k i} U_{k j}\right|^{2}}{\lambda^N_{i}-\lambda^N_{j}}\frac{2}{N} \mathrm{d}\left \langle B_{kk}  \right \rangle _{t}\Big )  \\
&=\sum_{1\le k<  h\le N}2U_{ki}U_{hi}\frac{1}{\sqrt{N} }\mathrm{d}B_{kh}+\sum_{k=1}^{N}U_{ki}^{2}\frac{\sqrt{2} }{\sqrt{N} }\mathrm{d}B_{kk}\\
&\quad +\frac{1}{N}  \sum_{j: j \neq i} \frac{1}{\lambda^N_{i}-\lambda^N_{j}} \Big ( \sum_{1\le k< h\le N} \left|U_{k i} U_{h j}+U_{h i} U_{k j}\right|^{2}\mathrm{d}\left \langle B_{kh} \right \rangle _{t} +\sum_{k=1}^{N} 2\left|U_{k i} U_{k j}\right|^{2} \mathrm{d}\left \langle B_{kk}  \right \rangle _{t}\Big).
 \end{align*}
\end{proof}



Consider the martingale part on the right-hand side of \eqref{eq:dyson} by,  
\[ \mathrm{d}M_{i}  \left ( t \right ) :=\sum_{1\le k<  h\le N}\sqrt{2}U_{ki}U_{hi}\mathrm{d}B_{kh}+\sum_{k=1}^{N}U_{ki}^{2}\mathrm{d}B_{kk}, \quad i=1, \dots, N. 
\]
We have 
 \begin{align*}
\sum_{1\le k<  h\le N}(\sqrt{2}U_{ki}U_{hi})^2+\sum_{k=1}^{N}U_{ki}^{4}   =\Big(\sum_{k=1}^{N}U_{ki}^{2}\Big)^2=1,
\end{align*}
 and thus Corollary \ref{G property} yields that,  under the diagonal structure assumption \ref{assum:G'}, each component $M_{i}(t)$ is a $\tilde{G}$-Brownian motion, for $i=1,2,\dots,N$.  However,  $M=(M_1,\ldots,M_N)^{\top}$ may fail to be a $G$-Brownian motion for any $G$. By contrast, under Assumption~\ref{assum:G} for $G$, by Corollary~\ref{corollary:G-brownian motion}  $M$ is a $\bar G$-Brownian motion with $\bar G$ still satisfying Assumption~\ref{assum:G}. See Remark~\ref{rem:Levy} for detailed discussions. Moreover, under Assumption~\ref{assum:G}, the $G$-Dyson Brownian motion given by \eqref{eq:dyson} admits a more concise representation~\eqref{e:G-Dyson}, as shown in Theorem~\ref{thm:dyson non-collision}.

\subsection{Orthogonal invariance and the GOE property}\label{sec: GOE discussion} In the classical setting, a Gaussian Orthogonal Ensemble (GOE) matrix $H$ satisfies $O^{\top} H O \overset{d}{=} H$ for any orthogonal matrix $O$. Let $\{\xi_{ij}, 1\le i\le j\le N\}$ be a family of i.i.d.\ standard  Gaussian random variables. Let $P_N$ be the law of the symmetric random matrix 
\begin{equation*}
H = \begin{bmatrix}
\sqrt{2}\xi_{11} & \xi_{12} & \cdots & \xi_{1N} \\
\xi_{12} & \sqrt{2}\xi_{22} & \cdots & \xi_{2N} \\
\vdots & \vdots & \ddots & \vdots \\
\xi_{1N} & \xi_{2N} & \cdots & \sqrt{2}\xi_{NN}
\end{bmatrix}.
\end{equation*}
Then, $O^{\top} H O$ still has the law $P_N$ for any $N \times N$ orthogonal matrix $O$ and thus $H$ is a GOE matrix (see, e.g., \cite{agz10}).

In this subsection, we investigate the GOE property under $G$-expectation.
We first identify a symmetric matrix $H = (H_{ij})_{1 \le i,j \le N}$ with the vector \begin{equation}\label{e:H'}
\widetilde{H}
=
\bigl(H_{11},H_{12},\ldots,H_{1N},H_{22},\ldots,H_{2N},
\ldots,H_{NN}\bigr)^{\top}.
\end{equation}
Let $D$ be a $d\times d$ diagonal matrix, recalling that $d:=\frac{N(N+1)}2$,   whose diagonal entries corresponding to the diagonal entries of $H$ are $1/\sqrt{2}$, and whose remaining diagonal entries are equal to $1$. Thus, we have 
$$D\widetilde H= \bigl(\xi_{11},\xi_{12},\ldots,\xi_{1N},\xi_{22},\ldots,\xi_{2N},
\ldots,\xi_{NN}\bigr)^{\top}:=\xi.$$
For any $N\times N$ orthogonal matrix $O$, there exists a linear map $Q:\mathbb{R}^{d}\longrightarrow\mathbb{R}^{d}$
such that $D\widetilde{O^{\top}HO}=Q\xi,$ where $\widetilde{O^{\top}HO}$ denotes  the vector consisting of  the upper-triangular entries of $O^{\top}HO$, as in~\eqref{e:H'}. We show that the matrix $Q$ is orthogonal.

Given any two  symmetric matrices $H^1$ and $H^2$, we have 
\begin{equation*}
\begin{aligned}
\big\langle D\widetilde{H}^1,D\widetilde{H}^2\big\rangle
&=
\frac{1}{2}\sum_{i=1}^N H^1_{ii}H^2_{ii}
+
\sum_{1\leq i<j\leq N}H^1_{ij}H^2_{ij} =
\frac{1}{2}\sum_{i=1}^N\sum_{j=1}^N H^1_{ij}H^2_{ij} =
\frac{1}{2}\operatorname{tr}(H^1H^2).
\end{aligned}
\end{equation*}
 Using the above identity, we obtain
\begin{align*}
\|Q\xi\|^2
&=
\big\langle D\widetilde{O^{\top}HO},
D\widetilde{O^{\top}HO}\big\rangle 
=
\frac{1}{2}\operatorname{tr}\bigl((O^{\top}HO)^{2}\bigr) \notag=
\frac{1}{2}\operatorname{tr}(H^2) =
\big\langle D\widetilde{H},D\widetilde{H}\big\rangle
=
\|\xi\|^2.
\end{align*}
This implies that $Q$ is an orthogonal matrix.

Now let $\xi_{ij}=B_{ij}(1)$, 
where $B=(B_{11},B_{12},\ldots,B_{1N},B_{22},\ldots,B_{2N},
\ldots,B_{NN})^T$ is a $d$-dimensional $G$-Brownian motion. By the construction of $D$, we have
$D\widetilde{H}=B(1)$. Our goal is to determine whether $O^{\top}HO$ and $H$ are identically distributed under the $G$-expectation for any orthogonal matrix $O$. That is, the question is whether $\widetilde{O^{\top}HO}\overset{d}{=}\widetilde{H},$ which is equivalent to $D\widetilde{O^{\top}HO}
\overset{d}{=}
D\widetilde{H}$, or, equivalently, to $QB(1)\overset{d}{=}B(1).$   Thus, the orthogonal invariance of the matrix $H$ reduces to the invariance in distribution of the underlying $G$-Brownian motion under  orthogonal transformations $Q$.

By Proposition~\ref{prop: transform-G} we know that $QB$ is a $G_Q$-Brownian motion, where $G_Q(A)=G(Q^\top A Q)$.  Assuming that $G$ satisfies Assumption~\ref{assum:G},  it is easy to verify that $G_Q=G$. Thus, $QB\overset{d}{=}B$ and hence the corresponding matrix $H$ is a GOE matrix under $G$-expectation. 

However, the orthogonal invariance may fail under  Assumption~\ref{assum:G'} if $\underline{\sigma}<\overline{\sigma}$, as indicated by the following example. 
\begin{example}
When $N=2$, for an orthogonal $2\times2$ matrix $O
=
\frac{1}{\sqrt{2}}
\begin{pmatrix}
1 & -1\\
1 & 1
\end{pmatrix}$,  the corresponding $3\times 3$ orthogonal matrix $Q$ is given by  $$Q=\begin{pmatrix}
\frac{1}{2} & \frac{1}{\sqrt 2}&\frac{1}{2}\\
-\frac{1}{\sqrt 2} & 0&\frac{1}{\sqrt 2}\\
\frac{1}{2} &- \frac{1}{\sqrt 2}&\frac{1}{2}
\end{pmatrix}.$$
Under Assumption~\ref{assum:G'} with $\underline{\sigma}<\overline{\sigma}$,  it can be verified that $$G_Q(A)=G(Q^{\top}AQ)=0\neq \frac{1}{2}
\left(
\overline{\sigma}^{\,2}
-
\underline{\sigma}^{\,2}
\right)=G(A),$$ 
where \begin{equation*}
A
=
\begin{pmatrix}
1 & 0&0\\
0 & 0&0\\
0&0&-1
\end{pmatrix}, \text{ and } Q^\top AQ
=
\begin{pmatrix}
0 & \frac{1}{\sqrt 2} &0\\
\frac{1}{\sqrt 2} & 0&\frac{1}{\sqrt 2}\\
0&\frac{1}{\sqrt 2}&0
\end{pmatrix}.
\end{equation*}
This indicates that $QB$ and $B$ have different distributions under Assumption~\ref{assum:G'}.
\end{example}

\subsection{A colliding \texorpdfstring{$G$}{}-Dyson Brownian motion: an example}\label{sec:collion} In this subsection, we provide an example demonstrating  that the non-collision property may fail for $G$-Dyson Brownian motion.



When $\Omega=C_0^d(\R^+)$, the $G$-expectation $\E$ coincides with the upper expectation $\hat \E[\cdot]=\sup_{P\in\mathscr P}E_P[\cdot]$ on $L_G^1(\Omega)$, where $\mathscr P=\overline{\mathcal P}$ with $\mathcal P=\{P_\theta: \theta \in \mathcal A_{0,\infty}^\Theta\}$. Here, $P_\theta$ denotes the law of $\int_0^\cdot \theta \dd W_s$, where $W$ is a $d$-dimensional Wiener process on $\Omega$, and $A_{0,\infty}^\Theta$ denotes the space of $\Theta$-valued adapted processes, with $\Theta$ being a given bounded closed subset of $\R^{d\times d}$. We refer to Section~\ref{sec:G-exp} for more details.

Consider the case $N=2$ (i.e., $d=3$). Assume that the generator $G$ of the underlying  $G$-Brownian motion is given by~\eqref{e:G-Theta}, with 
\begin{equation}\label{e:Theta}
\Theta=\Big\{\theta=(\theta_{ij})_{3\times 3}\in \R^{3\times 3}:0<\underline{\sigma}^{2} \leq \theta_{i 1}^{2}+\theta_{i 2}^{2}+\theta_{i 3}^{2} \leq \bar{\sigma}^{2}, i=1,2,3 \Big\}.
\end{equation}
Assuming  $\lambda^{2}_{1}(0)<\lambda^{2}_{2}(0)$,  we show below that there exists $\theta^* \in \mathcal{A}_{0, \infty}^{\Theta}$, such that  $\lambda_2^2(t)-\lambda_1^2(t)$ hits zero in finite time  $P_{\theta^*}$-almost surely. This indicates that the eigenvalues collide in finite time $P_{\theta^*}$-almost surely. Consequently, quasi-sure noncollision fails.

Let $$X_t = \begin{pmatrix} X_t^{11} & X_t^{12} \\ X_t^{12} & X_t^{22} \end{pmatrix}=\begin{pmatrix} \sqrt 2B_t^{11}+X_0^{11} & B_t^{12}+X_0^{12} \\ B_t^{12}+X_0^{12} & \sqrt 2B_t^{22}+X_0^{22} \end{pmatrix}$$ be a symmetric matrix-valued  process with distinct initial eigenvalues, where $B=(B^{11}, B^{12}, B^{22})$ is a $3$-dimensional $G$-Brownian motion. Define 
\begin{equation*}
Z_t := X_t^{11} - X_t^{22}, \quad Y_t := 2X_t^{12},\quad R_t := \sqrt{Z_t^2 + Y_t^2},\quad \tau^* := \inf\{t \ge 0 : R_t = 0\}.
\end{equation*}
Then, two eigenvalues of $X_t$ collide if and only if 
$R_t = 0$, and $\tau^*$ is the first spectral collision time. 


Given the bounds $0 < {\underline{\sigma}} < {\overline{\sigma}} < \infty$, we introduce 
\begin{equation*} \alpha= 4{\overline{\sigma}}^2, \quad \beta = 4{\underline{\sigma}}^2.
\end{equation*} We now construct the  process $\theta^*(t) = (\theta_1^*(t), \theta_2^*(t), \theta_3^*(t))^\top$:
\begin{equation}\label{eq:control_def}
\begin{cases}
\theta_1^*(t) &:= \left( \frac{\sqrt{\alpha}}{2\sqrt 2}\frac{Z_t}{R_t}, \, -\frac{\sqrt{\beta}}{2\sqrt 2}\frac{Y_t}{R_t}, \, \frac{\sqrt{\beta}}{2\sqrt 2} \right)^\top, \\
\theta_2^*(t) &:= \left( \frac{\sqrt{\alpha}}{2}\frac{Y_t}{R_t}, \, \frac{\sqrt{\beta}}{2}\frac{Z_t}{R_t}, \, 0 \right)^\top, \\
\theta_3^*(t) &:= \left( -\frac{\sqrt{\alpha}}{2\sqrt 2}\frac{Z_t}{R_t}, \, \frac{\sqrt{\beta}}{2\sqrt 2}\frac{Y_t}{R_t}, \, \frac{\sqrt{\beta}}{2\sqrt 2} \right)^\top.
\end{cases}
\end{equation}
We have ${\underline{\sigma}}^2=\frac{\beta}{4}\leq |\theta_2^*(t)|^2 = \frac{\alpha}{4}\frac{Y_t^2}{R_t^2} + \frac{\beta}{4}\frac{Z_t^2}{R_t^2}\leq \frac{\alpha}{4}= {\overline{\sigma}}^2$ 
and ${\underline{\sigma}}^2=\frac{\beta }{4}\leq|\theta_1^*(t)|^2 = |\theta_3^*(t)|^2 = \frac{\alpha}{8}\frac{Z_t^2}{R_t^2} + \frac{\beta}{8}\frac{Y_t^2}{R_t^2} + \frac{\beta}{8}\leq\frac{\alpha + \beta}{8}< {\overline{\sigma}}^2$. Thus, we get that $|\theta_i^*(t)|^2 \in [{\underline{\sigma}}^2, {\overline{\sigma}}^2] ,i = 1, 2, 3$ for all $t\in[0,\tau^*)$.

Then, under $P_{\theta^*}$ we have, 
 \begin{align*}
\mathrm{d}Z_t &= \sqrt2(\theta_1^*(t) - \theta_3^*(t))^\top \mathrm{d}W_{t} = \sqrt{\alpha}\frac{Z_t}{R_t} \mathrm{d}W^{11}_{t} - \sqrt{\beta}\frac{Y_t}{R_t} \mathrm{d}W^{12}_{t}, \\
\mathrm{d}Y_t &= 2{\theta_2^*(t)}^\top \mathrm{d}W_{t} = \sqrt{\alpha}\frac{Y_t}{R_t}\mathrm{d}W^{11}_{t}+ \sqrt{\beta}\frac{Z_t}{R_t} \mathrm{d}W^{12}_{t},
\end{align*}
where   $W=(W^{11}, W^{12}, W^{22})$ is a $3$-dimensional Wiener process. Thus, we get
\begin{equation*}
\begin{aligned}
\mathrm{d}\langle Z \rangle_t = (\alpha\frac{Z_t^2}{R_t^2} + \beta\frac{Y_t^2}{R_t^2})\mathrm{d}t,\quad  
\mathrm{d}\langle Y \rangle_t = (\alpha\frac{Y_t^2}{R_t^2} + \beta\frac{Z_t^2}{R_t^2})\mathrm{d}t, \quad 
\mathrm{d}\langle Z, Y \rangle_t = (\alpha - \beta)\frac{Z_tY_t}{R_t^2} \mathrm{d}t.
\end{aligned}
\end{equation*}

Applying It\^{o}'s formula to $R_t = \sqrt{Z_t^2 + Y_t^2}$, we obtain for $t<\tau^*$,
\begin{align*}
    \mathrm{d}R_t 
    =& \frac{Z_t}{R_t}\mathrm{d}Z_t + \frac{Y_t}{R_t}\mathrm{d}Y_t + \frac{1}{2}\left( \frac{\partial^2 R_t}{\partial Z_t^2}\mathrm{d}\langle Z \rangle_t + 2\frac{\partial^2 R_t}{\partial Z_t \partial Y_t}\mathrm{d}\langle Z,Y \rangle_t + \frac{\partial^2 R_t}{\partial Y_t^2}\mathrm{d}\langle Y \rangle_t \right)\\
    =&\frac{Z_t}{R_t} \mathrm{d}Z_t + \frac{Y_t}{R_t} \mathrm{d}Y_t+ \frac{1}{2}\left(\frac{Y_t^2}{R_t^3}\mathrm{d}\langle Z \rangle_t -2\frac{Z_tY_t}{R_t^3}\mathrm{d}\langle Z,Y \rangle_t + \frac{Z_t^2}{R_t^3}\mathrm{d}\langle Y \rangle_t \right)\\
    =& \frac{Z_t}{R_t} \left(\sqrt{\alpha}\frac{Z_t}{R_t} \mathrm{d}W^{11}_{t} - \sqrt{\beta}\frac{Y_t}{R_t} \mathrm{d}W^{12}_{t}\right)  + \frac{Y_t}{R_t} \left(\sqrt{\alpha}\frac{Y_t}{R_t} \mathrm{d}W^{11}_{t}+ \sqrt{\beta}\frac{Z_t}{R_t} \mathrm{d}W^{12}_{t}\right)\\
    &+ \frac{1}{2R_t} \Big[ \frac{Y_t^2}{R_t^2}(\alpha\frac{Z_t^2}{R_t^2} + \beta\frac{Y_t^2}{R_t^2}) + \frac{Z_t^2}{R_t^2}(\alpha\frac{Y_t^2}{R_t^2} + \beta\frac{Z_t^2}{R_t^2}) - 2\frac{Y_t^2}{R_t^2}\frac{Z_t^2}{R_t^2}(\alpha - \beta) \Big] \mathrm{d}t \\
    =&\sqrt{\alpha}\mathrm{d}W^{11}_{t}+\frac{\beta}{2R_t} \mathrm{d}t.
    \end{align*}

Define the scaled process $\rho_s := R_{\frac{s}{\alpha}}$.
By the Brownian scaling property, $\widetilde{W}_s := \sqrt{\alpha} W^{11}_{\frac{s}{\alpha}}$ is a standard 1-dimensional Brownian motion. Then we have
\begin{equation*}
\mathrm{d}\rho_s = \mathrm{d}R_{\frac{s}{\alpha}} = \sqrt{\alpha} \mathrm{d}W^{11}_{\frac{s}{\alpha}} + \frac{\beta}{2R_{\frac{s}{\alpha}}} \mathrm{d}\left(\frac{s}{\alpha}\right) = \mathrm{d}\widetilde{W}_s + \frac{\beta/\alpha}{2\rho_s} \mathrm{d}s.
\end{equation*}
Thus, $\rho_s$ is a $(1 + \frac{\beta}{\alpha})$-dimensional  Bessel process. Since $\beta<\alpha$, the dimension $1 + \frac{\beta}{\alpha} < 2$ and hence $P_{\theta^*}(\tau^* < \infty) = 1$.   Therefore, the quasi-sure non-collision property fails.


\subsection{A concluding remark}\label{rem:choice1}
The analysis in this section shows that
the GOE property may fail under Assumption~\ref{assum:G'}. Moreover, if the generator $G$ loses the diagonal structure, as in the example constructed in Section~\ref{sec:collion}, then q.s.\ noncollision property of the corresponding $G$-Dyson Brownian motion fails. These observations motivate our choice of the diagonal structure in Assumption~\ref{assum:G} for the underlying $G$-Brownian motion when investigating the spectral behavior of the generalized Wishart process.

\section{Spectral non-collision property for generalized \texorpdfstring{$G$}{}-Wishart processes}\label{sec:non-collision}

In this section, we prove that under Assumption~\ref{assum:G}, the eigenvalue processes of the generalized $G$-Wishart process~\eqref{e:g-Wishart} do not collide q.s. and the $G$-SDE~\eqref{eq:eigenvalue-dynamics} admits a unique strong solution. 
For $G$-Dyson Brownian motion, we obtain the comparison principle and further
 establish the existence and uniqueness of a strong solution to the $G$-SDE~\eqref{e:G-Dyson} satisfied by the eigenvalue processes even when the initial eigenvalues are allowed to coincide. For $G$-Wishart process,
we derive the $G$-SDE~\eqref{G-N-wishart SDE} satisfied by an $N\times N$ $G$-Wishart process and establish the existence and uniqueness of a positive-definite solution to the $G$-SDE~\eqref{e:G-Wishart} satisfied by $G$-Wishart process with  parameter  $\alpha$ when $\alpha\ge N+1$.

We will  work on the extended $G$-expectation space $(\Omega, L_{*}^{1}(\Omega),(L_{*}^{1}(\Omega_{t}))_{t \geq 0}
)$, throughout the rest of the paper.  This allows us to employ  It\^o's formula for $C^2$ functions. See Section~\ref{sec:extended-G} for details.


\subsection{Non-collision for generalized Wishart process}\label{subsection:Generalized Wishart}

In this subsection, we study the eigenvalue process~\eqref{eq:eigenvalue-dynamics} of the generalized Wishart process $X^{N}$ given by~\eqref{e:g-Wishart}.

Recall that the functions $g_{N}, h_{N}, b_{N}: \mathbb{R} \rightarrow \mathbb{R}$ act on the spectrum of $X^{N}_{t}$. That is, $f(X^{N}) = {U} f(\Lambda^{N}) U^{\top}$ for $f=g_{N}, h_{N}$ and $b_{N}$,
where $X^{N}={U}\Lambda^{N} U^\top$, with $\Lambda^N$ and $U$  defined as  in Remark~\ref{remark:smooth eigenvalue and eigenvector}.  It is clear that we have 
\begin{equation}\label{e:UfU}
\sum_{1\le k,h\le N}U_{ki}f(X^N)_{kh}U_{hj} =f(\lambda_i^N)\mathbf I_{\{i=j\}} ~ \text{ and } ~ (f(X^N))^n= f^n(X^N)  \text{ for any } n\in\mathbb N.
\end{equation}


\begin{proof}[Proof of Theorem~\ref{thm:eigenvalue-dynamics}]\makeatletter\def\@currentlabelname{The proof}\makeatother\label{proof:thm:eigen-dyn}
Applying It\^o's formula (Proposition~\ref{Prop:18}) to $\Phi_i(X^N(t))$ together with~ \eqref{e:derivative-matrix}, we get for $t<\tau$, 
 \begin{align*}
 &\mathrm{d}\lambda^{N} _{i} \left ( t \right )  =\mathrm{d}\Phi _{i} \left ( X^{N} \left ( t \right )  \right )\\
&=\sum_{1\le k\le h\le N}\frac{\partial\Phi _{i} \left ( X^{N}  \right ) }{\partial x_{kh} }\mathrm{d}X^{N}_{kh} +\frac{1}{2}  \sum_{1\le k\le h\le N,1\le k^{\prime}\le h^{\prime}\le N} \frac{\partial^{2} \Phi _{i} \left ( X^{N}  \right ) }{\partial x_{kh}\partial x_{k^{\prime}h^{\prime}}}  \mathrm{d}\left \langle X^{N}_{kh},X^{N}_{k^{\prime}h^{\prime}}  \right \rangle _{t} \\
&=\sum_{1\le k\le  h\le N}(2 U_{k i} U_{h i} \mathbf{I}_{[k \neq h]}+U_{k i}^{2} \mathbf{I}_{[k=h]})\Big(\sum_{p, q} g_{N}^{k p} \mathrm{d} B_{t}^{p q} h_{N}^{q h}+\sum_{p, q} h_{N}^{k p} \mathrm{d} B_{t}^{q p} g_{N}^{q h}+b_{N}^{kh} \mathrm{d}\langle B^{hh} \rangle_t\Big)\\
  &\qquad + \frac{1}{2}  \sum_{1\le k\le h\le N,1\le k^{\prime}\le h^{\prime}\le N} \frac{\partial^{2} \Phi _{i} \left ( X^{N}  \right ) }{\partial x_{kh}\partial x_{k^{\prime}h^{\prime}}}  \mathrm{d}\left \langle X^{N}_{kh},X^{N}_{k^{\prime}h^{\prime}}  \right \rangle _{t}\\
&=\mathrm{d} M_{i}+\sum_{k,h} U_{k i} b_{N}^{kh} U_{hi}\mathrm{d}\langle B^{hh} \rangle_t
+\sum_{j: j \neq i} \frac{1}{\lambda^{N}_{i}-\lambda^{N}_{j}}\sum_{ k,k',h, h'}U_{k i} U_{h j}U_{k^{\prime} i} U_{h^{\prime} j}\mathrm{d}\left \langle X^{N}_{kh},X^{N}_{k^{\prime}h^{\prime}}  \right \rangle_{t},
\end{align*}
where
\begin{equation}\label{e:dMi}
\mathrm{d} M_{i}=\sum_{k,h} U_{k i} \Big(\sum_{p, q} g_{N}^{kp} \mathrm{d} B_{t}^{pq} h_{N}^{qh}+\sum_{p, q} h_{N}^{kp} \mathrm{d} B_{t}^{qp} g_{N}^{qh}\Big) U_{hi}.
\end{equation}

Now we calculate the quadratic co-variation of $X^N$. Noting that 
$$\mathrm{d}X^{N}_{kh}=\sum_{p, q} g_{N}^{k p} \mathrm{d} B_{t}^{p q} h_{N}^{q h}+\sum_{p, q} h_{N}^{k p} \mathrm{d} B_{t}^{q p} g_{N}^{q h}+b_{N}^{kh} \mathrm{d}\langle B_{hh} \rangle_t,$$
and similarly for $\dd X^{N}_{k'h'}$, we get 
\begin{align*}
&\mathrm{d}\left \langle X^{N}_{kh},X^{N}_{k^{\prime}h^{\prime}}  \right \rangle _{t}\\
= & \sum_{p, q}\left(g_{N}^{k p} h_{N}^{q h} g_{N}^{k^{\prime}p} h_{N}^{q h^{\prime}}+g_{N}^{k p} h_{N}^{q h} g_{N}^{p h^{\prime}} h_{N}^{k^{\prime} q}\right) \mathrm{d}\left \langle B^{pq} \right \rangle _{t} +\sum_{p, q}\left(h_{N}^{k p} g_{N}^{q h} g_{N}^{k^{\prime} q} h_{N}^{p h^{\prime}}+h_{N}^{k p} g_{N}^{q h} h_{N}^{k^{\prime} p} g_{N}^{q h^{\prime}}\right) \mathrm{d}\left \langle B^{qp} \right \rangle _{t}  \\
= &\Big[ (g_{N}^{2})_{k k'} (h_{N}^{2})_{hh'}+(g_{N}^{2})_{kh'} (h_{N}^{2})_{h k'} +(g_{N}^{2})_{h k'} (h_{N}^{2})_{k h'}+(g_{N}^{2})_{h h'} (h_{N}^{2})_{k k'}\Big]  \mathrm{d}\langle B^{11} \rangle _{t}.
\end{align*}
Using the above identity together with~\eqref{e:UfU}, we get for $t<\tau$, 
\begin{align*}
d\lambda_i^N(t)
=&\mathrm{d} M_{i}+ b_{N}(\lambda^{N}_i)\mathrm{d}\langle B_{hh} \rangle_t\\
&+\sum_{j: j \neq i} \frac{1}{\lambda^{N}_{i}-\lambda^{N}_{j}}\Big[g_{N}^{2}(\lambda_i^{N}) h_{N}^{2}(\lambda_j^{N}) +g_{N}^{2}(\lambda_j^{N}) h_{N}^{2}(\lambda_i^{N})\Big]\mathrm{d}\langle B^{11}  \rangle _{t}\\
=&\mathrm{d} M_{i}+ b_{N}(\lambda^{N}_i)\mathrm{d}\langle B^{11} \rangle_t+\sum_{j:j \neq i} \frac{H_{N}\big(\lambda^{N}_{i}, \lambda^{N}_{j}\big)}{\lambda^{N}_{i}-\lambda^{N}_{j}} \mathrm{d}\langle B^{11}  \rangle _{t}.
 \end{align*}
Furthermore,  $dM_i$ given by~\eqref{e:dMi} can be rewritten as
\begin{align*}
    \mathrm{d} M_{i}
    &=(U^{\top}g_N(X_t^N)\mathrm{d} B_{t} h_N(X_t^N)U)_{ii}+(U^{\top}h_N(X_t^N)\mathrm{d} B^{\top}_{t} g_N(X_t^N)U)_{ii}\\
    &=( U^{\top}g_N(X_t^N)UU^{\top}\mathrm{d}B_t UU^{\top}h_N(X_t^N)U)_{ii}+(U^{\top}h_N(X_t^N)UU^{\top}\mathrm{d}B_t^\mathrm{T}UU^{\top} g_N(X_t^N)U)_{ii}\\
    &=g_N(\lambda_i^N)h_N(\lambda_i^N)(U^{\top}\mathrm{d}B_t U)_{ii}+g_N(\lambda_i^N)h_N(\lambda_i^N)(U^{\top}\mathrm{d}B^{T}_t U)_{ii}\\
    &=g_N(\lambda_i^N)h_N(\lambda_i^N)\sum_{k,h=1}^{N}U^{\top}_{ik}(\mathrm{d}B^{kh}_t+\mathrm{d}B^{hk}_t)U_{hi}.
    \end{align*}

Denote for $i=1, \dots, N$,
\begin{equation}\label{e:barB}
\mathrm{d}\bar{B}_{i}(t):=\frac{1}{2}\sum_{k,h=1}^{N}U^{\top}_{ik}(\mathrm{d}B^{kh}_t+\mathrm{d}B^{hk}_t)U_{hi}.
\end{equation}
Then, we have
\begin{equation}\label{e:d<B>}
\begin{aligned}
    \mathrm{d}\langle\bar{B}_{i},\bar{B}_{j}\rangle_t
    &=\frac{1}{ 4}\sum_{k,h=1}^{N}U^{\top}_{ik}U^{\top}_{jk}U_{hi}U_{hj}\mathrm{d}\langle B^{kh}\rangle_t+\frac{1}{ 4}\sum_{k,h=1}^{N}U^{\top}_{ik}U^{\top}_{jh}U_{hi}U_{kj}\mathrm{d}\langle B^{hk}\rangle_t\\
    &\quad +\frac{1}{ 4}\sum_{k,h=1}^{N}U^{\top}_{ik}U^{\top}_{jh}U_{hi}U_{kj}\mathrm{d}\langle B^{kh}\rangle_t+\frac{1}{ 4}\sum_{k,h=1}^{N}U^{\top}_{ik}U^{\top}_{jk}U_{hi}U_{hj}\mathrm{d}\langle B^{hk}\rangle_t\\
    &=\frac{1}{ 4}(U^{\top}U)^2_{ij}\mathrm{d}\langle B^{11}\rangle_t+\frac{1}{ 4}(U^{\top}U)_{ij}(U^{\top}U)_{ji}\mathrm{d}\langle B^{11}\rangle_t\\
    &\quad +\frac{1}{ 4}(U^{\top}U)_{ij}(U^{\top}U)_{ji}\mathrm{d}\langle B^{11}\rangle_t+\frac{1}{ 4}(U^{\top}U)^2_{ij}\mathrm{d}\langle B^{11}\rangle_t\\
    &=\mathbf{1}_{\{i=j\}}\mathrm{d}\langle B^{11}\rangle_t.
    \end{aligned}
    \end{equation}
Consequently, we get that for $t<\tau$,  and $i=1,2,\dots,N$,
\begin{equation*}
\mathrm{d} \lambda^{N}_{i}=2 g_{N}\left(\lambda^{N}_{i}\right) h_{N}\left(\lambda^{N}_{i}\right) \mathrm{d} \bar B_i+b_{N}\left(\lambda^{N}_{i}\right) \mathrm{d}\langle \bar B_i \rangle _{t}+\sum_{j:j \neq i} \frac{H_{N}\big(\lambda^{N}_{i}, \lambda^{N}_{j}\big)}{\lambda^{N}_{i}-\lambda^{N}_{j}} \mathrm{d}\langle \bar B_i  \rangle _{t}.
\end{equation*}

Now we prove that the first collision time $\tau$ is infinite  q.s.   We first introduce some auxiliary functions.

We consider the following Lyapunov function,  which is inspired by \cite{agz10}, 
\begin{equation*}f(x)=f\left(x_{1}, \ldots, x_{N}\right)=\frac{1}{N} \sum_{i=1}^{N} x_{i}^{2}-\frac{1}{N^{2}} \sum_{i \neq j} \log \left|x_{i}-x_{j}\right|,
\end{equation*}
defined for   $x=\left(x_{1}, \ldots, x_{N}\right) \in \Delta_N$. It is not hard to verify the following property of $f$, for $i\neq j$,
\begin{align}\label{e:property-f}
&f\left(x_{1}, \ldots, x_{N}\right)\geq -4,\quad-\frac{1}{N^{2}}\log \left|x_{i}-x_{j}\right|\leq f(x)+4, \quad\frac{1}{N}\sum_{k=1}^N x_k^2 \le 2 f(x) + 4.\end{align}

For any  $M>\max\{0,f(\lambda(0))\}$, set $R:=R(N, M)=e^{N^{2}(4+M)}$ and define the stopping times:
\begin{equation} \label{T_M} 
\tau_R:=\inf \big\{t>0: \min _{i \neq j}\left|\lambda^{N}_{i}(t)-\lambda^{N}_{j}(t)\right|< R^{-1}\big\} \text{ and }
T_{M}:=\inf \left\{t \geq 0: f\left(\lambda^{N}(t)\right) \geq M\right\}.
\end{equation}
Then, by~\eqref{e:property-f} we have for any $i\neq j$,
\begin{equation*}\log\left|\lambda^{N}_{i}(t)-\lambda^{N}_{j}(t)\right| \geq -N^{2}(f(\lambda^{N}(t))+4)\geq-N^{2}(M+4), \text{ for } t\le T_M.\end{equation*}
This implies that $T_{M}\leq\tau_R$. Note that when $t\leq T_{M}$,  $\lambda^{N}_{i}(t)$ is  bounded  by~\eqref{e:property-f}, and hence $g_{N}(\lambda^{N}_{i}), h_{N}(\lambda^{N}_{i}), b_{N}(\lambda^{N}_{i})$ and $H_{N}(\lambda^{N}_{i},\lambda^{N}_{j})$ are also bounded.

 For  $i=1, \ldots, N$,  define the functions
\begin{equation*}u_{i, 1}(x)=\sum_{k: k \neq i} \frac{1}{x_{i}-x_{k}}, \quad u_{i, 2}(x)=\sum_{k: k \neq i} \frac{1}{\left(x_{i}-x_{k}\right)^{2}}, \quad x\in \Delta_N.
\end{equation*}
Note that for  $x \in \Delta_{N} $,
$$A:=\sum_{
\substack{i\neq k, i\neq l,\\ k\neq l}} \frac{1}{(x_i-x_l)(x_i-x_k)}= \sum_{
\substack{i\neq k, i\neq l,\\ k\neq l}}\frac{1}{x_l-x_k}\left(\frac1{x_i-x_l} -\frac1{x_i-x_k}\right) =-2A,$$
and hence $A=0$. Thus, we get 
$$\sum_{i=1}^N\left(u_{i, 1}(x)^{2}-u_{i, 2}(x)\right)=0.$$
Similarly, noting that 
$$ B:= \sum_{i=1}^N\sum_{k\neq i} \frac{x_i}{x_i-x_k}=\sum_{i=1}^N\sum_{k\neq i} \left(1+ \frac{x_k}{x_i-x_k}\right)=N(N-1)-B,$$
we get
\begin{equation*} \sum_{i=1}^{N} u_{i, 1}(x) x_{i}=\frac{N(N-1)}{2}.
\end{equation*}

 Noting that $$\frac{\partial}{\partial x_i} f(x) = \frac{2}{N}x_i - \frac{2}{N^2} u_{i,1}(x), \frac{\partial^2}{\partial x_i^2} f(x) = \frac{2}{N} + \frac{2}{N^2} u_{i,2}(x),$$
we get for $t\leq T_{M}$, 
\begin{align*}
\mathrm{d} f\left(\lambda^{N}(t)\right)= &  \sum_{i=1}^{N}\frac{2}{N}\lambda^{N}_{i}(t)\Big(\sum_{j:j \neq i} \frac{H_{N}\big(\lambda^{N}_{i}, \lambda^{N}_{j}\big)}{\lambda^{N}_{i}-\lambda^{N}_{j}} \mathrm{d}\langle \bar B_i \rangle _{t}+b_{N}\left(\lambda^{N}_{i}\right) \mathrm{d}\langle \bar B_i \rangle _{t}\Big)\notag\\
&-\sum_{i=1}^{N}\frac{2}{N^{2}} u_{i, 1}\left(\lambda^{N}(t)\right)\Big(\sum_{j:j \neq i} \frac{H_{N}\big(\lambda^{N}_{i}, \lambda^{N}_{j}\big)}{\lambda^{N}_{i}-\lambda^{N}_{j}} \mathrm{d}\langle \bar B_i \rangle _{t}+b_{N}\left(\lambda^{N}_{i}\right) \mathrm{d}\langle \bar B_i \rangle _{t}\Big) \notag\\
 &+\frac{1}{2} \sum_{i=1}^{N}\Big(\frac{2}{N}+\frac{2}{N^{2}} u_{i, 2}\left(\lambda^{N}(t)\right)\Big)4 g_{N}^{2}\left(\lambda^{N}_{i}\right) h_{N}^{2}\left(\lambda^{N}_{i}\right)\mathrm{d}\langle \bar B_i \rangle _{t}+\mathrm{d}\tilde M(t),
\end{align*} 
where 
\begin{equation*}\mathrm{d} \tilde M(t)= \sum_{i=1}^{N}\Big(\frac{2}{N}\lambda_{i}^{N}(t)-\frac{2}{N^{2}}  u_{i, 1}\left(\lambda^{N}(t)\right)\Big) 2 g_{N}\left(\lambda^{N}_{i}\right) h_{N}\left(\lambda^{N}_{i}\right) \mathrm{d} \bar B_i,
\end{equation*}
for which we have  
\begin{equation}\label{e:EM}
    \hat{\mathbb{E}}\Big[ \int_{0}^{t \wedge T_{M}} \mathrm{d}\tilde M(s)\Big]=\hat{\mathbb{E}}\Big[\int_{0}^{t \wedge T_{M}} \sum_{i=1}^{N}\Big(\frac{2}{N}\lambda^{N}_{i}(s)-\frac{2}{N^{2}}  u_{i, 1}\left(\lambda^{N}(s)\right)\Big) 2 g_{N}\left(\lambda^{N}_{i}\right) h_{N}\left(\lambda^{N}_{i}\right) \mathrm{d} \bar B_i\Big]=0.
\end{equation}

We write $\mathrm{d}f(\lambda^N(t))$  as follows 
\begin{equation}\label{e:df}
 \mathrm{d}f(\lambda^N(t)) =\Big[ B(\lambda^N(t))+\sum_{i=1}^3 S_i(\lambda^N(t)) \Big] \dd \langle  B^{11} \rangle_t + \dd \tilde M_t
 \end{equation} 
where 
\begin{align*}
B(\lambda^{N}(t))&=\sum_{i=1}^{N}\Big(\frac{2}{N}\lambda^{N}_{i}(t)-\frac{2}{N^{2}}  u_{i, 1}\left(\lambda^{N}(t)\right)\Big)b_{N}\left(\lambda^{N}_{i}\right);\\
S_{1}(\lambda^{N}(t))&=\sum_{i=1}^{N}\frac{2}{N}\lambda^{N}_{i}(t)\sum_{j:j \neq i} \frac{H_{N}\big(\lambda^{N}_{i}, \lambda^{N}_{j}\big)}{\lambda^{N}_{i}-\lambda^{N}_{j}};\\
S_{2}(\lambda^{N}(t))&=-\sum_{i=1}^{N}\frac{2}{N^{2}} u_{i, 1}\left(\lambda^{N}(t)\right)\sum_{j:j \neq i} \frac{H_{N}\big(\lambda^{N}_{i}, \lambda^{N}_{j}\big)}{\lambda^{N}_{i}-\lambda^{N}_{j}} + \sum_{i=1}^{N}\frac{4}{N^{2}} u_{i, 2}\left(\lambda^{N}(t)\right) g_{N}^{2}\left(\lambda^{N}_{i}\right) h_{N}^{2}\left(\lambda^{N}_{i}\right);\\
S_{3}(\lambda^{N}(t))&= \sum_{i=1}^{N}\frac{4}{N} g_{N}^{2}\left(\lambda^{N}_{i}\right) h_{N}^{2}\left(\lambda^{N}_{i}\right).
\end{align*}

We first estimate $B(\lambda^{N}(t)).$
Noting that \begin{align*}
\sum_{i=1}^{N}  u_{i, 1}\left(\lambda^{N}(t)\right)b_{N}\left(\lambda^{N}_{i}\right)=\sum_{i\neq k}\frac{b_{N}\left(\lambda^{N}_{i}\right)}{\lambda^{N}_{i}-\lambda^{N}_{k}}=\sum_{i<k}\frac{b_{N}\left(\lambda^{N}_{i}\right)-b_{N}\left(\lambda^{N}_{k}\right)}{\lambda^{N}_{i}-\lambda^{N}_{k}},
\end{align*}
we get
$$\Big|\sum_{i=1}^{N}  u_{i, 1}\left(\lambda^{N}(t)\right)b_{N}\left(\lambda^{N}_{i}\right)\Big|\leq\frac{N(N-1)}{2}L_{b_{N}},$$ where we use  $L_f$ to denote the Lipschitz constant for a Lipschitz-continuous function $f$.  We also have $$\bigskip
\Big|\sum_{i=1}^{N}\lambda^{N}_{i}(t)b_{N}\left(\lambda^{N}_{i}\right)\Big|\leq\sum_{i=1}^{N}|\lambda^{N}_{i}(t)\max\{L_{b_{N}},|b_{N}(0)|\}(|\lambda^{N}_{i}(t)|+1)|,$$
and thus 
\begin{equation}\label{e:est-B}
|B(\lambda^{N}(t))|\leq C \sum_{i=1}^{N}(1+(\lambda^{N}_{i}(t))^{2}) \text{ for some } C>0.
 \end{equation}

For the term $S_{1}(\lambda^{N}(t))$, noting that 
$$\sum_{i=1}^{N}\lambda^{N}_{i}\sum_{j:j \neq i} \frac{H_{N}\big(\lambda^{N}_{i}, \lambda^{N}_{j}\big)}{\lambda^{N}_{i}-\lambda^{N}_{j}}=\sum_{i<j}H_{N}\big(\lambda^{N}_{i}, \lambda^{N}_{j}\big)\big(\frac{\lambda^{N}_{i}}{\lambda^{N}_{i}-\lambda^{N}_{j}}+\frac{\lambda^{N}_{j}}{\lambda^{N}_{j}-\lambda^{N}_{i}}\big)=\sum_{i<j}H_{N}\big(\lambda^{N}_{i}, \lambda^{N}_{j}\big),$$ and that 
$H_{N}\big(\lambda^{N}_{i}, \lambda^{N}_{j}\big)\leq C(1+(\lambda^{N}_{i}(t))^{2}+(\lambda^{N}_{j}(t))^{2})$ for $C>0$, we get 
\begin{equation}\label{e:S1}
S_{1}(\lambda^{N}(t))= \frac{2}{N}\sum_{i<j}H_{N}\big(\lambda^{N}_{i}, \lambda^{N}_{j}\big)\leq C \sum_{i=1}^{N}(1+(\lambda^{N}_{i}(t))^{2})
\end{equation}
for some $C>0$. For $S_3(\lambda^N(t))$, we have the same estimate:
\begin{equation}\label{e:S3}
S_{3}(\lambda^{N}(t))\leq C \sum_{i=1}^{N}(1+(\lambda^{N}_{i}(t))^{2}). 
\end{equation}

Now, we estimate  $S_{2}(\lambda^{N}(t))$. 
\begin{align*}
& \frac{1}{2}S_{2}(\lambda^{N}(t))\\
&=-\frac{1}{N^{2}}\sum_{i=1}^{N}\sum_{k:k \neq i}\frac{1}{\lambda^{N}_{i}-\lambda^{N}_{k}}\sum_{j:j \neq i} \frac{H_{N}\big(\lambda^{N}_{i}, \lambda^{N}_{j}\big)}{\lambda^{N}_{i}-\lambda^{N}_{j}}+\frac{2}{N^{2}}\sum_{i=1}^{N}\sum_{k:k \neq i}\frac{ g_{N}^{2}\left(\lambda^{N}_{i}\right) h_{N}^{2}\left(\lambda^{N}_{i}\right)}{(\lambda^{N}_{i}-\lambda^{N}_{k})^{2}}\\
&=-\frac{1}{N^{2}}\sum_{i=1}^{N}\sum_{k,j: k \ne i\ne j} \frac{H_{N}\big(\lambda^{N}_{i}, \lambda^{N}_{j}\big)}{(\lambda^{N}_{i}-\lambda^{N}_{k})(\lambda^{N}_{i}-\lambda^{N}_{j})}-\frac{1}{N^{2}}\sum_{i=1}^{N}\sum_{k:k \neq i}\frac{ H_{N}\big(\lambda^{N}_{i}, \lambda^{N}_{k}\big)}{(\lambda^{N}_{i}-\lambda^{N}_{k})^{2}}+\frac{2}{N^{2}}\sum_{i=1}^{N}\sum_{k:k \neq i}\frac{ g_{N}^{2}\left(\lambda^{N}_{i}\right) h_{N}^{2}\left(\lambda^{N}_{i}\right)}{(\lambda^{N}_{i}-\lambda^{N}_{k})^{2}}\\
&=-\frac{1}{N^{2}}\sum_{i=1}^{N}\sum_{k,j: k \ne i\ne j} \frac{H_{N}\big(\lambda^{N}_{i}, \lambda^{N}_{j}\big)}{(\lambda^{N}_{i}-\lambda^{N}_{k})(\lambda^{N}_{i}-\lambda^{N}_{j})}-\frac{1}{N^{2}}\sum_{i=1}^{N}\sum_{k:k \neq i}\frac{ H_{N}\big(\lambda^{N}_{i}, \lambda^{N}_{k}\big)-2g_{N}^{2}\left(\lambda^{N}_{i}\right) h_{N}^{2}\left(\lambda^{N}_{i}\right)}{(\lambda^{N}_{i}-\lambda^{N}_{k})^{2}},
\end{align*}
where in the above  $k\neq i\neq j$ actually means $k\neq i, i\neq j, k\neq j$.

The second term on the right-hand side equals
\begin{align*}
&-\frac{1}{N^{2}}\sum_{i<k}\frac{ H_{N}\big(\lambda^{N}_{i}, \lambda^{N}_{k}\big)-2g_{N}^{2}\left(\lambda^{N}_{i}\right) h_{N}^{2}\left(\lambda^{N}_{i}\right)+H_{N}\big(\lambda^{N}_{i}, \lambda^{N}_{k}\big)-2g_{N}^{2}\left(\lambda^{N}_{k}\right) h_{N}^{2}\left(\lambda^{N}_{k}\right)}{(\lambda^{N}_{i}-\lambda^{N}_{k})^{2}}\\
&=\frac{2}{N^{2}}\sum_{i<k}\frac{ \big(g_{N}^{2}\left(\lambda^{N}_{i}\right)-g_{N}^{2}\left(\lambda^{N}_{k}\right)\big)\big(h_{N}^{2}\left(\lambda^{N}_{i}\right)-h_{N}^{2}\left(\lambda^{N}_{k}\right)\big)}{(\lambda^{N}_{i}-\lambda^{N}_{k})^{2}},
\end{align*}
whose absolute value is bounded by 
$\frac{2}{N^{2}}\frac{N(N-1)}{2}L_{g^2_{N}}L_{h^2_{N}}\leq L_{g^2_{N}}L_{h^2_{N}}.$ 

For the first term of $ \frac{1}{2}S_{2}(\lambda^{N}(t))$, by symmetry we have
\begin{align*}&-\frac{1}{N^{2}}\sum_{i=1}^{N}\sum_{k,j: k \ne i\ne j} \frac{H_{N}\big(\lambda^{N}_{i}, \lambda^{N}_{j}\big)}{(\lambda^{N}_{i}-\lambda^{N}_{k})(\lambda^{N}_{i}-\lambda^{N}_{j})}\\
&=-\frac{1}{N^{2}}\sum_{1\leq i<j<k\leq N}\Big[\frac{H_{N}\big(\lambda^{N}_{i}, \lambda^{N}_{j}\big)+H_{N}\big(\lambda^{N}_{i}, \lambda^{N}_{k}\big)}{(\lambda^{N}_{i}-\lambda^{N}_{k})(\lambda^{N}_{i}-\lambda^{N}_{j})}+\frac{H_{N}\big(\lambda^{N}_{i}, \lambda^{N}_{j}\big)+H_{N}\big(\lambda^{N}_{j}, \lambda^{N}_{k}\big)}{(\lambda^{N}_{j}-\lambda^{N}_{k})(\lambda^{N}_{j}-\lambda^{N}_{i})}\\
&\qquad \qquad \qquad\qquad \qquad  +\frac{H_{N}\big(\lambda^{N}_{i}, \lambda^{N}_{k}\big)+H_{N}\big(\lambda^{N}_{j}, \lambda^{N}_{k}\big)}{(\lambda^{N}_{k}-\lambda^{N}_{i})(\lambda^{N}_{k}-\lambda^{N}_{j})}\Big]\\
&=-\frac{1}{N^{2}}\sum_{1\leq i<j<k\leq N}\frac{H_{N}\big(\lambda^{N}_{i}, \lambda^{N}_{k}\big)(\lambda^{N}_{k}-\lambda^{N}_{i})-H_{N}\big(\lambda^{N}_{i}, \lambda^{N}_{j}\big)(\lambda^{N}_{j}-\lambda^{N}_{i})-H_{N}\big(\lambda^{N}_{j}, \lambda^{N}_{k}\big)(\lambda^{N}_{k}-\lambda^{N}_{j})}{(\lambda^{N}_{j}-\lambda^{N}_{i})(\lambda^{N}_{k}-\lambda^{N}_{i})(\lambda^{N}_{k}-\lambda^{N}_{j})}.
\end{align*}
For the numerator, we have
\begin{align*}
&H_{N}\big(\lambda^{N}_{i}, \lambda^{N}_{k}\big)(\lambda^{N}_{k}-\lambda^{N}_{i})-H_{N}\big(\lambda^{N}_{i}, \lambda^{N}_{j}\big)(\lambda^{N}_{j}-\lambda^{N}_{i})-H_{N}\big(\lambda^{N}_{j}, \lambda^{N}_{k}\big)(\lambda^{N}_{k}-\lambda^{N}_{j})\\
&=\left[g_{N}^{2}\left(\lambda^{N}_{i}\right)h_{N}^{2}\left(\lambda^{N}_{i}\right)+g_{N}^{2}\left(\lambda^{N}_{k}\right)h_{N}^{2}\left(\lambda^{N}_{k}\right)-\left(g_{N}^{2}\left(\lambda^{N}_{i}\right)-g_{N}^{2}\left(\lambda^{N}_{k}\right)\right)\left(h_{N}^{2}\left(\lambda^{N}_{i}\right)-h_{N}^{2}\left(\lambda^{N}_{k}\right)\right)\right]\left(\lambda^{N}_{k}-\lambda^{N}_{i}\right)
\\
&\qquad -\big[g_{N}^{2}(\lambda^{N}_{i})h_{N}^{2}(\lambda^{N}_{i})+g_{N}^{2}(\lambda^{N}_{j})h_{N}^{2}(\lambda^{N}_{j})-\big(g_{N}^{2}(\lambda^{N}_{i})-g_{N}^{2}(\lambda^{N}_{j})\big)\big(h_{N}^{2}(\lambda^{N}_{i})-h_{N}^{2}(\lambda^{N}_{j})\big)\big]\big(\lambda^{N}_{j}-\lambda^{N}_{i}\big)\\
&\qquad -\big[g_{N}^{2}(\lambda^{N}_{k})h_{N}^{2}(\lambda^{N}_{k})+g_{N}^{2}(\lambda^{N}_{j})h_{N}^{2}(\lambda^{N}_{j})-\big(g_{N}^{2}(\lambda^{N}_{k})-g_{N}^{2}(\lambda^{N}_{j})\big)\big(h_{N}^{2}\left(\lambda^{N}_{k}\right)-h_{N}^{2}(\lambda^{N}_{j})\big)\big]\big(\lambda^{N}_{k}-\lambda^{N}_{j}\big)\\
&=g_{N}^{2}\left(\lambda^{N}_{i}\right)h_{N}^{2}\left(\lambda^{N}_{i}\right)(\lambda^{N}_{k}-\lambda^{N}_{j})-g_{N}^{2}(\lambda^{N}_{j})h_{N}^{2}(\lambda^{N}_{j})(\lambda^{N}_{k}-\lambda^{N}_{i})+g_{N}^{2}\left(\lambda^{N}_{k}\right)h_{N}^{2}\left(\lambda^{N}_{k}\right)(\lambda^{N}_{j}-\lambda^{N}_{i})\\
&\qquad-\big(g_{N}^{2}\left(\lambda^{N}_{i}\right)-g_{N}^{2}\left(\lambda^{N}_{k}\right)\big)\big(h_{N}^{2}\left(\lambda^{N}_{i}\right)-h_{N}^{2}\left(\lambda^{N}_{k}\right)\big)(\lambda^{N}_{k}-\lambda^{N}_{i})
\\
&\qquad+\big(g_{N}^{2}(\lambda^{N}_{i})-g_{N}^{2}(\lambda^{N}_{j})\big)\big(h_{N}^{2}(\lambda^{N}_{i})-h_{N}^{2}(\lambda^{N}_{j})\big)(\lambda^{N}_{j}-\lambda^{N}_{i})\\
&\qquad+\big(g_{N}^{2}\left(\lambda^{N}_{k}\right)-g_{N}^{2}(\lambda^{N}_{j})\big)\big(h_{N}^{2}\left(\lambda^{N}_{k}\right)-h_{N}^{2}(\lambda^{N}_{j})\big)(\lambda^{N}_{k}-\lambda^{N}_{j})\\
&:= I+II+III+IV. 
\end{align*}

For the term $I$ we have the following estimate: for $i<j<k$, 
\begin{align*}
&\frac{I}{(\lambda^{N}_{j}-\lambda^{N}_{i})(\lambda^{N}_{k}-\lambda^{N}_{i})(\lambda^{N}_{k}-\lambda^{N}_{j})}\\
=&\frac{\big[g_{N}^{2}\left(\lambda^{N}_{i}\right)h_{N}^{2}\left(\lambda^{N}_{i}\right)-g_{N}^{2}(\lambda^{N}_{j})h_{N}^{2}(\lambda^{N}_{j})\big](\lambda^{N}_{k}-\lambda^{N}_{j})+\big[g_{N}^{2}\left(\lambda^{N}_{k}\right)h_{N}^{2}\left(\lambda^{N}_{k}\right)-g_{N}^{2}(\lambda^{N}_{j})h_{N}^{2}(\lambda^{N}_{j})\big](\lambda^{N}_{j}-\lambda^{N}_{i})}{(\lambda^{N}_{j}-\lambda^{N}_{i})(\lambda^{N}_{k}-\lambda^{N}_{i})(\lambda^{N}_{k}-\lambda^{N}_{j})}\\
=&\frac{1}{\lambda^{N}_{k}-\lambda^{N}_{i}}\left[\frac{g_{N}^{2}\left(\lambda^{N}_{k}\right)h_{N}^{2}\left(\lambda^{N}_{k}\right)-g_{N}^{2}(\lambda^{N}_{j})h_{N}^{2}(\lambda^{N}_{j})}{\lambda^{N}_{k}-\lambda^{N}_{j}}-\frac{g_{N}^{2}(\lambda^{N}_{j})h_{N}^{2}(\lambda^{N}_{j})-g_{N}^{2}\left(\lambda^{N}_{i}\right)h_{N}^{2}\left(\lambda^{N}_{i}\right)}{\lambda^{N}_{j}-\lambda^{N}_{i}}\right]\geq - C,
\end{align*}
 for  some constant $C>0$, noting that  \(g_{N}^2h_{N}^2\) is either convex or of class \(C^{1,1}\) and $\lambda^N_i<\lambda^N_j<\lambda^N_k$. 
 
For the remaining terms, we have for $i<j<k$, 
\begin{align*}
&\left|\frac{II+III+IV}{(\lambda^{N}_{j}-\lambda^{N}_{i})(\lambda^{N}_{k}-\lambda^{N}_{i})(\lambda^{N}_{k}-\lambda^{N}_{j})}\right|\\
&=\left|-\frac{\big[\big(g_{N}^{2}\left(\lambda^{N}_{k}\right)- g_{N}^{2}(\lambda^{N}_{j})\big)+\big(g_{N}^{2}(\lambda^{N}_{j})- g_{N}^{2}\left(\lambda^{N}_{i}\right)\big)\big]\big[\big(h_{N}^{2}\left(\lambda^{N}_{k}\right)-h_{N}^{2}(\lambda^{N}_{j})\big)+\big(h_{N}^{2}(\lambda^{N}_{j})- h_{N}^{2}\left(\lambda^{N}_{i}\right)\big)\big]}{(\lambda^{N}_{j}-\lambda^{N}_{i})(\lambda^{N}_{k}-\lambda^{N}_{j})}\right.\\
&\quad +\frac{\big(g_{N}^{2}(\lambda^{N}_{j})- g_{N}^{2}\left(\lambda^{N}_{i}\right)\big)\big(h_{N}^{2}(\lambda^{N}_{j})- h_{N}^{2}\left(\lambda^{N}_{i}\right)\big)}{(\lambda^{N}_{k}-\lambda^{N}_{i})(\lambda^{N}_{k}-\lambda^{N}_{j})}\left.+\frac{\big(g_{N}^{2}\left(\lambda^{N}_{k}\right)- g_{N}^{2}(\lambda^{N}_{j})\big)\big(h_{N}^{2}\left(\lambda^{N}_{k}\right)- h_{N}^{2}(\lambda^{N}_{j})\big)}{(\lambda^{N}_{j}-\lambda^{N}_{i})(\lambda^{N}_{k}-\lambda^{N}_{i})}\right|\\
&=\left|-\frac{\big(g_{N}^{2}(\lambda^{N}_{j})- g_{N}^{2}\left(\lambda^{N}_{i}\right)\big)\big(h_{N}^{2}(\lambda^{N}_{j})- h_{N}^{2}\left(\lambda^{N}_{i}\right)\big)}{(\lambda^{N}_{j}-\lambda^{N}_{i})(\lambda^{N}_{k}-\lambda^{N}_{i})}\right.-\left.\frac{\big(g_{N}^{2}\left(\lambda^{N}_{k}\right)- g_{N}^{2}(\lambda^{N}_{j})\big)\big(h_{N}^{2}(\lambda^{N}_{j})- h_{N}^{2}\left(\lambda^{N}_{i}\right)\big)}{(\lambda^{N}_{j}-\lambda^{N}_{i})(\lambda^{N}_{k}-\lambda^{N}_{j})}\right.\\
&\quad \left.- \frac{\big(g_{N}^{2}(\lambda^{N}_{j})- g_{N}^{2}\left(\lambda^{N}_{i}\right)\big)\big(h_{N}^{2}\left(\lambda^{N}_{k}\right)- h_{N}^{2}(\lambda^{N}_{j})\big)}{(\lambda^{N}_{j}-\lambda^{N}_{i})(\lambda^{N}_{k}-\lambda^{N}_{j})}\right.\left.-\frac{\big(g_{N}^{2}\left(\lambda^{N}_{k}\right)- g_{N}^{2}(\lambda^{N}_{j})\big)\big(h_{N}^{2}\left(\lambda^{N}_{k}\right)- h_{N}^{2}(\lambda^{N}_{j})\big)}{(\lambda^{N}_{k}-\lambda^{N}_{i})(\lambda^{N}_{k}-\lambda^{N}_{j})}\right|\\
&\leq L_{g^2_{N}}L_{h^2_{N}}\frac{\lambda^{N}_{j}-\lambda^{N}_{i}}{\lambda^{N}_{k}-\lambda^{N}_{i}}+2L_{g^2_{N}}L_{h^2_{N}}+L_{g^2_{N}}L_{h^2_{N}}\frac{\lambda^{N}_{k}-\lambda^{N}_{j}}{\lambda^{N}_{k}-\lambda^{N}_{i}}=3L_{g^2_{N}}L_{h^2_{N}}.
\end{align*}

By the above analysis for $I$ and $II+III+IV$,  we get the following estimate for  the first term of $ \frac{1}{2}S_{2}(\lambda^{N}(t))$:
$$-\frac{1}{N^{2}}\sum_{1\le i\neq k\neq j\le N}\frac{H_{N}\big(\lambda^{N}_{i}, \lambda^{N}_{j}\big)}{(\lambda^{N}_{i}-\lambda^{N}_{k})(\lambda^{N}_{i}-\lambda^{N}_{j})}\leq \frac{1}{N^{2}}\binom{N}{3}C,$$
 for  some $C>0$.   Then, combining the estimates for the two terms of $ \frac{1}{2}S_{2}(\lambda^{N}(t))$, we get 
 \begin{equation}\label{e:S2}
     S_2(\lambda^N(t))\le C, 
 \end{equation}
 for some constant $C$ that may depend on $N$.

Therefore, by the equation~\eqref{e:df} of $\dd f(\lambda^N(t))$ and the estimates \eqref{e:est-B},  \eqref{e:S1}, \eqref{e:S3} and \eqref{e:S2}, we obtain 
\begin{align*}
    \mathrm{d} f\left(\lambda^{N}(t)\right)&\leq C \sum_{i=1}^{N}(1+(\lambda^{N}_{i}(t))^{2})\mathrm{d}t+\mathrm{d}\tilde M(t)= [CN+C\sum_{i=1}^{N}(\lambda^{N}_{i}(t))^{2}]\mathrm{d}t+\mathrm{d}\tilde M(t)\\
    &\leq C( f\left(\lambda^{N}(t)\right)+1)\mathrm{d}t+\mathrm{d}\tilde M(t),\end{align*}
where the last step follows from~\eqref{e:property-f}. Taking the integral and then the expectation, by~\eqref{e:EM} we obtain  
\begin{align*}    \hat{\mathbb{E}}\left[f\left(\lambda^{N}\left(t \wedge T_{M}\right)\right)\right]& \leq C\hat{\mathbb{E}}\Big[\int_{0}^{t \wedge T_{M}}f\left(\lambda^{N}\left(s\right)\right)\mathrm{d}s\Big]+Ct+f\left(\lambda^{N}(0)\right)\\
&\leq C\int_{0}^{t}\hat{\mathbb{E}}\left[f\left(\lambda^{N}\left(s \wedge T_{M}\right)\right)\right]\mathrm{d}s+Ct+f\left(\lambda^{N}(0)\right).
\end{align*}
Then, Gronwall's inequality yields
$\hat{\mathbb{E}}\left[f\left(\lambda^{N}\left(t \wedge T_{M}\right)\right)\right] \leq C_{t}$, where $C_{t}$ is a finite number only depending on $t,N$ and $\lambda^N(0)$. Noting that $0\le f\left(\lambda^{N}\left(t \wedge T_{M}\right)\right)+4\le M+4$ q.s.\ and that $f\left(\lambda^{N}\left(t \wedge T_{M}\right)\right)=M$ on the event $[T_M\le t]$, we get
\begin{align*} 
 c \left(T_{M} \leq t\right) (M+4) 
& \leq \hat{\mathbb{E}}\left[f\left(\lambda^{N}\left(t \wedge T_{M}\right)\right)+4\right] \leq C_{t}+4,
\end{align*}
and hence
\begin{equation*}c\left(T_{M}\leq t\right) \leq \frac{C_{t}+4}{M+4} .\end{equation*}

Denote $A_{n}^{t}=\left\{T_{n^2}\leq t\right\}$ for $n\in \mathbb{N}.$  Then, we have
$$\sum_{n=1}^{\infty}c(A_{n}^{t})=\sum_{n=1}^{\infty}c(T_{n^2}\leq t)\leq\sum_{n=1}^{\infty} \frac{C_{t}+4}{n^2+4}<\infty,$$
and Lemma~\ref{lem:B-C} yields   that $\limsup_{n\to \infty}A_{n}^{t}=\bigcap_{n=1}^{\infty}A_{n}^{t}$  is a polar set, i.e., $c (\bigcap_{n=1}^{\infty}A_{n}^{t})=0$.  Thus, we have $c(\bigcup_{k=1}^{\infty}\bigcap_{n=1}^{\infty}A_{n}^{k})=\lim_{k\to\infty}c(\bigcap_{n=1}^{\infty}A_{n}^{k})=0$. Consequently, we get
$$c\Big(\Big\{\lim_{R\to\infty}\tau_R=\infty\Big\}^{c}\Big)\leq c\Big(\left\{\lim_{n\to\infty}T_{n^2}=\infty\right\}^{c}\Big)=c(\bigcup_{k=1}^{\infty}\bigcap_{n=1}^{\infty}A_{n}^{k})=0,$$
since $T_{n^{2}}\leq \tau_R$ with $R=e^{N^{2}(4+n^2)}$.  This implies $\tau=\infty$ q.s., since   $\tau_{R}\leq\tau$ for all $R>0$. 

 Finally, noting that $\tau=\infty$ q.s.,  
by assumption, $X_t^N$ has a q.c.\ version which is still denoted by  ${X}_t^N$, then 
$U^{\top}_{ik}(X_t^N)U_{hi}(X_t^N)=(P_i(X_t^N))_{kh}$ is also q.c.,
where $P_i(X)$ is the orthogonal spectral projection onto the
eigenspace of $\lambda_i$ defined by
$P_i(X)=\prod_{j\ne i}\frac{X-\lambda_j(X)I_N}{\lambda_i(X)-\lambda_j(X)}$. Thus, it is easy to check that $\bar B( t )=(\bar{B_{i}}  (t), i=1, \dots, N)$ given by~\eqref{e:barB} is a symmetric martingale. Then, by Corollary~\ref{corollary:G-brownian motion} and \eqref{e:d<B>},  we get that $\bar B( t )$ is an $N$-dimensional  $\bar G$-Brownian motion  satisfying Assumption~\ref{assum:G} with $d=N$. 
\end{proof}


Now we are ready to show that the $G$-SDE~\eqref{eq:eigenvalue-dynamics} admits a unique strong solution.

\begin{proof}[Proof of Theorem~\ref{thm:strong-sol}]\makeatletter\def\@currentlabelname{the proof}\makeatother\label{proof:thm:strong-sol}

    
Denote $A_t:=\langle \bar B^1\rangle_t$. Then $0<\underline\sigma^2(t-s)\le A_t-A_s\le \overline\sigma^2(t-s),\, 0\le s<t$.
For $u\ge0$, define $T_u:=\inf\{t\ge0:A_t>u\}$, $W_u:=\bar B_{T_u}$.
Then it is easy to know that $T_u<\infty$, and $A_{T_u}=u$. For any fixed $P\in\mathscr P$, one can verify that
$\langle W^i,W^j\rangle_u=\langle \bar B^i,\bar B^j\rangle_{T_u}=
\mathbf{1}_{\{i=j\}}A_{T_u}=\mathbf{1}_{\{i=j\}}u$.
Thus, $W$ is a standard $N$-dimensional Brownian motion under
$P$, with respect to the time-changed filtration $\mathcal F_{T_u}$.

Denote $F_i(x):=\sum_{j:j\neq i}\frac{H_{N}(x_i,x_j)}{x_i-x_j}$.
Now we consider the following classical SDE:
\begin{equation}\label{eq:classical-eigenvalue-sde}
\mathrm{d}Y_i=2g_N(Y_i)h_N(Y_i)\mathrm{d}W^{i}_u+(b_N(Y_i)+F_i(Y))\mathrm{d}u,\, Y(0)=\lambda(0)\in\Delta_N.
\end{equation}

Applying Theorem 2 in \cite{GraczykMalecki2013} through the same localization argument as in the proof of Theorem 6  in \cite{GraczykMalecki2013}, and combining weak
existence with the Yamada--Watanabe theorem, we obtain  a global pathwise unique strong
solution.
Equivalently, there is a nonanticipative Borel map
$\Phi:C_0([0,\infty);\mathbb R^N)\longrightarrow C([0,\infty);\Delta_N)$ such that $Y=\Phi(W)$ is the strong solution of
\eqref{eq:classical-eigenvalue-sde}. 

Define $\lambda(t):=Y(A_t)=\Phi(W)(A_t)$.
Then the  time-change formula for stochastic integrals and ~\eqref{eq:classical-eigenvalue-sde}
yields the following identity
\begin{align*}
    \lambda_i(t)=Y_i(A_t)&=\lambda_i(0)+
  \int_0^{A_t}2g_N(Y_i(u))h_N(Y_i(u))\mathrm{d}W^{i}_u+
  \int_0^{A_t}[b_N(Y_i(u))+F_i(Y(u))]\mathrm{d}u\notag\\
&=\lambda_i(0)+
  \int_0^t2g_N(Y_i(A_s))h_N(Y_i(A_s))\mathrm{d}\bar B^{i}_s+
  \int_0^t [b_N(Y_i(A_s))+F_i(Y(A_s))]\mathrm{d}A_s\notag\\
  &=\lambda_i(0)+
  \int_0^t2g_N(\lambda_i(s))h_N(\lambda_i(s))\mathrm{d}\bar B^{i}_s+
  \int_0^t [b_N(\lambda_i(s))+F_i(\lambda(s))]\mathrm{d}A_s,
\end{align*}
holds $P$-a.s.  Since $P\in\mathscr P$ is arbitrary, we obtain a strong solution to ~\eqref{eq:eigenvalue-dynamics}.


Finally, let $\tilde\lambda$ be another strong solution of
~\eqref{eq:eigenvalue-dynamics} with the same initial condition.  Define
$\tilde Y_u:=\widetilde\lambda_{T_u}$.
Applying the inverse time-change formula shows that $\tilde Y$ solves
the same classical SDE \eqref{eq:classical-eigenvalue-sde}. Then 
$\tilde Y_u=Y_u,\,u\ge0$ outside a $P$-null set.  Since $T_{A_t}=t$, it follows that
$\tilde\lambda(t)=\tilde Y(A_t)=Y(A_t)=\lambda(t),\,t\ge0$,
outside a $P$-null set.  As this argument applies to every
$P\in\mathscr P$, the solution is pathwise unique
quasi-surely.
\end{proof}

\subsection{Discussion of \texorpdfstring{$G$}{}-Dyson Brownian motion}
 \label{subsection:Discussion of \texorpdfstring{$G$}{}-Dyson Brownian motion}

By Theorem~\ref{thm:strong-sol}, $G$-Dyson Brownian motion~\eqref{e:G-Dyson} with $\lambda^N(0)\in\Delta_N$  has a unique strong solution. Below we will prove that, similar to classical Dyson Brownian motion (see~\cite{agz10}),  it also admits a unique strong solution when  $\lambda^N(0)\in\overline{\Delta}_N$.  To this end, we first
establish the following comparison result.

\begin{lemma}\label{lem:com}
Assume Assumption~\ref{assum:G} for the underlying $G$-Brownian motion. Let  $\lambda^{N}(t)$ and  $\eta^{N}(t)$  be two solutions to the G-SDE~\eqref{e:G-Dyson} starting from  $\lambda^{N}(0)\in \Delta_N$ and $\eta^{N}(0) \in \Delta_{N}$, respectively.  Assume   $\lambda_{i}^{N}(0)< \eta_{i}^{N}(0)$  for each $i\in\{1, \dots, N\}$. Then,
\begin{equation*}c\left(\{\lambda_{i}^{N}(t)<\eta_{i}^{N}(t) \text{ for all } t \geq 0 \text { and } i=1, \ldots, N\}^c\right)=0 .
\end{equation*}
\end{lemma}
\begin{proof} The argument essentially follows from the proof of \cite[Lemma 4.3.6]{agz10}, which is provided here for the reader's convenience. 

First, we have
$\mathrm{d}\big(\sum_{i} \lambda_{i}^{N}(t)-\sum_{i} \eta_{i}^{N}(t)\big)=0,$ since the interaction terms cancel pairwise in the summation, noting that $\langle \bar B_i\rangle_t=\langle \bar B_1\rangle_t$ for all $i=1, \dots, N$.
Consequently, we get
\begin{equation}\label{eq:two}\sum_{i}\left(\lambda_{i}^{N}(t)-\eta_{i}^{N}(t)\right)=\sum_{i}\left(\lambda_{i}^{N}(0)-\eta_{i}^{N}(0)\right)<0 .\end{equation}
On the other hand, for each $i=1,\ldots,N$, we have
\begin{equation}
\begin{aligned}
\mathrm{d}\bigl(\lambda_i^N(t)-\eta_i^N(t)\bigr)
={}&
\frac1N\sum_{j:j\ne i}
\frac{
\eta_i^N(t)-\lambda_i^N(t)
-\eta_j^N(t)+\lambda_j^N(t)
}{
\bigl(\eta_i^N(t)-\eta_j^N(t)\bigr)
\bigl(\lambda_i^N(t)-\lambda_j^N(t)\bigr)
}
\,\mathrm{d}\langle \bar{B}_1\rangle_t.
\end{aligned}
\label{eq:difference-equation2}
\end{equation}
Define the first contact time
$$
T
:=
\inf\left\{
t>0:
\lambda_i^N(t)-\eta_i^N(t)=0
\text{ for some }i\in\{1,\ldots,N\}
\right\}.
$$
If $T<\infty$, then by the continuity of the paths, we get $\lambda_i^N(t)-\eta_i^N(t)<0, t<T$ and for all $i=1,\ldots,N$.
Let  $i_0$ be such that
$\lambda_{i_0}^N(T)-\eta_{i_0}^N(T)=0$.
Then
$\lambda_j^N(T)-\eta_j^N(T)\le0,
\, j=1,\ldots,N$.
Moreover, by \eqref{eq:two}, these differences cannot  be all zero.
Hence, there exists $j_0\ne i_0$ such that
$\lambda_{j_0}^N(T)-\eta_{j_0}^N(T)<0$.
Since both $\lambda^N$ and $\eta^N$ remain in the strictly ordered region $\Delta_N$, the
denominator in \eqref{eq:difference-equation2} is strictly positive. At $t=T$, the
numerator corresponding to $j_0$ is given by
\[
\eta_{i_0}^N(T)-\lambda_{i_0}^N(T)
-\eta_{j_0}^N(T)+\lambda_{j_0}^N(T)
=
\lambda_{j_0}^N(T)-\eta_{j_0}^N(T)<0.
\]
Therefore, by continuity we have, for $t<T$ sufficiently close to $T$,
$\dd(\lambda_{i_0}^N(t)-\eta_{i_0}^N(t))<0$ and thus $\lambda_{i_0}^N(t)-\eta_{i_0}^N(t)>0$. This contradicts the definition of $T$. Thus, $T=\infty$ q.s., and the desired result follows.
\end{proof}

We now establish the existence and uniqueness of a solution to \eqref{e:G-Dyson} with initial condition $\lambda^N(0)\in\overline{\Delta}_N$, thereby relaxing the assumption $\lambda^N(0)\in\Delta_N$ in Theorem~\ref{thm:dyson non-collision}. It turns out that, as in the classical Dyson Brownian motion, the particles of the $G$-Dyson Brownian motion separate instantaneously after the initial time.

\begin{proof}[Proof of Proposition~\ref{prop:closed-chamber}] \makeatletter\def\@currentlabelname{the proof}\makeatother\label{proof:prop:closed-chamber}The proof essentially follows the proof of \cite[Proposition 4.3.5]{agz10}.
For $n\in\mathbb{N}$, denote
$$\lambda_i^{N,n}(0)
:=
\lambda_i^N(0)+\frac{i}{n},\, \lambda_i^{N,-n}(0)
:=
\lambda_i^N(0)-\frac{N+1-i}{n}\text{ for } 
 i=1,\ldots,N.$$
Then \(\lambda^{N,n}(0),\lambda^{N,-n}(0)\in\Delta_N\). Moreover, for every \(n\ge1\),
$$\lambda_i^{N,-n}(0)<\lambda_i^{N,-n-1}(0)
<\lambda_i^{N,n+1}(0)<\lambda_i^{N,n}(0),
\quad i=1,\ldots,N.$$
Let \(\lambda^{N,n}\) and  \(\lambda^{N,-n}\) denote the solution to \eqref{e:G-Dyson} starting from these respective initial conditions. By Lemma~\ref{lem:com}, we have  
$$\lambda_i^{N,-n}(t)<\lambda_i^{N,-n-1}(t)
<\lambda_i^{N,n+1}(t)<\lambda_i^{N,n}(t),
\quad i=1,\ldots,N,
\quad\text{q.s.,}$$
for each $t>0$. Furthermore, since the interaction terms cancel after summation, we
have
\begin{equation}
\sum_{i=1}^N
\bigl(\lambda_i^{N,n}(t)-\lambda_i^{N,-n}(t)\bigr)
=
\sum_{i=1}^N
\bigl(\lambda_i^{N,n}(0)-\lambda_i^{N,-n}(0)\bigr),
\quad t\ge0,
\label{eq:approx-sum}
\end{equation}
the right-hand side of which converges to zero as
\(n\to\infty\). Since
$$0\le
\lambda_i^{N,n}(t)-\lambda_i^{N,-n}(t)
\le
\sum_{j=1}^N
\bigl(\lambda_j^{N,n}(t)-\lambda_j^{N,-n}(t)\bigr),$$
it follows that both \(\{\lambda^{N,n}\}_{n\ge1}\) and
\(\{\lambda^{N,-n}\}_{n\ge1}\) are Cauchy sequences on any compact time interval, converging uniformly to the same limit denoted by
\(\lambda^N\). In particular,
$\lambda^N\in C\bigl(\mathbb{R}_+;\overline{\Delta}_N\bigr)$.

 We next show that $\lambda^N$ is a solution of ~\eqref{e:G-Dyson}, for which it suffices to prove
\begin{equation}\label{eq:instantaneous-separation}
\lambda^N(t)\in\Delta_N \text{ for all } t>0 \text{ q.s.}
\end{equation}
 Indeed, if \eqref{eq:instantaneous-separation} holds, then starting at any  time
$s_0>0$,  the solution of ~\eqref{e:G-Dyson} starting from $\lambda^N(s_0)$ is bounded  by $\lambda^{N,n}$ and $\lambda^{N,-n}$ for all  $n$ sufficiently large, which forces the solution to coincide with the limit $\lambda^N(t),\,t \ge s_0$.

Suppose on the contrary that there exists \(t'>0\) such that
$\lambda^N(s)\in\overline{\Delta}_N\setminus\Delta_N,
 0<s\le t'.$  Then, we can find $i_0\in\{2,\dots,N\}$ and $t\in(0,t']$ such that $\lambda_{k}^N(s)<\lambda_{k+1}^N(s)$ for $k\geq i_0$ and
$\lambda_{i_0-1}^N(s)=\lambda_{i_0}^N(s),
 0<s\le t$.   More precisely, there exist a constant \(C>0\) and a sequence
\(\varepsilon_n\downarrow0\) such that, for \(n\) sufficiently large,
$$\big|\lambda_k^{N,n}(s)-\lambda_{k+1}^{N,n}(s)\big|
\ge C,\quad k\ge i_0,\quad 0<s\le t,$$
and
$$\big|
\lambda_{i_0}^{N,n}(s)-\lambda_{i_0-1}^{N,n}(s)\big|\le\varepsilon_n,
\quad 0<s\le t.$$
By \eqref{e:G-Dyson} we get
\[
\begin{aligned}
\lambda_{i_0}^{N,n}(s)
\ge{}
\lambda_{i_0}^{N,n}(0)
+\sqrt{\frac{2}{N}}\,\bar{B}_{i_0}(s) +\frac1N\big(\frac1{\varepsilon_n}
-\frac{N-i_0}{C}\big) \langle \bar{B}_{i_0}\rangle_s.
\end{aligned}
\]
Letting \(n\to\infty\), we obtain
$\lambda_{i_0}^{N,n}(s)\rightarrow\infty,
\,\text{q.s.}$,
which contradicts the uniform bound~\eqref{eq:uniform-energy-bound} established below.

Denote
$C_N(n,t)
:=
\frac1N\sum_{i=1}^N\bigl(\lambda_i^{N,n}(t)\bigr)^2$.
We claim 
\begin{equation}\label{eq:uniform-energy-bound}
\sup_{n\ge1}\sup_{s\in[0,t]}
\sqrt{C_N(n,s)}<\infty,
\quad\text{q.s.}
\end{equation}

Let \(n\ge n'\). By Lemma~\ref{lem:com},
$\lambda_i^{N,n}(s)\le\lambda_i^{N,n'}(s),
 \,i=1,\ldots,N.$
Moreover, using the triangle inequality in \(\mathbb{R}^N\), together
with \eqref{eq:approx-sum}, we obtain
\[
\begin{aligned}
\Big|
\sqrt{C_N(n,s)}-\sqrt{C_N(n',s)}
\Big|
&\le
\left(
\frac1N\sum_{i=1}^N
\bigl(
\lambda_i^{N,n}(s)-\lambda_i^{N,n'}(s)
\bigr)^2
\right)^{1/2} \\
&\le
\frac1{\sqrt N}
\sum_{i=1}^N
\bigl(
\lambda_i^{N,n'}(s)-\lambda_i^{N,n}(s)
\bigr) \\
&=
\frac1{\sqrt N}
\sum_{i=1}^N
\bigl(
\lambda_i^{N,n'}(0)-\lambda_i^{N,n}(0)
\bigr).
\end{aligned}
\]
Consequently,
$$\sup_{n\ge n'}\sup_{s\in[0,t]}
\sqrt{C_N(n,s)}
\le
\sup_{s\in[0,t]}\sqrt{C_N(n',s)}
+
\frac1{\sqrt N}
\sup_{n\ge n'}
\sum_{i=1}^N
\bigl(
\lambda_i^{N,n'}(0)-\lambda_i^{N,n}(0)
\bigr).$$
Thus, it remains to bound
\(\sup_{s\in[0,t]}C_N(n',s)\) for a fixed \(n'\).

$G$-It\^{o}'s formula yields
$$C_N(n,t)=D_N(n,t)+\frac{2\sqrt{2}}{N\sqrt N}\sum_{i=1}^N
\int_0^t\lambda_i^{N,n}(s)\,\mathrm{d}\bar{B}_i(s),$$
where
\[
\begin{aligned}
D_N(n,t)
:=
C_N(n,0)
+
\frac{2}{N}
\sum_{i=1}^N
\int_0^t
\lambda_i^{N,n}(s)
\frac1N
\sum_{j:j\ne i}
\frac{1}{
\lambda_i^{N,n}(s)-\lambda_j^{N,n}(s)
}
\,\mathrm{d}\langle \bar{B}_i\rangle_s +
\frac1N
\sum_{i=1}^N
\int_0^t\frac{2}{N}\,\mathrm{d}\langle \bar{B}_i\rangle_s.
\end{aligned}
\]
By the symmetry of the interaction terms, we get
\[
D_N(n,t)
=
C_N(n,0)
+
\frac{N+1}{N}
\langle \bar{B}_1\rangle_t
\le
C_N(n,0)
+
\frac{N+1}{N}\overline{\sigma}^{\,2}t.
\]

For \(R>0\), define
$S_R:=\inf\{s\ge0:C_N(n,s)\ge R\}$.
By the \(G\)-Burkholder--Davis--Gundy inequality, there exists a
constant \(C>0\), independent of \(R\), such that
\[
\begin{aligned}
\hat{\mathbb E}\big[
\sup_{s\in[0,t]}C_N(n,s\wedge S_R)^2
\big]
&\le
2\hat{\mathbb E}\left[D_N(n,t)^2\right] +
2N^{-2}C
\int_0^t
\hat{\mathbb E}\big[
\sup_{r\in[0,u]}
C_N(n,r\wedge S_R)
\big]\mathrm{d}u\\
&\le 2\hat{\mathbb E}\left[D_N(n,t)^2\right] + N^{-2}C t+N^{-2}C
\int_0^t
\hat{\mathbb E}\big[
\sup_{r\in[0,u]}
C_N(n,r\wedge S_R)^2
\big]\mathrm{d}u.
\end{aligned}
\]
Gronwall's inequality therefore yields
that $\hat{\mathbb E}\big[
\sup_{s\in[0,t]}C_N(n,s\wedge S_R)^2
\big]$ has a finite upper bound which is independent of $R$. Letting \(R\to\infty\), we conclude that
$\hat{\mathbb E}\big[
\sup_{s\in[0,t]}C_N(n,s)^2
\big]
<\infty.$ This proves \eqref{eq:uniform-energy-bound}, and hence \eqref{eq:instantaneous-separation} follows.

Finally, we show that $\lambda^N$ is a unique solution to~\eqref{e:G-Dyson}. For any  fixed $s>0$, 
\(\lambda^N(s)\in\Delta_N\) by~\eqref{eq:instantaneous-separation}, and thus Theorem~\ref{thm:dyson non-collision} implies that there is a unique solution to \eqref{e:G-Dyson} starting from $\lambda^N(s)$. By the uniform convergence of $\lambda^{N,n}(t)\to \lambda^N(t)$ on any compact interval, this solution coincides with
\((\lambda^N(t))_{t\ge s}\). Since \(s>0\) is arbitrary,
\(\lambda^N\) is a solution on \((0,\infty)\).
 Moreover, by the preceding uniform convergence
result, we have
$\lambda^N\in C(\mathbb{R}_+;\overline{\Delta}_N)$,
and hence
$\lim_{t\downarrow0}\lambda^N(t)=\lambda^N(0)
$. Thus, $\lambda^N$ is a solution with initial condition $\lambda^N(0)$. The uniqueness follows directly from Theorem~\ref{thm:dyson non-collision} and~\eqref{eq:instantaneous-separation}. 
\end{proof}

\subsection{Discussion of \texorpdfstring{$G$}{}-Wishart process} 
\label{subsection:Discussion of \texorpdfstring{$G$}{}-Wishart process}

\begin{lemma}\label{lem:G-Wishart}
 Let $X^{N}_{t}$ be given by~\eqref{wishart definition}. 
 Assume that the underlying $\tilde G$-Brownian motion $Y$ satisfies Assumption \ref{assum:G} with $d=p\times N$  and $p\ge N$.  Then $X^{N}_{t}$ satisfies the  $G$-SDE~\eqref{G-N-wishart SDE}.
\end{lemma}
\begin{proof}
Noting $X^{N}_{ij}(t) = \sum_{k=1}^p Y_{ki}(t) Y_{kj}(t)$, $G$-It\^o's formula yields 
\begin{align*}
    \mathrm{d}X^{N}_{ij}(t)& = \sum_{k=1}^p \Big( Y_{ki}(t) \mathrm{d} Y_{kj}(t) +Y_{kj}(t) \mathrm{d}Y_{ki}(t) + \mathrm{d}\langle Y_{ki},Y_{kj} \rangle_t \Big)\\
    &=\sum_{k=1}^p \Big( Y_{ki}(t) \mathrm{d}Y_{kj}(t) + Y_{kj}(t) \mathrm{d}Y_{ki}(t)\Big)+p \mathbf{1}_{\{i=j\}}\mathrm{d}\langle Y_{ki} \rangle_t,
   \end{align*}
which in the matrix form is 
\begin{equation}\label{e:dX}
    \mathrm{d}X^{N}(t)=Y_t^{\top} \mathrm{d}Y_t + \mathrm{d}Y_t^{\top}Y_t+p \mathbf{I}_N \mathrm{d}\langle Y_{11} \rangle_t.
\end{equation}

Define an $N \times N$ matrix-valued process $B_t$ by
\begin{equation}\label{e:dhatB}
\mathrm{d}B_t := \Big( \sqrt{X^{N}_t}\Big)^{-1} Y_t^{\top}\mathrm{d}Y_t,
\end{equation}
where $(\sqrt{X_t^N})^{-1}$ can be well defined by the following argument.

As in the proof of Theorem~\ref{thm:strong-sol}, we can transform $Y$ into a Brownian motion $W$ under each $P\in\mathcal{P}$ by a time change. Denote $A_t:=\langle  Y_{11}\rangle_t$. 
Then $0<\underline\sigma^2(t-s)\le A_t-A_s\le \overline\sigma^2(t-s),\, 0\le s<t$.
For $u\ge0$, define $T_u:=\inf\{t\ge0:A_t>u\}$ and  $W_u:=Y_{T_u}$.
Then $T_u<\infty$ and $A_{T_u}=u$.  
For each fixed $P\in\mathscr P$,  it is easy to verify that
$\langle W_{ij},W_{hl}\rangle_u=\langle Y_{ij},Y_{hl}\rangle_{T_u}=
\mathbf{1}_{\{(i,j)=(h,l)\}}A_{T_u}=\mathbf{1}_{\{(i,j)=(h,l)\}}u$.
Thus, $W$ is a $p\times N$-dimensional Brownian motion under each
$P$, with respect to the time-changed filtration $\mathcal F_{T_u}$.

For $p\ge N$, let
$ D:=\left\{z\in\mathbb R^{p\times N}:\operatorname{rank}(z)<N
\right\}$.
Since $D$ is contained in the zero set of a nonzero polynomial,
it has Lebesgue measure zero. Moreover, for every fixed $u>0$, the
random matrix $W_u$ has a nondegenerate Gaussian density. Consequently, $P(W_u\in D)=0$ for $u>0$.
Therefore,
$ E_P\left[
\int_0^\infty \mathbf 1_{ D}(W_u)\mathrm du
\right]=\int_0^\infty P(W_u\in D)\mathrm du=0$.
Hence,
$\int_0^\infty \mathbf 1_{ D}(W_u)\mathrm du=0,
\, P\text{-a.s.}$
Thus, we can obtain that 
$\int_0^r \mathbf 1_{ D}(Y_t)\mathrm dt=
\int_0^r \mathbf 1_{ D}(W_{A_t})\mathrm dt\le \frac{1}{\underline{\sigma}^2}
\int_0^{A_r}\mathbf 1_{ D}(W_u)\mathrm du
=0,\, P\text{-a.s.}$
for every $r>0$.
Consequently,
$\int_0^\infty \mathbf 1_{ D}(Y_t)\mathrm dt=0,
\, P\text{-a.s.}$ Since
$X_t^N=Y_t^{\mathsf T}Y_t$,
we have
$\operatorname{rank}(X_t^N)
=
\operatorname{rank}(Y_t)$.
It follows that
$\operatorname{rank}(X_t^N)=N,
\mathrm d t\times P\text{-a.e.}$
for every \(P\in\mathscr P\).


 Therefore, we get 
\begin{equation}\label{e:projection1}
\sqrt{X^{N}_t}\mathrm{d}B_t=\sqrt{X^{N}_t}\Big(\sqrt{X^{N}_t}\Big)^{-1} Y_t^{\top}\mathrm{d}Y_t=Y_t^{\top}\mathrm{d}Y_t,\quad\mathrm{d}B_t^{\top} \sqrt{X^{N}_t}=\mathrm{d}Y_t^{\top}Y_t.
\end{equation}
Noting that \eqref{e:dhatB} yields
$$\mathrm{d}B_{ki}=\sum_{s=1}^{N}\sum_{n=1}^{p}\Big( \sqrt{X^{N}_t}\Big)^{-1}_{ks}Y_{ns}\mathrm{d}Y_{ni}$$ and that  the underlying $G$-Brownian motion satisfies Assumption \ref{assum:G} with $d=p\times N$, we obtain
\begin{equation}\label{e:d<hatB>}
\begin{aligned}
    \mathrm{d}\langle B_{ij},B_{hl} \rangle_t&=\sum_{s=1}^{N}\sum_{n=1}^{p}\Big( \sqrt{X^{N}_t}\Big)^{-1}_{is}Y_{ns}\sum_{r=1}^{N}\sum_{q=1}^{p}\Big(\sqrt{X^{N}_t}\Big)^{-1}_{hr}Y_{qr}\mathrm{d}\langle Y_{nj}, Y_{ql}\rangle_t\\
&=\sum_{s=1}^{N}\sum_{r=1}^{N}\Big(\sum_{n=1}^{p}Y_{ns}Y_{nr}\Big)\Big( \sqrt{X^{N}_t}\Big)^{-1}_{hr}\Big( \sqrt{X^{N}_t}\Big)^{-1}_{is}\mathbf{1}_{\{j=l\}}\mathrm{d}\langle Y_{nj}, Y_{nl}\rangle_t\\
    &=\sum_{s=1}^{N}\sum_{r=1}^{N}\Big( \sqrt{X^{N}_t}\Big)^{-1}_{is}X_{sr}^{N}\Big( \sqrt{X^{N}_t}\Big)^{-1}_{hr}\mathbf{1}_{\{j=l\}}\mathrm{d}\langle Y_{1j}, Y_{1l}\rangle_t\\
    &=\mathbf{1}_{\{(i,j)=(h,l)\}}\mathrm{d}\langle Y_{11}\rangle_t.
\end{aligned}
\end{equation}

This together with \eqref{e:dX} implies the desired \eqref{G-N-wishart SDE},
and by Corollary \ref{corollary:G-brownian motion}, $B$ is an $N\times N$  $G$-Brownian matrix with $G$ satisfying Assumption~\ref{assum:G} with $d=N\times N$. This completes the proof.
\end{proof}

The following result shows that when $\alpha\ge N+1$,  the solution to the $G$-SDE~\eqref{e:G-Wishart} with positive  initial eigenvalues is positive definite. 

\begin{proposition}\label{theorem: wishart away from zero} Let $X_t^N$ be a solution to the $G$-SDE~\eqref  {e:G-Wishart} with  $x_0 \in \mathbb{S}_N^{++}$.  When $\alpha \ge N + 1$, we have $c\left( \exists t \in [0, T] : \det X^{N}_t \leq 0 \right) = 0$ for any $T>0$.  
\end{proposition}
\begin{proof}
Let  $\kappa_m := \inf \big\{ t \ge 0 : \det X^{N}_t \le \frac{1}{m} \, \text{or} \, \operatorname{Tr}(X^{N}_{t}) \ge m \big\} \wedge T.$ Then, to prove the desired result,  it suffices to prove $\lim_{m \to \infty} c(\kappa_m < T) = 0.$

Noting that
$$\mathrm{d}X^{N}_{ij}  =\sum_{k=1}^N\sqrt{X^{N}_{ik}} \mathrm{d} B_{kj}+\sum_{k=1}^N\mathrm{d} B_{ki} \sqrt{X^{N}_{kj}}+\alpha \mathbf{1}_{\{i=j\}}  \mathrm{d}\langle B_{11} \rangle_t,$$
we have
\begin{align*}
    \mathrm{d}\langle X^{N}_{ij}, X^{N}_{kl} \rangle_t 
    =& \sum_{m,n=1}^N \left\langle \sqrt{X^{N}_{im}} \mathrm{d}B_{mj} + \sqrt{X^{N}_{jm}} \mathrm{d}B_{mi} ,\, \sqrt{X^{N}_{kn}} \mathrm{d}B_{nl} + \sqrt{X^{N}_{ln}} \mathrm{d}B_{nk} \right\rangle_t\\
 = &\sum_{m,n=1}^N \bigg[\sqrt{X^{N}_{im}} \sqrt{X^{N}_{kn}} \mathrm{d}\langle B_{mj}, B_{nl} \rangle_t +  \sqrt{X^{N}_{im}} \sqrt{X^{N}_{ln}} \mathrm{d}\langle B_{mj}, B_{nk} \rangle_t \nonumber \\
    &\qquad +\sqrt{X^{N}_{jm}}\sqrt{X^{N}_{kn}} \mathrm{d}\langle B_{mi}, B_{nl} \rangle_t + \sqrt{X^{N}_{jm}} \sqrt{X^{N}_{ln}} \mathrm{d}\langle B_{mi}, B_{nk} \rangle_t\bigg] \nonumber \\
    =& \sum_{m=1}^N \bigg[\sqrt{X^{N}_{im}}\sqrt{X^{N}_{km}} \mathbf{1}_{\{j=l\}} + \sqrt{X^{N}_{im}}\sqrt{X^{N}_{lm}} \mathbf{1}_{\{j=k\}} \nonumber \\
    &\qquad + \sqrt{X^{N}_{jm}}\sqrt{X^{N}_{km}}\mathbf{1}_{\{i=l\}} + \sqrt{X^{N}_{jm}}\sqrt{X^{N}_{lm}} \mathbf{1}_{\{i=k\}} \bigg]\mathrm{d}\langle B_{11} \rangle_t\\
=&\big( X^{N}_{ik}\mathbf{1}_{\{j=l\}} + X^{N}_{il}\mathbf{1}_{\{j=k\}} + X^{N}_{jk}\mathbf{1}_{\{i=l\}} + X^{N}_{jl}\mathbf{1}_{\{i=k\}} \big) \mathrm{d}\langle B_{11} \rangle_t.
\end{align*}

Denote 
\begin{equation*}
    V(X^{N}):= \operatorname{Tr}(X^{N}) - \ln(\det X^{N})=\sum_{i=1}^N \left( \lambda^{N}_i - \ln \lambda^{N}_i \right).
\end{equation*}
We have
\begin{align*}
    \frac{\partial V}{\partial (X^{N})_{ij}} = \mathbf 1_{\{i=j\}} - (X^{N})^{-1}_{ji}, \quad 
    \frac{\partial^2 V}{\partial X^{N}_{ij} \partial X^{N}_{kl}} = (X^{N})^{-1}_{il} (X^{N})^{-1}_{kj}.
\end{align*}
 Applying  $G$-It\^{o}'s formula to $V(X^{N}_t)$, we get
\begin{align*} \label{Ito_expansion}
    \mathrm{d}V(X^{N}_t) =& \sum_{i,j=1}^N \frac{\partial V}{\partial X^{N}_{ij}}\big(\sum_{k=1}^N\sqrt{X^{N}_{ik}} \mathrm{d} B_{kj}+\sum_{k=1}^N\mathrm{d} B_{ki} \sqrt{X^{N}_{kj}}+\alpha \mathbf{1}_{\{i=j\}} \mathrm{d}\langle B_{ii} \rangle_t \big)\\
    &+ \frac{1}{2} \sum_{i,j,k,l=1}^N \frac{\partial^2 V}{\partial X^{N}_{ij} \partial X^{N}_{kl}} \big( X^{N}_{ik}\mathbf{1}_{\{j=l\}} + X^{N}_{il}\mathbf{1}_{\{j=k\}} + X^{N}_{jk}\mathbf{1}_{\{i=l\}} + X^{N}_{jl}\mathbf{1}_{\{i=k\}} \big) \mathrm{d}\langle B_{11} \rangle_t\\
    =&\sum_{i,j=1}^N \frac{\partial V}{\partial X^{N}_{ij}}\big(\sum_{k=1}^N\sqrt{X^{N}_{ik}} \mathrm{d} B_{kj}+\sum_{k=1}^N\mathrm{d} B_{ki} \sqrt{X^{N}_{kj}} \big)+ \alpha \sum_{i=1}^N \left( 1 - ((X^{N})^{-1})_{ii} \right) \mathrm{d}\langle B_{11} \rangle_t\\
    &+ \frac{1}{2} \sum_{i,j,k,l=1}^N(X^{N})^{-1}_{il} (X^{N})^{-1}_{kj} \big( X^{N}_{ik}\mathbf{1}_{\{j=l\}} + X^{N}_{il}\mathbf{1}_{\{j=k\}} + X^{N}_{jk}\mathbf{1}_{\{i=l\}} + X^{N}_{jl}\mathbf{1}_{\{i=k\}} \big) \mathrm{d}\langle B_{11} \rangle_t. 
    \end{align*}
The second-order term simplifies to 
\begin{align*}
&\frac{1}{2}\sum_{i,j,k=1}^N (X^{N})^{-1}_{ij} (X^{N})^{-1}_{kj} X^{N}_{ik} +\frac{1}{2}\sum_{i,j,l=1}^N (X^{N})^{-1}_{il} (X^{N})^{-1}_{jj}X^{N}_{il} \\
   &+\frac{1}{2}\sum_{i,j,k=1}^N (X^{N})^{-1}_{ii} (X^{N})^{-1}_{kj} X^{N}_{jk}+\frac{1}{2}
\sum_{i,j,l=1}^N (X^{N})^{-1}_{il} (X^{N})^{-1}_{ij} X^{N}_{jl}\\
 =&\frac{1}{2}\sum_{j,k=1}^N  (\mathbf{I}_N)_{jk} (X^{N})^{-1}_{kj}+\frac{1}{2}\operatorname{Tr}(\mathbf{I}_N) \operatorname{Tr}((X^{N})^{-1}) +\frac{1}{2}\operatorname{Tr}((X^{N})^{-1}) \operatorname{Tr}(\mathbf{I}_N)+\frac{1}{2}\sum_{i,l=1}^N (\mathbf{I}_N)_{il} (X^{N})^{-1}_{li}\\
    =    &(N+1) \operatorname{Tr}((X^{N})^{-1}). 
\end{align*}
Thus, we get 
    \begin{align*}
    \mathrm{d}V(X_t^N)=\sum_{i,j=1}^N \frac{\partial V}{\partial X^{N}_{ij}}\big(\sum_{k=1}^N\sqrt{X^{N}_{ik}} \mathrm{d} B_{kj}+\sum_{k=1}^N\mathrm{d} B_{ki} \sqrt{X^{N}_{kj}} \big)+ \left[ \alpha N -(\alpha - (N + 1))\operatorname{Tr}((X^{N})^{-1})  \right]\mathrm{d}\langle B_{11} \rangle_t.
\end{align*}

Note that  $X^{N}_t$ is positive definite for $t < \kappa_m$ and so is its inverse $(X^{N}_t)^{-1}$.  Thus, when $\alpha \ge N+1$, we have for $t < \kappa_m$,
$$\mathrm{d}V(X^{N}_t) \le \mathrm{d}M_t^V + \alpha N \, \mathrm{d}\langle B_{11} \rangle_t,$$
where $$\mathrm{d}M_t^V=\sum_{i,j=1}^N \frac{\partial V}{\partial X^{N}_{ij}}\big(\sum_{k=1}^N\sqrt{X^{N}_{ik}} \mathrm{d}B_{kj}+\sum_{k=1}^N\mathrm{d} B_{ki} \sqrt{X^{N}_{kj}} \big).$$
Thus, we have
\begin{equation}\label{e:EV}
 \mathbb{\hat E}\left[V(X_{\kappa_m}^N)\right]\le \mathbb{\hat E}\left[ V(X^{N}_0) \right] + \alpha N \overline{\sigma}^2 T.
\end{equation}

Now we study the behavior of $V(X_t^N)$ on the event $\{\kappa_m<T\}$, where either $\det X^{N}_t= \frac1m$ or $\operatorname{Tr}(X^{N}_t)= m$ for some $t=\kappa_m<T$. 
If  $\det X^{N}_t  = \frac{1}{m}$, then we have
$V(X^{N}_t) = \operatorname{Tr}(X^{N}_t) - \ln \det X^{N}_t  > \ln m$, since $X^{N}_t \in \mathcal{S}_N^{++}$. 
If $\operatorname{Tr}(X^{N}_t) = m$, noting that 
$\frac1N \ln \det X_t^N =\frac{1}{N}\sum_{i=1}^N \ln \lambda^{N}_i \le \ln\left( \frac{1}{N}\sum_{i=1}^N \lambda^{N}_i \right) = \ln \left( \frac{\operatorname{Tr}(X^{N}_t)}{N} \right)$, we have
$V(X^{N}_t) \ge  m -N\ln\frac mN$.

Thus, we have 
 $V(X^{N}_{\kappa_m}) \ge g(m)$ on the event $\{\kappa_m< T \}$,  where $g(m) := \min\{\ln m ,\, m - N\ln(m/N) \}$. Consequently, by~\eqref{e:EV} we get 
\begin{equation*}
    c(\kappa_m < T) \le c\big(  V(X^{N}_{\kappa_m}) \ge g(m) \big) \le \frac{\mathbb{\hat E}\left[ V(X^{N}_0) \right] + \alpha N \overline{\sigma}^2 T }{g(m)}.
\end{equation*}
Since  $\lim\limits_{m \to \infty} g(m) = \infty$, taking the limit yields 
$\lim\limits_{m \to \infty} c(\kappa_m < T) = 0$. This completes the proof.
\end{proof}

Now, we are ready to prove Theorem \ref{thm:wishart original existence}.
\begin{proof}[Proof of Theorem~\ref{thm:wishart original existence}]\makeatletter\def\@currentlabelname{the proof}\makeatother\label{proof:thm:wishart original existence}
For $m \ge 1$, define
\begin{equation*}
    D_m := \Big\{ Y^{N} \in \mathbb{S}_N^{++} : \det Y^{N} > \frac{1}{m}  \text{ and } \operatorname{Tr}(Y^{N}) < m \Big\}.
\end{equation*}
Let $\psi_m: \mathbb{S}_N \to [0, 1]$ be a Lipschitz continuous cutoff function satisfying:
\begin{equation*}
    \psi_m(Y^{N}) = 
    \begin{cases} 
    1, & Y^{N} \in D_{m}, \\
    0, & Y^{N} \notin D_{m+1}.
    \end{cases}
\end{equation*}
We define  $\Phi_m(Y^{N}) := \psi_m(Y^{N})\sqrt{Y^{N}}$. Clearly, $\Phi_m$ is globally Lipschitz continuous on $\mathbb{S}_N$. Then the following $G$-SDE
\begin{equation} \label{e:G-Wishart lipschitz}
    \mathrm{d}X^{N}_{(m)} = \Phi_m(X^{N}_{(m)}) \mathrm{d}B_t + \mathrm{d}B_t^{\top} \Phi_m(X^{N}_{(m)}) + \alpha \mathbf{I}_{N} \mathrm{d}\langle B_{11}\rangle_t, \quad X^{N}_{(m)}(0) = x_0,
\end{equation}
admits a unique solution $(X^{N}_{(m)}(t))_{t \in [0, T]}\in (M_G^2(0,T))^d$ with $d=N\times N$ by Theorem 5.1.3 in \cite{peng2019nonlinear}.  For each fixed 
$t\in [0,T]$, the integral representation of the solution to ~\eqref{e:G-Wishart lipschitz}, by the standard Picard estimate, yields that $X^{N}_{(m)}(t)\in (L_G^2(\Omega_t))^d$.

For \(R>0\), define
$\tau_R:=\inf\bigl\{s\ge 0:|X_{(m)}^N(s)|\ge R\bigr\}\wedge T,
\,X_{(m)}^{N,R}(s):=X_{(m)}^N(s\wedge\tau_R),
\, 0\le s\le T$.
Let
$Y_m(t):= \hat {\mathbb E}\!\left[  \sup_{0\leq r\leq t}\left|X_{(m)}^N(r)\right|^4\right],\,Y_{m,R}(t):=\hat {\mathbb E}\!\left[  \sup_{0\leq r\leq t}\left|X_{(m)}^{N,R}(r)\right|^4\right]$.
Since $X_{(m)}^N(s)\in\mathbb S_N^+$ and
$0\leq\psi_m(X_{(m)}^N(s))\leq 1$, we have
\begin{align*}
\left|\Phi_m\bigl(X_{(m)}^N(s)\bigr)\right|^2=
\psi_m\bigl(X_{(m)}^N(s)\bigr)^2
\left|\sqrt{X_{(m)}^N(s)}\right|^2 \leq \operatorname{Tr}\bigl(X_{(m)}^N(s)\bigr) \leq \sqrt N\,\left|X_{(m)}^N(s)\right|,
\end{align*}
where $|\cdot|$ denotes the
Frobenius norm defined by $|A|^2
:= \operatorname{Tr}(A^{\top}A)
= \sum_{i,j=1}^{N} |A_{ij}|^2$  for a matrix \(A=(A_{ij})\in\mathbb{R}^{N\times N}\).

Then
\begin{align*}
\left(
  \int_0^t \mathbf{1}_{\{s\le \tau_R\}}
  \left|\Phi_m\bigl(X_{(m)}^N(s)\bigr)\right|^2\,\mathrm d\langle B\rangle_s
\right)^2
&\leq
C_N {\overline{\sigma}^4}\left(\int_0^t\mathbf{1}_{\{s\le \tau_R\}}\left|X_{(m)}^N(s)\right|\,\mathrm ds\right)^2 \\
&\leq
C_N {\overline{\sigma}^4}t\int_0^t\left|X_{(m)}^{N,R}(s)\right|^2\,\mathrm ds \leq
C_N {\overline{\sigma}^4}t^2\sup_{0\leq r\leq t}
       \left|X_{(m)}^{N,R}(r)\right|^2.
\end{align*}
Thus,
\begin{align*}
\hat{\mathbb E}\!\left[
  \left(
    \int_0^t \mathbf{1}_{\{s\le \tau_R\}}
    \left|\Phi_m\bigl(X_{(m)}^{N}(s)\bigr)\right|^2\,\mathrm d\langle B\rangle_s
  \right)^2
\right]\leq C_N {\overline{\sigma}^4}t^2\,
\hat{\mathbb E}\!\left[
  \sup_{0\leq r\leq t}
  \left|X_{(m)}^{N,R}(r)\right|^2
\right]\leq C_N {\overline{\sigma}^4}t^2\,Y_{m,R}(t)^{1/2}.
\end{align*}

By the $G$-Burkholder-Davis-Gundy inequality, we can obtain that
\begin{align*}
Y_{m,R}(t)
&\leq C_0|x_0|^4
+C_0\widehat{\mathbb E}\!\left[
  \left(
    \int_0^t\mathbf{1}_{\{s\le \tau_R\}}
    \left|\Phi_m\bigl(X_{(m)}^{N}(s)\bigr)\right|^2\,\mathrm d\langle B\rangle_s
  \right)^2
\right]+C_0\alpha^4N^2{\overline{\sigma}^8}t^4\\
&\leq C_0|x_0|^4+C_0C_N{\overline{\sigma}^4}t^2\,Y_{m,R}(t)^{1/2}+C_0\alpha^4N^2{\overline{\sigma}^8}t^4\\
&\leq C_0|x_0|^4+\frac12Y_{m,R}(t)+C_1t^4,
\end{align*}
where $C_0$ and $C_1$ are independent of $m$. Consequently,
we get the uniform estimate
$$\sup_{m\geq1}
\hat{\mathbb E}\!\left[
  \sup_{0\leq r\leq T}
  \left|X_{(m)}^{N,R}(r)\right|^4
\right]
\leq
2C_0|x_0|^4+2C_1T^4
\leq C_T(1+|x_0|^4),$$
where $C_T$ is independent of $m$.

Noting that, by Fatou's lemma, for any $P\in\mathscr P$,
$$E_P\left[
  \sup_{0\leq r\leq T}
  \left|X_{(m)}^{N}(r)\right|^4
\right]\le \liminf_{R\to \infty}E_P\left[
  \sup_{0\leq r\leq T}
  \left|X_{(m)}^{N,R}(r)\right|^4
\right]\le \liminf_{R\to \infty}\hat{\mathbb E}\!\left[
  \sup_{0\leq r\leq T}
  \left|X_{(m)}^{N,R}(r)\right|^4
\right],$$
we get for all $m\in\mathbb N$,
$$\hat{\mathbb E}\left[
  \sup_{0\leq r\leq T}
  \left|X_{(m)}^{N}(r)\right|^4
\right]\le C_T(1+|x_0|^4).$$

Denote  $\kappa^{(m)}_{m}:=\inf\{t\ge0:X^{N}_{(m)}(t)\notin D_m\}\wedge T$. Then, $X^{N}_{(m)}(t)=X^{N}_{(m+k)}(t), t\in[0,\kappa_m^{(m)}]$ for all $k\in\mathbb N$ by the uniqueness of the solution, and hence  $\kappa_m^{(m)}$ is increasing in $m$.  Set $X^{N}_t := \lim_{m\to\infty}X^{N}_{(m)}(t)$. By a similar argument as in the proof of Proposition~\ref{theorem: wishart away from zero}, we can prove that $\lim\limits_{m \to \infty} \kappa^{(m)}_m = T$ q.s., and hence $X^{N}_t$  is defined on the total interval $[0,T]$. 
Then,   we get 
\begin{align*}
    \hat{\mathbb E}\left[
  \left|X_{(m)}^{N}(t)-X^{N}(t)\right|^2
\right]&\le 2\Big[\big(\hat{\mathbb E}\big[
  \big|X_{(m)}^{N}(t)\big|^4
\big]\big)^{\frac 1 2}+\big(\hat{\mathbb E}\big[
  \big|X^{N}(t)\big|^4
\big]\big)^{\frac 1 2}\Big](c(\kappa^{(m)}_m<T))^{\frac{1}{2}}\\
&\le 4 C_T(1+|x_0|^4)(c(\kappa^{(m)}_m<T))^{\frac{1}{2}}.\end{align*}
This implies that $X_{(m)}^N(t)$ converges to $X^N(t)$ in $L^2$. 

Finally, the uniqueness of the solution follows directly from the uniqueness of the solution on each $[0, \kappa_m^{(m)}]$ and the fact that $\kappa_m^{(m)}\to T$ q.s.\ as $m\to\infty$. 
\end{proof}

\begin{proof}[Proof of Proposition~\ref{thm: wishart initial zero}]\makeatletter\def\@currentlabelname{the proof}\makeatother\label{proof:thm: wishart initial zero}
As in the proof of Lemma \ref{lem:G-Wishart},  define
$A_t:=\langle  B^{11}\rangle _t,\,
T_u:=\inf\{t\ge0:A_t>u\},\,
W_u:=B_{T_u}$ and $Z_u:=X_{T_u}$.
Under each $P\in\mathscr P$, $W$ is an $N\times N$ standard
Brownian matrix and
\begin{equation}\label{eq:classical wishart}
\mathrm{d}Z_u=\sqrt{Z_u}\,\mathrm{d}W_u+\mathrm{d}W_u^{\top}\sqrt{Z_u}+ p I_N\,\mathrm{d}u,\,
   Z_0=x_0.
\end{equation}

Let $0\leq \kappa_1(u)\leq \cdots \leq \kappa_N(u)$
be the ordered eigenvalues of $Z_u$. For each fixed $u>0$, the law of $Z_u$ is absolutely continuous
with respect to Lebesgue measure  on $\mathbb S_N$, and we have $P(0<\kappa_1(u)<\cdots<\kappa_N(u))=1$. By \cite{Bru1991}, for each $u>0$, we have $P(E_u)=1$ where $E_u:=\{0<\kappa_1(t)<\kappa_2(t)<\cdots<\kappa_N(t), t\ge u\}$.  
Let $\Omega_*:=\bigcap_{r\in\mathbb Q_+}E_r$ and we have $P(\Omega_*)=1$. Since every $t>0$ exceeds some
$r\in\mathbb Q_+$, on $\Omega_*$ we have
$0<\kappa_1(t)<\cdots<\kappa_N(t)$ 
for all $t>0$.
Since $P\in\mathscr P$ is arbitrary, the ordered eigenvalues of
$X^N$ satisfy
$0<\lambda_1^N(t)<\cdots<\lambda_N^N(t)$
 for all $t>0$ q.s. Then, by Theorem~\ref{thm:wishart non-collision}, the differential equation~\eqref{eq:wishart eigenvalue zero} holds q.s.\ for $t>0$. We still need to show that equation~\eqref{eq:wishart eigenvalue zero} starting at  time 0 holds, for which it suffices to prove the existence of the integral of the equation on $[0,T]$.


Set
$R_u:=\operatorname{tr}(Z_u)$ and $R_0:=\operatorname{tr}(x_0)$.
Taking traces in \eqref{eq:classical wishart} yields
$R_u=R_0+M_u+pNu$,
where
$M_u
:=2\int_0^u\operatorname{tr}\bigl(\sqrt{Z_v}\,\mathrm dW_v\bigr)$
with
$\mathrm d\langle M\rangle_u
=4R_u\,\mathrm du$.
Denoting    $\tau_n:=\inf\{u\ge0:R_u\ge n\}$ 
, we have $E_P[R_{u\wedge\tau_n}]\le R_0+pNu$.
Thus, by Fatou's lemma, $E_P[R_u]\le R_0+pNu$. 
Moreover, It\^o's isometry yields, for $0\le s\le t\le U$,
\begin{align*}
E_P\Big[\|Z_t-Z_s\|_F^2\Big]
&\le 8E_P\Big[\big\|\int_s^t\sqrt{Z_v}\,\mathrm dW_v\big\|_F^2\Big]
      +2p^2N(t-s)^2 \\
      &=8N\int_s^t E_P[R_v]\,\mathrm dv
      +2p^2N(t-s)^2\le C_U(t-s).
\end{align*}
where $C_U$ depends only on $U,N,p$, and $x_0$.

For $0\le x_1<\cdots<x_N$, set
$F_i(x):=\sum_{k:k\ne i}\frac{x_i+x_k}{x_i-x_k},
\,
V(x):=\sum_{1\le i<j\le N}
(j-i)\frac{x_i+x_j}{x_j-x_i}$.

Then $|F_i(x)|\le V(x)$ and
$\sum_{j=1}^N jF_j(x)=V(x)$.
By $G$-It\^o's formula, for $t\ge r>0$,
$$\lambda^N_i(t)-\lambda^N_i(r)=
2\int_r^t\sqrt{\lambda^N_i(s)}\,\mathrm d\bar B_s^i
+p(\langle \bar B^i\rangle _t-\langle \bar B^i\rangle _r)
+\int_r^t F_i(\lambda^N(s))\,\mathrm d \langle \bar B^i\rangle _s,$$
where
$\mathrm d\bar B^i_s
=\frac{1}{2}
\sum_{k,h=1}^N
\bigl(P_i(X^N_s)\bigr)_{kh}\,(\mathrm dB_s^{kh}+\mathrm dB_s^{hk})$ and $\mathrm d \langle \bar B^i,\bar B^j\rangle _s=\mathbf{1}_{\{i=j\}}\mathrm d \langle  B^{11}\rangle _s=\mathbf{1}_{\{i=j\}}\mathrm dA_s$.
 Here $P_i(X_s^N)$ denotes the orthogonal spectral projection onto the
eigenspace of $\lambda^N_i(s)$ defined by
$$P_i(X_s^N)=\prod_{j\ne i}\frac{X_s^N-\lambda^N_j(X_s^N)I_N}{\lambda^N_i(X_s^N)-\lambda^N_j(X_s^N)},\,s>0.$$
Since
$P_i(X^N)$ is a smooth function of $X^N$ on $[r,t]$ and 
 $
\hat{\mathbb E}\!\left[
\int_0^r \|P_i(X_s^N)\|_F^2\,\mathrm ds
\right]=r$, letting $r\to 0$, $\int_{r}^t\mathrm d\bar B^i_s$ can be extended to $[0,t]$. By Corollary \ref{corollary:G-brownian motion}, $\bar B=(\bar B^1, \dots, \bar B^N)$ is an  $N$-dimensional  $\bar{G}$-Brownian motion with $\bar{G}$ satisfying Assumption~\ref{assum:G} with $d=N$.

Denote
$S(x):=\sum_{j=1}^N jx_j$.
Applying It\^o's formula gives
\[
S(\kappa(v))-S(\kappa(r))
=
\tilde{M}_v-\tilde{M}_r+p\frac{N(N+1)}{2}(v-r)
+\int_r^vV(\kappa(u))\,\mathrm du,
\]
where \(\tilde{M}\) is the martingale part.
Therefore,
\begin{align*}
E_{P}\Big[\int_r^v V(\kappa(u))\,\mathrm du\Big]
&\le E_{P}\bigl[S(\kappa(v))-S(\kappa(r))\bigr]\le\big(\sum_{j=1}^N j^2\big)^{1/2}E_{P}\bigl[|\kappa(v)-\kappa(r)|\bigr]\\
&\leq\big(\sum_{j=1}^N j^2\big)^{1/2}\big(E_{P}\bigl[\|Z_v-Z_r\|^2_{\mathrm F}\bigr]\big)^{1/2}\le C\sqrt v.
\end{align*}
Letting $r\downarrow0$ gives
$0\le E_{P}\Big[\int_0^v V(\kappa(u))\,\mathrm du\Big]\le C\sqrt v$.

For each $P\in\mathscr P$, the time-change formula gives $P$-a.s.
\[
\int_0^r |F_i(\lambda^N(s))|\,dA_s
=
\int_0^{A_r} |F_i(\kappa(u))|\,du
\le
\int_0^{\bar{\sigma}^2r} V(\kappa(u))\,du.
\]
Then
\[
\hat{\mathbb E}\!\left[
  \int_0^r |F_i(\lambda^N(s))|\,\mathrm d \langle \bar B^i\rangle _s
\right]
\le
\sup_{P\in\mathscr P} E_P\!\left[\int_0^{\bar{\sigma}^2r} V(\kappa(u))\,\mathrm d u\right]
\le C\sqrt{\bar{\sigma}^2r}.
\]
And we can obtain that 
\[
\begin{aligned}
&\hat{\mathbb E}\Big[|2\int_0^r
\sqrt{\lambda^N_i(s)}\,\mathrm d\bar B_s^i|^2\Big]
\le 4\,\hat{\mathbb E}
\Big[\int_0^r\lambda^N_i(s)\,\mathrm d \langle \bar B^i\rangle _s\Big]\le 4\,
\hat{\mathbb E}\Big[\int_0^r\operatorname{tr}(X^N_s)\,\mathrm d A_s\Big]\\
&\le 4\,
\sup_{P\in\mathscr P} E_P
\Big[\int_0^{ {\overline \sigma}^2 r}
\operatorname{tr}(Z_u)\,\mathrm du\Big]
\le 4\int_0^{ {\overline \sigma}^2 r}(R_0+pNu) \mathrm du=4R_0\bar{\sigma}^2 r+2pN\bar{\sigma}^4 r^2.
\end{aligned}
\]
Letting $r\to 0$ in the identity on $[r,t]$ yields ~\eqref{eq:wishart eigenvalue zero}.
\end{proof}

\section{High-dimensional spectral limits}\label{sec:hign limits}

In this section, we investigate the high-dimensional limits of the empirical spectral measures of generalized Wishart processes, providing proofs for Theorem~\ref{thm: pre-compact}, Theorem~\ref{thm:limit-equation},  Proposition~\ref{dyson converge}, and Proposition~\ref{wishart converge}.

\subsection{Generalized Wishart Processes}
\label{subsection:Generalized Wishart Processes}
First, we prove the tightness of the empirical spectral measures.
\begin{proof}[Proof of Theorem~\ref{thm: pre-compact}]\makeatletter\def\@currentlabelname{the proof}\makeatother\label{proof:thm: pre-compact}
Let  $f\in C^2(\mathbb{R})$ be fixed.   For each
$t\in[0,T]$, we have
$$\left\langle f, L_{N}(t)\right\rangle=\int_\R f(x) L_{N}(t)(\mathrm{d}x)=\frac{1}{N} \sum_{i=1}^{N} \int_\R f(x) \delta_{\lambda_{i}^{N}(t)}(\mathrm{d}x)=\frac{1}{N} \sum_{i=1}^{N} f\left(\lambda_{i}^{N}(t)\right).$$
By $G$-Itô's formula, we get for $i=1,\ldots,N$,
\begin{align*}
   & f\left(\lambda_{i}^{N}(t)\right)-f\left(\lambda_{i}^{N}(0)\right) \\
= & \int_{0}^{t} f^{\prime}\left(\lambda_{i}^{N}(s)\right) 2 g_{N}\left(\lambda_{i}^{N}\right) h_{N}\left(\lambda_{i}^{N}\right) \mathrm{d} \bar B_i+\int_{0}^{t} f^{\prime}\left(\lambda_{i}^{N}(s)\right)b_{N}\left(\lambda_{i}^{N}\right)  \mathrm{d} \langle \bar B_i  \rangle _{s}\\
&+\int_{0}^{t} f^{\prime}\left(\lambda_{i}^{N}(s)\right) \sum_{k:k \neq i} \frac{H_{N}\left(\lambda_{i}^{N}, \lambda_{k}^{N}\right)}{\lambda_{i}^{N}-\lambda_{k}^{N}} \mathrm{d} \langle \bar B_i  \rangle _{s}+ 2\int_{0}^{t} f^{\prime \prime}\left(\lambda_{i}^{N}(s)\right)g_{N}\left(\lambda_{i}^{N}\right)^{2} h_{N}\left(\lambda_{i}^{N}\right)^{2}\mathrm{d} \langle \bar B_i  \rangle _{s}.
\end{align*}
Summing over $i$ and dividing by $N$, we obtain 
\begin{align}\label{eq:pre-symmetrization}
&\langle f,  \left.L_{N}(t)\right\rangle-\left\langle f, L_{N}(0)\right\rangle\notag\\
=& \frac{1}{N} \sum_{i=1}^{N} \int_{0}^{t} f^{\prime}\left(\lambda_{i}^{N}(s)\right) 2 g_{N}\left(\lambda_{i}^{N}\right) h_{N}\left(\lambda_{i}^{N}\right) \mathrm{d} \bar B_i+\frac{1}{N} \sum_{i=1}^{N}\int_{0}^{t} f^{\prime}\left(\lambda_{i}^{N}(s)\right)b_{N}\left(\lambda_{i}^{N}\right) \mathrm{d} \langle \bar B_i  \rangle _{s}\notag\\
&+\frac{1}{N}\sum_{k \neq i}  \int_{0}^{t} f^{\prime}\left(\lambda_{i}^{N}(s)\right)  \frac{H_{N}\left(\lambda_{i}^{N}, \lambda_{k}^{N}\right)}{\lambda_{i}^{N}-\lambda_{k}^{N}} \mathrm{d} \langle \bar B_i  \rangle _{s}+\frac{2}{N} \sum_{i=1}^{N}\int_{0}^{t} f^{\prime \prime}\left(\lambda_{i}^{N}(s)\right)g_{N}\left(\lambda_{i}^{N}\right)^{2} h_{N}\left(\lambda_{i}^{N}\right)^{2}\mathrm{d} \langle \bar B_i  \rangle _{s}\notag\\
=& M_{f}^{N}(t) +\int_{0}^{t}\left\langle f^{\prime} b_{N}, L_{N}(s)\right\rangle\mathrm{d} \langle \bar B_i  \rangle _{s}+2 \int_{0}^{t}\left\langle f^{\prime \prime} g_{N}^{2} h_{N}^{2}, L_{N}(s)\right\rangle \mathrm{d} \langle \bar B_i  \rangle _{s}\notag\\ &+\frac{1}{N} \sum_{k \neq i} \int_{0}^{t} f^{\prime}\left(\lambda_{i}^{N}(s)\right) \frac{H_{N}\left(\lambda_{i}^{N}, \lambda_{k}^{N}\right)}{\lambda_{i}^{N}-\lambda_{k}^{N}} \mathrm{d} \langle \bar B_i  \rangle _{s},
\end{align} 
where we denote
$$M_{f}^{N}(t):=\frac{1}{N} \sum_{i=1}^{N} \int_{0}^{t} f^{\prime}\left(\lambda_{i}^{N}(s)\right) 2 g_{N}\left(\lambda_{i}^{N}\right) h_{N}\left(\lambda_{i}^{N}\right)  \mathrm{d}\bar B_i.$$
We adopt the convention
$\frac{f'(x)-f'(y)}{x-y}=f''(x),
\,\text{whenever }x=y$.
Using the symmetry of $H_N$, the interaction term in
\eqref{eq:pre-symmetrization} can be symmetrized as follows:
\begin{align*}
&\frac{1}{N} \sum_{k \neq i} \int_{0}^{t} f^{\prime}\left(\lambda_{i}^{N}(s)\right) \frac{H_{N}\left(\lambda_{i}^{N}, \lambda_{k}^{N}\right)}{\lambda_{i}^{N}-\lambda_{k}^{N}} \mathrm{d} \langle \bar B_i  \rangle _{s}\\
=&\frac{1}{2N} \sum_{k \neq i} \int_{0}^{t} \frac{f^{\prime}\left(\lambda_{i}^{N}(s)\right)-f^{\prime}\left(\lambda_{k}^{N}(s)\right)}{\lambda_{i}^{N}(s)-\lambda_{k}^{N}(s)} H_{N}\left(\lambda_{i}^{N}, \lambda_{k}^{N}\right)\mathrm{d} \langle \bar B_i  \rangle _{s} \\
= & \frac{1}{2N} \sum_{k \neq i} \int_{0}^{t} \iint \frac{f^{\prime}(x)-f^{\prime}(y)}{x-y}H_{N}\left(x,y\right) \delta_{\lambda_{i}^{N}(s)}(\mathrm{d}x) \delta_{\lambda_{k}^{N}(s)}(\mathrm{d}y)\mathrm{d} \langle \bar B_i  \rangle _{s} \\ 
= &\frac{N}{2} \int_{0}^{t} \iint \frac{f^{\prime}(x)-f^{\prime}(y)}{x-y}H_{N}\left(x,y\right) L_{N}(s)(\mathrm{d}x) L_{N}(s)(\mathrm{d}y)  \mathrm{d} \langle \bar B_i  \rangle _{s}\\
&-\frac{1}{2N} \sum_{i=1}^{N} \int_{0}^{t} \iint \frac{f^{\prime}(x)-f^{\prime}(y)}{x-y}H_{N}\left(x,y\right) \delta_{\lambda_{i}^{N}(s)}(\mathrm{d}x) \delta_{\lambda_{i}^{N}(s)}(\mathrm{d}y) \mathrm{d} \langle \bar B_i  \rangle _{s}\\
=&\frac{N}{2} \int_{0}^{t} \iint \frac{f^{\prime}(x)-f^{\prime}(y)}{x-y}H_{N}\left(x,y\right) L_{N}(s)(\mathrm{d}x) L_{N}(s)(\mathrm{d}y)  \mathrm{d} \langle \bar B_i  \rangle _{s}\\
&-\int_{0}^{t}\left\langle f^{\prime \prime} g_{N}^{2} h_{N}^{2}, L_{N}(s)\right\rangle \mathrm{d} \langle \bar B_i  \rangle _{s}.
\end{align*} 
Substituting this identity into \eqref{eq:pre-symmetrization} yields
\begin{equation}\label{eq:f,L}
\begin{aligned}
\left\langle f, L_{N}(t)\right\rangle  =&\left\langle f, L_{N}(0)\right\rangle + \frac{N}{2} \int_{0}^{t} \iint \frac{f^{\prime}(x)-f^{\prime}(y)}{x-y}H_{N}\left(x,y\right) L_{N}(s)(\mathrm{d}x) L_{N}(s)(\mathrm{d}y)  \mathrm{d} \langle \bar B_i  \rangle _{s}\\
&+\int_{0}^{t}\left\langle f^{\prime} b_{N}, L_{N}(s)\right\rangle \mathrm{d} \langle \bar B_i  \rangle _{s}+\int_{0}^{t}\left\langle f^{\prime \prime} g_{N}^{2} h_{N}^{2}, L_{N}(s)\right\rangle \mathrm{d} \langle \bar B_i  \rangle _{s}+M_{f}^{N}(t).
\end{aligned} 
\end{equation}

The quadratic variation of  $M_{f}^{N}(t)$  satisfies the following estimate:
\begin{equation}\label{e:<MfN>}
\left\langle M_{f}^{N}\right\rangle_{t}=\frac{4}{N} \int_{0}^{t}\langle\left(f^{\prime }\right)^{2}g_{N}^{2} h_{N}^{2}, L_{N}(s)\rangle \mathrm{d} \langle \bar B_i  \rangle _{s} \leq \frac{4{\overline{\sigma}}^{2} t\left\|f^{\prime }g_{N}h_{N}\right\|_{\infty}^{2}}{N}.
\end{equation}
By Assumption \ref{assumption:generalized-wishart}, the drift and interaction terms in
\eqref{eq:f,L}, with $f=\varphi$, are
bounded uniformly with respect to $N$. 
Let $C$ be a constant such that
$$\|\varphi'b_N\|_{\infty}
+
\|\varphi''g_N^2h_N^2\|_{\infty}
+
\left\|
N H_N(x,y)\frac{\varphi'(x)-\varphi'(y)}{x-y}
\right\|_{\infty}
\le C,$$
uniformly in $N\ge1$. Hence, by~\eqref{eq:f,L} there exists a constant
$C_1=C_1(T,\bar{\sigma}^2,C)$ such that
\begin{equation}\label{e:bound-varphi-L}
\sup_{0\leq t\leq T}
\langle\varphi,L_N(t)\rangle
\leq
C_0+C_1T+
\sup_{0\leq t\leq T}|M_\varphi^N(t)|.
\end{equation}

By \eqref{e:<MfN>} we get
\begin{align*}
c\Big(
\sup_{0\leq t\leq T}|M_\varphi^N(t)|\geq\varepsilon
\Big)
\leq
\frac{C_{l_1}}{\varepsilon^{2l_1}}
\hat{\mathbb{E}}
\big[
\langle M_\varphi^N\rangle_T^{l_1}
\big]
\leq
\frac{C_{l_1,T,\bar{\sigma}}}{\varepsilon^{2l_1}}
\Big(
\frac{
\|\varphi'g_Nh_N\|_{\infty}^2
}{N}
\Big)^{l_1},
\end{align*}
where the right-hand side is summable in $N$ by Assumption \ref{assumption:generalized-wishart}. 
Then, Lemma \ref{lem:B-C} together with the estimate~\eqref{e:bound-varphi-L} implies
that there exists a finite constant $M$ such that $\sup_{0\leq t\leq T}
\langle\varphi,L_N(t)\rangle\leq M$,
for all sufficiently large $N$. Thus,
\begin{equation}
c\Big(
\bigcap_{N\geq1}
\bigcup_{N_0\geq N}
\left\{
L_{N_0}(t)\in K_M
\text{ for every }t\in[0,T]
\right\}^c
\Big)=0,
\label{eq:compact-containment}
\end{equation}
where 
$$K_M:=
\left\{
\mu\in \mathbf{P}(\mathbb{R}):
\langle \varphi, \mu\rangle\leq M
\right\}$$ is a compact subset in $\mathbf{P}(\mathbb{R})$ by Lemma~B.2 in~\cite{song2024eigenvalue}.

For $f\in C^2(\mathbb{R})$, by \eqref{eq:f,L} we have  
\begin{equation*}
\left|
\langle f,L_N(t)\rangle-\langle f,L_N(s)\rangle
\right|
\leq
C_2\bar{\sigma}^2|t-s|
+
|M_f^N(t)-M_f^N(s)|,
\end{equation*}
where \begin{equation*}
 C_2:=\sup_{N>0}\Big\{\|f'b_N\|_{\infty}
+
\|f''g_N^2h_N^2\|_{\infty}
+
\frac{1}{2}\Big\|
N H_N(x,y)\frac{f'(x)-f'(y)}{x-y}
\Big\|_{\infty}\Big\},  
\end{equation*}
which is finite for each $\tilde f_k$ by assumption.
Thus,  for some  $0<\delta\le 1$ and $M>C_2\bar{\sigma}^2$, we have
\begin{equation}\label{e:est-f-LN}
\begin{aligned}
&c\Big(
\sup_{\substack{s,t\in[0,T]\\|t-s|\leq\delta}}
\left|
\langle f,L_N(t)\rangle-\langle f,L_N(s)\rangle
\right|
\geq M\delta^{1/8}
\Big)\\
&\leq c\Big(
\sup_{\substack{s,t\in[0,T]\\|t-s|\leq\delta}}
|M_f^N(t)-M_f^N(s)|
\geq (M-C_2\bar{\sigma}^2)\delta^{1/8}
\Big). 
\end{aligned}
\end{equation}

Let
$J:=\lceil T\delta^{-1}\rceil$ for some fixed $\delta\in(0,1)$. Set
$I_k:=[k\delta,(k+1)\delta]\cap[0,T]$, $k=0,\ldots,J-1$.  For $0\le t<s\le T$ with $|t-s|\le \delta$, assuming $t\in I_k$ for some $k$, we have either,  if $s$ also lies in $I_k$, 
\begin{align*}|M_f^N(t)-M_f^N(s)|&\le |M_f^N(t)-M_f^N(k\delta)|+|M_f^N(k\delta)-M_f^N(s)|\\
&\le 2\sup_{t\in I_k}
|M_f^N(t)-M_f^N(k\delta)|,
\end{align*}
 or, if $s\in I_{k+1}$,
\begin{align*}|M_f^N(t)-M_f^N(s)|&\le |M_f^N(t)-M_f^N(k\delta)|+|M_f^N(k\delta)-M_f^N((k+1)\delta)|+|M_f^N(s)-M_f^N((k+1)\delta)|\\
&\le 2\sup_{t\in I_k}
|M_f^N(t)-M_f^N(k\delta)|+\sup_{t\in I_{k+1}}
|M_f^N(t)-M_f^N((k+1)\delta)|.
\end{align*}
Thus, we get 
\begin{align*}c\Big(
\sup_{\substack{s,t\in[0,T]\\|t-s|\leq\delta}}
|M_f^N(t)-M_f^N(s)|
\geq (M-C_2\bar{\sigma}^2)\delta^{1/8}
\Big) \le \sum_{k=0}^{J-1}
c\Big(
\sup_{t\in I_k}
|M_f^N(t)-M_f^N(k\delta)|
\geq\frac{(M-C_2\bar{\sigma}^2)}{3}\delta^{1/8}
\Big).
\end{align*}
 Hence, by \eqref{e:est-f-LN} and the
$G$-Burkholder--Davis--Gundy inequality, we get
\begin{equation}\label{eq:modulus-estimate}
c\Big(
\sup_{\substack{s,t\in[0,T]\\|t-s|\leq\delta}}
\left|
\langle f,L_N(t)\rangle-\langle f,L_N(s)\rangle
\right|
\geq M\delta^{1/8}
\Big)\leq
\frac{
C_{T,\bar{\sigma},l_2}\,
\delta^{\,3l_2/4-1}
}{
(M-C_2\bar{\sigma}^2)^{2l_2}N^{l_2}
}
\|f'g_Nh_N\|_{\infty}^{2l_2}.
\end{equation}

For each $f\in C^2(\mathbb{R})$
and each $\varepsilon>0$, define
$$C_T(f,\varepsilon)
:=
\bigcap_{n=1}^{\infty}
\Big\{
\mu\in C([0,T],\mathbf{P}(\mathbb{R})):
\sup_{\substack{s,t\in[0,T]\\|t-s|\leq n^{-4}}}
|\mu_t(f)-\mu_s(f)|
\leq\frac{1}{\varepsilon\sqrt n}
\Big\}.$$

Let $(\varepsilon_k)_{k\in\mathbb N}$ be a sequence of positive numbers converging to 0 satisfying
$\varepsilon_k^{-1}>C_k\bar{\sigma}^2+1+k\psi(k)^{\frac{1}{2l_2}}$, $\delta=n^{-4}$ and $M=\varepsilon_k^{-1}$.
By \eqref{eq:modulus-estimate} we get
\begin{equation}\label{eq:CT-estimate'}
\begin{aligned}
    c\big(L_N\notin C_T(\tilde f_k,\varepsilon_k)\big)
&\leq \sum_{n=1}^{\infty}c\Big(\sup_{\substack{s,t\in[0,T]\\|t-s|\leq n^{-4}}}
\left|
\langle \tilde f_k,L_N(t)\rangle-\langle \tilde f_k,L_N(s)\rangle
\right|
\geq\frac{1}{\varepsilon_k\sqrt n} \Big)\\
&\le \sum_{n=1}^{\infty}\frac{
C_{T,\bar{\sigma},l_2}\,
n^{-3l_2+4}
}{
(\varepsilon_k^{-1}-C_k\bar{\sigma}^2)^{2l_2}N^{l_2}
}
\|\tilde f_k'g_Nh_N\|_{\infty}^{2l_2}.\\
&=\frac{
C_{T,\bar{\sigma},l_2}\,\|\tilde f_k'g_Nh_N\|_{\infty}^{2l_2}
}{
(\varepsilon_k^{-1}-C_k\bar{\sigma}^2)^{2l_2}N^{l_2}
}\sum_{n=1}^{\infty}n^{-3l_2+4}.
\end{aligned}
\end{equation}
Note that $l_2\geq2$, and hence 
$\sum_{n=1}^{\infty}n^{-3l_2+4}<\infty.$ 

Then by Assumption \ref{assumption:generalized-wishart} and~\eqref{eq:CT-estimate'}, we get
\begin{align}\label{eq:CT-estimate}    \sum_{k\ge1}\sum_{N=1}^{\infty}c(L_N\notin C_T(\tilde f_k,\varepsilon_k))&\le\sum_{k\ge1}\sum_{N=1}^{\infty}\frac{
C\,\|\tilde f_k'g_Nh_N\|_{\infty}^{2l_2}
}{
(\varepsilon_k^{-1}-C_k\bar{\sigma}^2)^{2l_2}N^{l_2}
}\le \sum_{k\ge1}\frac{C}{k^{2l_2}}<\infty.
\end{align}

Denote  $\mathcal{K}:=\{\mu\in C([0,T],\mathbf{P}(\mathbb{R})):\mu_t\in K_M,\forall t\in[0,T]\}\cap\bigcap_{k\geq1}
C_T(\widetilde f_k,\varepsilon_k)$, which is  a compact subset of $C([0,T],\mathbf{P}(\mathbb{R}))$ by \cite[Lemma~4.3.13]{agz10}. Combining \eqref{eq:compact-containment}, \eqref{eq:CT-estimate}, and
Lemma \ref{lem:B-C}, we conclude that
$$c\big(
\bigcap_{N\geq1}
\bigcup_{N_0\geq N}
\{L_{N_0}\in\mathcal{K}\}^c
\big)=0.$$
Therefore, the sequence
$\{L_N\}_{N\geq1}$ is relatively compact in
$C([0,T],\mathbf{P}(\mathbb{R}))$, quasi-surely.
\end{proof}

Now we derive the equation satisfied by the spectral limit. 
\begin{proof}[Proof of Theorem~\ref{thm:limit-equation}]\makeatletter\def\@currentlabelname{the proof}\makeatother\label{proof:thm:limit-equation}
For any limit point $\mu$, we can find a 
subsequence $\left\{N_i\right\}$, such that $L_{N_i}$ converges to $\mu$ in $C([0,T],\mathbf{P}(\mathbb{R}))$ as $N_i$ tends to infinity.

Fix \(z\in\mathbb{C}\setminus\mathbb{R}\) and define
$f_z(x):=\frac{1}{z-x},
\, x\in\mathbb{R}$.
Then
$f_z'(x)=\frac{1}{(z-x)^2},
\,
f_z''(x)=\frac{2}{(z-x)^3}.$
By~\eqref{eq:f,L} with $f=f_z, N=N_i$, we have
\begin{align}
&\langle f_z,L_{N_i}(t)\rangle\notag\\
=&
\langle f_z,L_{N_i}(0)\rangle
+M_{f_z}^{N_i}(t)+
\int_0^t
\int_{\mathbb{R}}
\frac{b_{N_i}(x)}{(z-x)^2}
\,L_{N_i}(s)(\mathrm{d}x)\,\mathrm{d} \langle \bar{B}_i  \rangle _{s}
\notag\\
&+
\frac{N_i}{2}
\int_0^t
\iint_{\mathbb{R}^2}
\frac{
(z-x)^{-2}-(z-y)^{-2}
}{x-y}
H_{N_i}(x,y)
\,L_{N_i}(s)(\mathrm{d}x)L_{N_i}(s)(\mathrm{d}y)
\,\mathrm{d}\langle\bar{B}_i\rangle_s
\notag\\
&+
\int_0^t
\int_{\mathbb{R}}
\frac{2g_{N_i}^2(x)h_{N_i}^2(x)}
{(z-x)^3}
\,L_{N_i}(s)(\mathrm{d}x)
\,\mathrm{d}\langle\bar{B}_i\rangle_s .\label{e:fL}
\end{align}
Here
$$M_{f_z}^{N_i}(t)
:=
\frac{2}{N_i}\sum_{j=1}^{N_i}
\int_0^t
\frac{
g_{N_i}\bigl(\lambda_j^{N_i}(s)\bigr)
h_{N_i}\bigl(\lambda_j^{N_i}(s)\bigr)
}{
\bigl(z-\lambda_j^{N_i}(s)\bigr)^2
}
\,\mathrm{d}\bar{B}_j(s).$$

By the $G$-Burkholder-Davis-Gundy inequality and ~\eqref{e:<MfN>}, we get
$$c\Big(
\sup_{0\leq t\leq T}|M_{f_z}^N(t)|\geq\varepsilon
\Big)
\leq
\frac{C}{\varepsilon^{2}}
\hat{\mathbb{E}}
\big[
\langle M_{f_z}^N\rangle_T
\big]
\leq \frac{4C{\overline{\sigma}}^{2} T\left\|f_z^{\prime }g_{N}h_{N}\right\|_{\infty}^{2}}{\varepsilon^{2}N}.$$
Therefore, 
$$\sum_{N=1}^{\infty}c\Big(
\sup_{0\leq t\leq T}|M_{f_z}^N(t)|\geq\varepsilon
\Big)\le \sum_{N=1}^{\infty}\frac{4C{\overline{\sigma}}^{2} T}{\varepsilon^{2}N}\left\|\frac{g_{N}h_{N}}{(z-x)^2}\right\|_{\infty}^{2}<\infty,$$
where the last step is due to the fact $\|\frac{g_{N}h_{N}}{(z-x)^2}\|_{\infty}^{2}\le C/N$ for some constant $C$ by \cite[Remark~2.3]{syy20}. Thus, by Lemma~\ref{lem:B-C}, we can verify that, for any $\varepsilon>0$, 
$$c\Big(\big[\liminf_{N\to \infty}\{\sup_{0\leq t\leq T}|M_{f_z}^N(t)|<\varepsilon\}\big]^c\Big)=0,$$
which yields the q.s.\ convergence of  $M_{f_z}^{N_{i}}$  towards zero, uniformly on compact time intervals.

Next, noting that
\begin{align*}
&\left|\int \frac{b_{N_{i}}(x)}{(z-x)^{2}} L_{N_{i}}(s)(\mathrm{d}x)-\int \frac{b(x)}{(z-x)^{2}} \mu_{s}(\mathrm{d}x)\right| \\
\leq&\left|\int \frac{b_{N_{i}}(x)-b(x)}{(z-x)^{2}} L_{N_{i}}(s)(\mathrm{d}x)\right|+\left|\int \frac{b(x)}{(z-x)^{2}}\left(L_{N_{i}}(s)(\mathrm{d}x)-\mu_{s}(\mathrm{d}x)\right)\right| \\
\leq& \frac{\sup _{x}\left|b_{N_{i}}(x)-b(x)\right|}{(\operatorname{Im}(z))^{2}}+\left|\int \frac{b(x)}{(z-x)^{2}}\left(L_{N_{i}}(s)(\mathrm{d}x)-\mu_{s}(\mathrm{d}x)\right)\right|\to 0,N_{i} \to \infty,
\end{align*}
 we have
\begin{equation*}
    \lim _{N_{i} \rightarrow \infty} \int_{0}^{t} \int \frac{b_{N_{i}}(x)}{(z-x)^{2}} L_{N_{i}}(s)(\mathrm{d}x) \mathrm{d} \langle \bar{B}_i  \rangle _{s}=\int_{0}^{t} \int \frac{b(x)}{(z-x)^{2}} \mu_{s}(\mathrm{d}x) \mathrm{d} \langle  \bar B^1  \rangle _{s}.\end{equation*}

For the last term on the right-hand side of \eqref{e:fL},  we have
\begin{align*}
&\left|\int \frac{2 g_{N_{i}}^{2}(x) h_{N_{i}}^{2}(x)}{(z-x)^{3}} L_{N_{i}}(s)(\mathrm{d}x)\right|\\
=&\frac{1}{N_{i}}\left|\int \frac{N_{i} H_{N_{i}}(x, x)}{(z-x)^{3}} L_{N_{i}}(s)(\mathrm{d}x)\right| \\
 \leq &\frac{1}{N_{i}}\left|\int \frac{N_{i} H_{N_{i}}(x, x)-H(x, x)}{(z-x)^{3}} L_{N_{i}}(s)(\mathrm{d}x)\right|+\frac{1}{N_{i}}\left|\int \frac{H(x, x)}{(z-x)^{3}} L_{N_{i}}(s)(\mathrm{d}x)\right| \\
\leq &\frac{\sup _{x, y}\left|N_{i} H_{N_{i}}(x, x)-H(x, x)\right|}{N_{i}|\operatorname{Im}(z)|^{3}}+\frac{C_{z}}{N_{i}}\left\|\frac{H(x, x)}{\left(1+|x|^{3}\right)}\right\|_{L^{\infty}(\mathrm{d}x)}\to 0,N_{i} \to \infty,
\end{align*}
where  $C_{z}$  is a constant depending only on  $z.$

Now, it remains to identify the limit of the interaction term. A direct calculation gives
\begin{align}
\frac{(z-x)^{-2}-(z-y)^{-2}}{x-y}
=&
\frac{(z-y)^2-(z-x)^2}
{(z-x)^2(z-y)^2(x-y)}
=
\frac{2z-x-y}{(z-x)^2(z-y)^2}
\notag\\
=&
\frac{1}{(z-x)(z-y)^2}
+
\frac{1}{(z-x)^2(z-y)}.
\label{eq:resolvent-difference-quotient}
\end{align}
Therefore, using the symmetry of $H_{N_i}$ and exchanging $x$ and $y$
in the second term on the right-hand side of
\eqref{eq:resolvent-difference-quotient}, we obtain
\begin{align*}
&\frac12
\iint_{\mathbb{R}^2}
\frac{(z-x)^{-2}-(z-y)^{-2}}{x-y}
N_iH_{N_i}(x,y)
\,L_{N_i}(s)(\mathrm{d}x)L_{N_i}(s)(\mathrm{d}y)
\notag\\
=&
\iint_{\mathbb{R}^2}
\frac{N_iH_{N_i}(x,y)}
{(z-x)(z-y)^2}
\,L_{N_i}(s)(\mathrm{d}x)L_{N_i}(s)(\mathrm{d}y).
\end{align*}

As $\frac{H(x,y)}{(1+|x|)(1+ y^{2})}$ is bounded and by the uniform convergence of  $N H_{N}(x, y)$  to  $H(x,y),$  we can obtain that  $\frac{N H_{N}(x, y)}{(z-x)(z-y)^{2}}$  is bounded and as $N_{i} \to \infty,$
\begin{align*}
&\left|\iint \frac{N_{i} H_{N_{i}}(x, y)}{(z-x)(z-y)^{2}} L_{N_{i}}(s)(\mathrm{d}x) L_{N_{i}}(s)(\mathrm{d}y)-\iint \frac{H(x,y)}{(z-x)(z-y)^{2}} \mu_{s}(\mathrm{d}x) \mu_{s}(\mathrm{d}y)\right| \\
\leq&\left|\iint \frac{N_{i} H_{N_{i}}(x, y)-H(x,y)}{(z-x)(z-y)^{2}} L_{N_{i}}(s)(\mathrm{d}x) L_{N_{i}}(s)(\mathrm{d}y)\right| \\&+\left|\iint \frac{H(x,y)}{(z-x)(z-y)^{2}}\left(L_{N_{i}}(s)(\mathrm{d}x) L_{N_{i}}(s)(\mathrm{d}y)-\mu_{s}(\mathrm{d}x) \mu_{s}(\mathrm{d}y)\right)\right| \\
\leq &\frac{\sup _{x, y}\left|N_{i} H_{N_{i}}(x, y)-H(x,y)\right|}{|\operatorname{Im}(z)|^{3}} +\left|\iint \frac{H(x,y)}{(z-x)(z-y)^{2}}\left(L_{N_{i}}(s)(\mathrm{d}x) L_{N_{i}}(s)(\mathrm{d}y)-\mu_{s}(\mathrm{d}x) \mu_{s}(\mathrm{d}y)\right)\right|\to 0.
\end{align*}
Then we have
 \begin{align*}
&\lim _{N_{i} \rightarrow \infty} \int_{0}^{t} \iint \frac{N_{i} H_{N_{i}}(x, y)}{(z-x)(z-y)^{2}} L_{N_{i}}(s)(\mathrm{d}x) L_{N_{i}}(s)(\mathrm{d}y) \mathrm{d} \langle \bar{B}_i  \rangle _{s} \\
=&\int_{0}^{t} \iint \frac{H(x,y)}{(z-x)(z-y)^{2}} \mu_{s}(\mathrm{d}x) \mu_{s}(\mathrm{d}y) \mathrm{d} \langle  \bar B^1  \rangle _{s}.
\end{align*}
This finishes the proof. 
\end{proof}

\subsection{Dyson Brownian motion} 
\label{subsection:Dyson Brownian motion}

\begin{proof}[Proof of Proposition~\ref{dyson converge}]
\makeatletter\def\@currentlabelname{the proof}\makeatother\label{proof:Proposition  dyson converge}
By Theorems \ref{thm: pre-compact} and~\ref{thm:limit-equation}, the sequence  $\left(L_{N}(t)\right)_{t \in[0, T]}$  is q.s.\ pre-compact in $C\left([0, T], \mathbf{P}(\mathbb{R})\right)$ and 
 any limit point  $\left(\mu_{t}, t \in[0, T]\right)$ of $L_{N} $ satisfies the equation
$$ \int f(x)  \mu_{t}(\mathrm{d}x)=\int f(x)  \mu_{0}(\mathrm{d}x)
+\frac{1}{2} \int_{0}^{t} \iint \frac{ f^{\prime}( x)-f^{\prime} ( y)}{x-y}  \mu_{s}(\mathrm{d}x)  \mu_{s}(\mathrm{d}y) \mathrm{d}\langle  \bar B^1 \rangle  _{s}.$$ 
Set  $f(x)=(z-x)^{-1}$  for some $ z \in \mathbb{C} \backslash \mathbb{R}$, and let $$V_{t}(z)= \int \frac1{z-x}  \mu_{t}(\mathrm{d}x), \text{ for } t\ge 0,$$ as given in~\eqref{eq:G}.  Then, 
\begin{align*}
    \frac{f'(x)-f'(y)}{x-y}&=\frac{ (z-x)^{-2}-(z-y)^{-2}}{x-y} =\frac{ 1}{(z-x)(z-y)^{2}}+\frac{ 1}{(z-y)(z-x)^{2}},
\end{align*}
and by the symmetry we get
\begin{align*}
 V_{t}(z)&=V_0(z)
+\frac{1}{2} \int_{0}^{t} \iint \Big[\frac{ 1}{(z-x)(z-y)^{2}}+\frac{ 1}{(z-y)(z-x)^{2}}\Big] \mu_{s}(\mathrm{d}x)  \mu_{s}(\mathrm{d}y) \mathrm{d}\langle  \bar B^1 \rangle  _{s}\\
&=V_0(z)
+\int_{0}^{t} \iint \frac{ 1}{(z-x)(z-y)^{2}} \mu_{s}(\mathrm{d}x)  \mu_{s}(\mathrm{d}y) \mathrm{d}\langle  \bar B^1 \rangle  _{s}\\
&= V_{0}(z)-\int_{0}^{t}V_{s}(z)\partial_{z} V_{s}(z)\mathrm{d}\langle  \bar B^1\rangle  _{s}.
\end{align*}
This proves \eqref{eq:Gt}.

Clearly, for each fixed $t$, $V_{t}(z)$  is an analytic function on $\mathbb{C}_{+}:=\{z\in\mathbb C:\operatorname{Im}(z)>0\}$. Now we prove the uniqueness of an analytic solution to~\eqref{eq:Gt}, which yields the uniqueness of  $\mu$. 

Fix a solution $V$ of~\eqref{eq:G}
and~\eqref{eq:Gt}. Consider the characteristic equation
\begin{equation}    \mathrm{d}z_t=V_t(z_t)\,\mathrm{d}\langle  \bar B^1 \rangle _t,
    \quad z_0\in\mathbb C_+.
    \label{eq:characteristic_equation}
\end{equation}
To justify the existence of a solution to~\eqref{eq:characteristic_equation},
we first introduce a regularization. For $\varepsilon>0$, define
$$V_t^\varepsilon(z) :=\int_{\mathbb R}\frac{\overline z-x}{|z-x|^2+\varepsilon}\,\mu_t(\mathrm{d}x),\, z\in\mathbb C_+,$$
and consider the regularized characteristic equation
\begin{equation}
    \mathrm{d}z_t^\varepsilon
    =V_t^\varepsilon(z_t^\varepsilon)
    \,\mathrm{d}\langle  \bar B^1 \rangle _t,
    \quad z_0^\varepsilon=z_0.
    \label{eq:regularized_characteristic_equation}
\end{equation}
Since $V_t^\varepsilon$ is smooth and uniformly Lipschitz on $\mathbb C_+$,
equation~\eqref{eq:regularized_characteristic_equation} admits a unique
solution.

Denote
$z_t^\varepsilon=u_t^\varepsilon+iv_t^\varepsilon$. The imaginary part of~\eqref{eq:regularized_characteristic_equation} yields
$$\frac{\mathrm{d}v_t^\varepsilon}{\mathrm{d}\langle  \bar B^1 \rangle _t}
=-v_t^\varepsilon\int_{\mathbb R}\frac{1}{|z_t^\varepsilon-x|^2+\varepsilon}\,\mu_t(\mathrm{d}x),$$
and consequently,
$$\frac{\mathrm{d}(v_t^\varepsilon)^2}{\mathrm{d}\langle  \bar B^1 \rangle _t} = -2(v_t^\varepsilon)^2\int_{\mathbb R}
\frac{1}{|z_t^\varepsilon-x|^2+\varepsilon}\,\mu_t(\mathrm{d}x)
\geq -2.$$
Noting $\langle \bar B^1 \rangle_t \leq \overline{\sigma}^{\,2}t$ q.s.,
we deduce that
\begin{equation}
    |\operatorname{Im}z_t^\varepsilon|^2
    =(v_t^\varepsilon)^2
    \geq
    (\operatorname{Im}z_0)^2-2\overline{\sigma}^{\,2}t.
    \label{eq:imaginary_part_lower_bound}
\end{equation}
In particular, for $0\leq t\leq\frac{(\operatorname{Im}z_0)^2}{4\overline{\sigma}^{\,2}}$, we have $\operatorname{Im}z_t^\varepsilon
\geq\frac{\operatorname{Im}z_0}{\sqrt{2}}>0,$  and  $$|V_t^{\varepsilon}(z_t^\varepsilon)|\le \int_{\mathbb R}\frac{1}{\operatorname{Im}z_t^\varepsilon}\,\mu_t(\mathrm{d}x)\le \frac{\sqrt 2}{\operatorname{Im}z_0},$$
which implies that $\frac{\mathrm{d}v_t^{\varepsilon}}{\mathrm{d}\langle B\rangle_t}\le \frac{\sqrt 2}{\operatorname{Im}z_0}$.  For the real part,   it follows from the definition of
$V_t^\varepsilon$ that
$$\begin{aligned}
    \left|
        \frac{\mathrm{d}u_t^\varepsilon}{\mathrm{d}\langle  \bar B^1 \rangle _t}
    \right|
    \leq
    \int_{\mathbb R}
    \frac{|u_t^\varepsilon-x|}
    {|u_t^\varepsilon-x|^2+(v_t^\varepsilon)^2+\varepsilon}
    \,\mu_t(\mathrm{d}x)\leq
    \frac{1}{2\sqrt{(v_t^\varepsilon)^2+\varepsilon}}.
\end{aligned}$$


The above estimates imply
that the family $\{z^\varepsilon\}_{\varepsilon>0}$ is equicontinuous on $[0,
    \frac{(\operatorname{Im}z_0)^2}{4\overline{\sigma}^{\,2}}]$. Hence, by the Arzel\`a--Ascoli theorem, there exists a
subsequence which converges uniformly to
a continuous process $z_t$. Passing to the limit in
\eqref{eq:regularized_characteristic_equation}, we find that $z_t$ solves
\eqref{eq:characteristic_equation}. In addition,
$\operatorname{Im}z_t\geq\frac{\operatorname{Im}z_0}{\sqrt 2}>0$ for 
$t\in[0,
    \frac{(\operatorname{Im}z_0)^2}{4\overline{\sigma}^{\,2}}]$. On the domain where $\operatorname{Im}z_t$ is bounded away from zero, the
function $V_t$ is uniformly Lipschitz. Therefore, the solution of
\eqref{eq:characteristic_equation} is unique on the above time interval.
Applying the chain rule to $V_t(z_t)$ and using
\eqref{eq:Gt} and \eqref{eq:characteristic_equation}, we obtain
\[
\begin{aligned}
    \mathrm{d}\bigl(V_t(z_t)\bigr)
    &=-V_t(z_t)\partial_zV_t(z_t)
    \,\mathrm{d}\langle  \bar B^1 \rangle _t
    +\partial_zV_t(z_t)\,\mathrm{d}z_t \\
    &=-V_t(z_t)\partial_zV_t(z_t)
    \,\mathrm{d}\langle  \bar B^1 \rangle _t
    +\partial_zV_t(z_t)V_t(z_t)
    \,\mathrm{d}\langle  \bar B^1 \rangle _t \\
    &=0.
\end{aligned}
\]
Thus, $ V_t(z_t)=V_0(z_0)$. 
Substituting this identity into~\eqref{eq:characteristic_equation} gives $\mathrm{d}z_t=V_0(z_0)\,\mathrm{d}\langle  \bar B^1 \rangle _t$.
Hence, we get 
\begin{equation}\label{e:zt}
z_t=z_0+V_0(z_0)\langle  \bar B^1 \rangle_t =:\Phi_t(z_0), 
\end{equation}
and    $V_t(z_t)=V_0(z_0)$.


Note that  $|V_0(z)|\leq\frac{1}{\operatorname{Im}z}$ and $|V'_0(z)|
\leq\frac{1}{(\operatorname{Im}z)^2}$. Then,  for $z_0\in D_t:= \big\{ z\in\mathbb C_+:\operatorname{Im}z>2\sqrt{\overline{\sigma}^{\,2}t}\big\}$,  eq.~\eqref{e:zt} yields 
$\operatorname{Im} z_t
\geq\operatorname{Im}z_0-\frac{\overline{\sigma}^2t}{\operatorname{Im}z_0}>0$
and
$|\partial_z\Phi_t(z_0)|
=
|1+V'_0(z_0)\langle  \bar B^1 \rangle _{t}|
\geq
1-\frac{\overline{\sigma}^2t}{(\operatorname{Im}z)^2}
>0$. Thus, $U_t:=\Phi_t(D_t)$
is a nonempty open subset of $\mathbb{C}_+$. Now suppose $V$ and $\widetilde V$ are two solutions with the same initial transform $V_0$.  Then 
$V_t\bigl(\Phi_t(z)\bigr)=V_0(z)=\widetilde V_t\bigl(\Phi_t(z)\bigr)$ for all $z\in D_t$. Therefore,  $V_t\bigl(w\bigr)=\widetilde V_t\bigl(w\bigr)$ for all $w\in U_t$. Since both functions are analytic on the connected domain $\mathbb C_+$, the identity theorem yields  $V_t\bigl(w\bigr)=\widetilde V_t\bigl(w\bigr)$, $w\in \mathbb{C}_+$.
This yields the uniqueness of  the solution to~\eqref{eq:G}
and~\eqref{eq:Gt}.
\end{proof}

\subsection{Wishart process}
\label{subsection:Wishart process}
\begin{proof}[Proof of Proposition~\ref{wishart converge}]
\makeatletter\def\@currentlabelname{The proof}\makeatother\label{proof:Proposition  wishart converge}
By Theorem \ref{thm: pre-compact} and Theorem \ref{thm:limit-equation}, we can obtain  that the sequence  $\left(L_{N}(t)\right)_{t \in[0, T]}$  is q.s.\ pre-compact in $C\left([0, T], \mathbf P(\mathbb{R})\right)$ and that
 any limit point  $\left(\mu_{t}, t \in[0, T]\right)$ of $L_{N} $ satisfies 
\begin{align*}
\int \frac{\mu_{t}(\mathrm{d} x)}{z-x}=&  \int \frac{\mu_{0}(\mathrm{d} x)}{z-x}+\int_{0}^{t}\left[\int \frac{c}{(z-x)^{2}} \mu_{s}(\mathrm{d} x)\right]\mathrm{d} \langle  \bar B^1 \rangle _{s} \\
&+\int_{0}^{t}\left[\iint \frac{x+y}{(z-x)(z-y)^{2}} \mu_{s}(\mathrm{d} x) \mu_{s}(\mathrm{d} y)\right] \mathrm{d}\langle  \bar B^1 \rangle _{s},
\end{align*}
for  $z \in \mathbb{C} \backslash \mathbb{R}.$ 
Denote $$V_t(z)=\int_{\mathbb R}\frac{1}{z-x}\,\mu_t(\mathrm{d}x),
    \quad z\in\mathbb C_+,$$ and we have
 $\partial_zV_t(z) =-\int_{\mathbb R}\frac{1}{(z-x)^2}\,\mu_t(\mathrm{d}x).$ 
Moreover,
\begin{align*}
    &\iint \frac{x+y}{(z-x)(z-y)^{2}} \mu_{s}(\mathrm{d} x) \mu_{s}(\mathrm{d} y)
    =\iint \frac{(x-z)+(y-z)+2z}{(z-x)(z-y)^{2}} \mu_{s}(\mathrm{d} x) \mu_{s}(\mathrm{d} y) \\
    =&\iint \frac{2z}{(z-x)(z-y)^{2}} \mu_{s}(\mathrm{d} x) \mu_{s}(\mathrm{d} y)  -\iint \frac{z-x}{(z-x)(z-y)^{2}} \mu_{s}(\mathrm{d} x) \mu_{s}(\mathrm{d} y)\\
    &-\iint \frac{z-y}{(z-x)(z-y)^{2}} \mu_{s}(\mathrm{d} x) \mu_{s}(\mathrm{d} y)\\
    =&-2 z V_{s}(z) \partial_{z} V_{s}(z)+\partial_{z} V_{s}(z)-V_{s}(z)^{2}.
\end{align*}
Consequently,  \(V_t\) satisfies
\begin{equation}
\begin{aligned}
    \mathrm{d}V_t(z)
    &=-c\,\partial_zV_t(z)\,\mathrm{d}\langle  \bar B^1 \rangle _{t}-
    \bigl(V_t(z)^2-\partial_zV_t(z)
    +2zV_t(z)\partial_zV_t(z)\bigr)
    \,\mathrm{d}\langle  \bar B^1 \rangle _{t} \\
    &=\partial_zV_t(z)\bigl(1-c-2zV_t(z)\bigr)
    \,\mathrm{d}\langle  \bar B^1 \rangle _{t}
    -V_t(z)^2\,\mathrm{d}\langle  \bar B^1 \rangle _{t}.
    \label{eq:rewritten_equation_5_19}
\end{aligned}
\end{equation}

 We now prove the
uniqueness of an analytic solution to~\eqref{eq:rewritten_equation_5_19}.  
The uniqueness of the Stieltjes transform $V_t(z)$ then implies the uniqueness of the
corresponding probability measure $\mu_t$.

Consider the characteristic equation
\begin{equation}
    \mathrm{d}z_t
    =
    \bigl(c-1+2z_tV_t(z_t)\bigr)
    \,\mathrm{d}\langle  \bar B^1 \rangle _{t},
    \quad z_0\in \mathbb C_+.
\label{eq:characteristic_equation_5_20}
\end{equation}
Denote $z=u+iv\in\mathbb C_+$. It is easy to know that
 $|V_t(z)|\leq\frac{1}{v}$ and $|\partial_zV_t(z)|\leq\frac{1}{v^2}$. 
Furthermore,
 $\operatorname{Im}V_t(z)=-\int_{\mathbb R}\frac{v}{(u-x)^2+v^2}\,\mu_t(\mathrm{d}x)\leq0$,  
and $|\operatorname{Im}V_t(z)|\leq\frac{1}{v}.$ We also have $|\operatorname{Re}V_t(z)|
    \leq\frac{1}{2v}$.  Let $y^*$ and $\rho^*$ be two positive numbers which are to be determined later and define 
\begin{equation}\label{e:D}
    D:=\Big\{  z\in\mathbb C_+:\operatorname{Im}z\geq y^*,
\ \frac{|\operatorname{Re}z|}{\operatorname{Im}z}
 \leq\rho^*\Big\}.
 \end{equation}    
Then we have for all $t\ge0$, 
\begin{equation}\label{e:onD}
|V_t(z)|\leq\frac{1}{y^*},\, |\partial_zV_t(z)|\leq\frac{1}{|y^*|^2} \text{ on  } D.
\end{equation}

We now construct a solution to \eqref{eq:characteristic_equation_5_20} by Picard iteration. Set
$z_t^0\equiv z_0=x_0+iy_0\in\mathbb C_+$,  and for $n\geq0$, define $z^{n+1}$ by
\begin{equation}
    \mathrm{d}z_t^{n+1}
    =   \bigl(c-1+2z_t^nV_t(z_t^n)\bigr)
    \,\mathrm{d}\langle  \bar B^1 \rangle _{t},
    \quad
    z_0^{n+1}=z_0.
    \label{eq:picard_iteration}
\end{equation}
Denote $ z_t^n=x_t^n+iy_t^n,\,\rho_t^n:=\frac{x_t^n}{y_t^n},
\, V_t(z_t^n)=a_t^n+ib_t^n$.
If $y_t^n>0$, we have
$b_t^n=\operatorname{Im}V_t(z_t^n)\leq0$,
and, by the estimates established above,
$|a_t^n|\leq\frac{1}{2y_t^n},\,|b_t^n|\leq\frac{1}{y_t^n}$.
Furthermore,
\[
\begin{aligned}
    z_t^nV_t(z_t^n)=
    (x_t^n+iy_t^n)(a_t^n+ib_t^n)=
    x_t^na_t^n-y_t^nb_t^n+i(x_t^nb_t^n+y_t^na_t^n).
\end{aligned}
\]
The real and imaginary parts in~\eqref{eq:picard_iteration} yield respectively, 
\begin{align}
    \mathrm{d}x_t^{n+1}
    &=
     \bigl(2x_t^na_t^n-2y_t^nb_t^n-1+c\bigr)
    \,\mathrm{d}\langle  \bar B^1 \rangle _{t},
    \label{eq:picard_real_part}\\
    \mathrm{d}y_t^{n+1}
    &=
    2\bigl(x_t^nb_t^n+y_t^na_t^n\bigr)
    \,\mathrm{d}\langle  \bar B^1 \rangle _{t}.
    \label{eq:picard_imaginary_part}
\end{align}
Assume that $y_t^n>0$ and $|\rho_t^n|\leq \rho^*$.
Then
$|x_t^na_t^n|\leq\frac{|x_t^n|}{2y_t^n}=\frac{1}{2}|\rho_t^n|
\leq\frac{\rho^*}{2}$,
and
$$0\leq -y_t^nb_t^n\leq 1,\,|x_t^nb_t^n|\leq\frac{|x_t^n|}{y_t^n}
=|\rho_t^n|,\,|y_t^na_t^n|\leq\frac12.$$
It follows from~\eqref{eq:picard_real_part} that
$$(c-\rho^*-1)
\,\mathrm{d}\langle  \bar B^1 \rangle _{t} \leq \mathrm{d}x_t^{n+1}\leq
(c+\rho^*+1)
\,\mathrm{d}\langle  \bar B^1 \rangle _{t}$$
and from~\eqref{eq:picard_imaginary_part} that
$$-(2\rho^*+1)\,\mathrm{d}\langle  \bar B^1 \rangle _{t}
\leq \mathrm{d}y_t^{n+1} \leq(2\rho^*+1)\,\mathrm{d}\langle  \bar B^1 \rangle _{t}.$$

Then we obtain that
\begin{align}
    |x_t^{n+1}|
    &\leq
    |x_0|+
    (c+\rho^*+1)\overline{\sigma}^{2}T,
    \label{eq:uniform_real_bound}\\
    y_t^{n+1}
    &\geq
    y_0-(2\rho^*+1)
    \overline{\sigma}^{2}T.
    \label{eq:uniform_imaginary_bound}
\end{align}

Given $\rho^*>0$, we can choose $y_0>0$ sufficiently large such that 
$ y^* := y_0-(2\rho^*+1)\overline{\sigma}^{\,2}T>0$ and $\frac{C_0}{y^*}\le \rho^*$, where 
$C_0:=|x_0|+c \overline{\sigma}^{2}T+(\rho^*+1)\overline{\sigma}^{\,2}T$.
Then~\eqref{eq:uniform_imaginary_bound} implies that
$ y_t^{n+1}\geq y^*$ for $t\in[0,T]$, and consequently,
 $|\rho_t^{n+1}| =\frac{|x_t^{n+1}|}{y_t^{n+1}}\leq\frac{C_0}{y^*}\le \rho^*$. Thus, we get
$ y_t^n\geq y^*, |\rho_t^n|\leq\rho^*$ for all $t\in[0,T]$ and $n\in\mathbb N$.  Moreover,  there exist constants $C_{x_0},C_{y_0}>0$, independent of
$n$, such that $|x_t^n|\leq C_{x_0},\,y^*\leq y_t^n\leq C_{y_0},\, t\in[0,T]$.
Therefore, we have $z_t^n\in D$ and 
$|z_t^n|\leq C_T:=\sqrt{C_{x_0}^2+C_{y_0}^2}$ for all $n\geq0$ and  $t\in[0,T]$.

We now prove the convergence of the Picard iteration. Define
$ F_t(w):=2wV_t(w)-1+c$.  For $w_1,w_2\in D$ with  $|w_2|\le C_T$, we have by \eqref{e:onD},
\begin{equation}
\begin{aligned}
    |F_t(w_1)-F_t(w_2)|
    &=
    2|w_1V_t(w_1)-w_2V_t(w_2)| \\
    &\leq
    2|w_1-w_2|\,|V_t(w_1)|
    +2|w_2|\,|V_t(w_1)-V_t(w_2)| \\
    &\leq
    2\Big(        \frac{1}{y^*}
        +\frac{C_T}{(y^*)^2}
    \Big)|w_1-w_2|.
\end{aligned}
\end{equation}
By~\eqref{eq:picard_iteration}, we get 
\begin{align*}
    |z_t^{n+1}-z_t^n|
    \leq
    \int_0^t
    |F_s(z_s^n)-F_s(z_s^{n-1})|
    \,\mathrm{d}\langle  \bar B^1 \rangle _s\leq
    L\int_0^t
    |z_s^n-z_s^{n-1}|
    \,\mathrm{d}\langle  \bar B^1 \rangle _s,
\end{align*}
where $L:=2\big( \frac{1}{y^*}+\frac{C_T}{(y^*)^2}\big)$. Thus, we have, denoting
$A_n(t):=\sup_{0\leq s\leq t}  |z_s^{n+1}-z_s^n|$, 
$$A_n(t)\leq L\int_0^tA_{n-1}(s)\,\mathrm{d}\langle  \bar B^1 \rangle _s\le {\overline\sigma}^2L\int_0^tA_{n-1}(s)\,\mathrm{d}s,$$
which yields
$$A_n(t) \leq \frac{(L\overline{\sigma}^{\,2}t)^n}{n!}A_0(T),
\, t\in[0,T].$$
Therefore, the sequence $\{z^n\}_{n\geq0}$ is uniformly Cauchy on $[0,T]$.
It converges uniformly to a continuous function $z_t$ satisfying
 $$z_t =z_0+\int_0^tF_s(z_s)
\,\mathrm{d}\langle  \bar B^1 \rangle _s.$$
Thus, the existence of the solution to~\eqref{eq:characteristic_equation_5_20} is established. 



 Let $z_t$ be a solution to~\eqref{eq:characteristic_equation_5_20} and $V_t(z)$ a solution to~\eqref{wishart:z-x,G}.  By the chain rule, we get
\begin{align*}
&\mathrm{d}\bigl(V_t(z_t)\bigr)\\
    &=(-c\partial_zV_t(z_t)-V_t(z_t)^2+\partial_zV_t(z_t)-2z_tV_t(z_t)\partial_zV_t(z_t))
    \,\mathrm{d}\langle  \bar B^1 \rangle _t
    +\partial_zV_t(z_t)\,\mathrm{d}z_t \\
    &=(-c\partial_zV_t(z_t)-V_t(z_t)^2+\partial_zV_t(z_t)-2z_tV_t(z_t)\partial_zV_t(z_t))
    \,\mathrm{d}\langle  \bar B^1 \rangle _t
    +\partial_zV_t(z_t)(2z_tV_t(z_t)-1+c)
    \,\mathrm{d}\langle  \bar B^1 \rangle _t \\
    &=-V_{t}(z_{t})^{2}\mathrm{d} \langle  \bar B^1 \rangle _{t}.
\end{align*}
Then,  $\frac{\mathrm{d}V_t(z_{t})}{V_t(z_{t})^2} = -\mathrm{d}\langle  \bar B^1 \rangle _{t}$, which yields
$-\frac{1}{V_t(z_{t})}+\frac{1}{V_0(z_0)}=-\langle  \bar B^1 \rangle _{t}$.
Therefore, we get 
\begin{equation}\label{e:Vz}
V_t(z_{t})=\frac{V_0(z_0)}{1+V_0(z_0)\langle  \bar B^1 \rangle_t}.
\end{equation}
Now, one can verify directly that
$$z_t = z_0 \left( 1 + V_0(z_0)\langle  \bar B^1 \rangle_t \right)^2 + (c-1) \langle  \bar B^1 \rangle_t \left( 1 + V_0(z_0)\langle  \bar B^1 \rangle_t\right)$$ 
is the  unique solution to~\eqref{eq:characteristic_equation_5_20}.

Denote
$$\Psi_t(z) := z \left( 1 + V_0(z)\langle  \bar B^1 \rangle_t \right)^2 + (c-1) \langle  \bar B^1 \rangle_t \left( 1 + V_0(z)\langle  \bar B^1 \rangle_t\right).$$
Then, we get $$\partial_z\Psi_t(z)
= \left( 1 + V_0(z)\langle  \bar B^1 \rangle_t \right)^2+\left[2z\left( 1 + V_0(z)\langle  \bar B^1 \rangle_t \right)\langle  \bar B^1 \rangle_t+(c-1)\langle  \bar B^1 \rangle_t^2\right] V_0'(z),$$ 
and for $z\in D$, by \eqref{e:onD},
\begin{align*}    |\partial_z\Psi_t(z)-1|&=\left|\big(V_0(z)\langle  \bar B^1 \rangle_t \big)^2+2V_0(z)\langle  \bar B^1 \rangle_t +\left[2z\big( 1 + V_0(z)\langle  \bar B^1 \rangle_t \big)\langle  \bar B^1 \rangle_t+(c-1)\langle  \bar B^1 \rangle_t^2\right] V_0'(z)\right|\\
    &\le\frac{\langle  \bar B^1 \rangle_t^2}{(\operatorname{Im} z)^2}+\frac{2\langle  \bar B^1 \rangle_t}{\operatorname{Im} z} +\frac{2|z|\langle  \bar B^1 \rangle_t}{(\operatorname{Im} z)^2}+\frac{2|z|\langle  \bar B^1 \rangle_t^2}{(\operatorname{Im} z)^3}+\frac{(c-1)\langle  \bar B^1 \rangle_t^2}{(\operatorname{Im} z)^2}\\
    &=\frac{c\langle  \bar B^1 \rangle_t^2}{(\operatorname{Im} z)^2}+\frac{2(|z|+\operatorname{Im} z)\langle  \bar B^1\rangle_t}{(\operatorname{Im} z)^2}+\frac{2|z|\langle  \bar B^1 \rangle_t^2}{(\operatorname{Im} z)^3}\\
    &\le\frac{c\langle  \bar B^1 \rangle_t^2}{(\operatorname{Im} z)^2}+\frac{2(\sqrt{1+(\rho^*)^2}\operatorname{Im} z+\operatorname{Im} z)\langle  \bar B^1\rangle_t}{(\operatorname{Im} z)^2}+\frac{2\sqrt{1+(\rho^*)^2}\operatorname{Im} z\langle  \bar B^1 \rangle_t^2}{(\operatorname{Im} z)^3}\\
    &=\frac{c\langle  \bar B^1 \rangle_t^2}{(\operatorname{Im} z)^2}+\frac{2(\sqrt{1+(\rho^*)^2}+1)\langle  \bar B^1 \rangle_t}{\operatorname{Im} z}+\frac{2\sqrt{1+(\rho^*)^2}\langle  \bar B^1 \rangle_t^2}{(\operatorname{Im} z)^2}\\
    &=: E(z).
    \end{align*}
    
We can choose an interior point $z^*\in D$ with $\operatorname{Im}z^*$ sufficiently large  such that  $E(z^*)<1$. Then,  $\partial_z\Psi_t(z^*)\neq 0$ and hence there is an open neighborhood
$U$ of $z^*$ satisfying $U\subset D$ such that $\Psi_t(U)\subset \mathbb C_+$ is a nonempty open set. Let $V$ and $\tilde V$ be two solutions with the same initial condition $V_0$.
 Then, by \eqref{e:Vz} we must have $V_t(w)=\tilde V_t(w)$ for all $w\in \Psi_t(U)$.
Noting that  $V_t$ and $\tilde V_t$ are analytic on $\mathbb C_+$, the identity theorem implies that
$V_t(w)=\tilde V_t(w),
\, w\in\mathbb C_+$.
This yields  the uniqueness of the solution of (\ref{wishart:z-x,G}).
\end{proof}


\subsection*{Acknowledgment} J.\ Song wishes to thank Mingshang Hu and Yongsheng Song for helpful discussions. J.\ Song is partially supported by  National Key R\&D Program of China (Grant No. 2024YFA1013501) and the Fundamental Research Funds for the Central Universities. S. Yang's work is supported by the National Natural Science Foundation of China (Grant No. 12471450) and Taishan Scholar Talent Project Youth Project.



\bibliographystyle{plain}
 \bibliography{ref}

\end{document}